\documentclass[11pt]{amsart}

\usepackage{amssymb,amsmath, enumitem}
\usepackage{amsthm}
\usepackage{graphicx}
\usepackage{mathtools}
\usepackage[usenames,dvipsnames]{xcolor}
\usepackage{multirow}
\usepackage{tikz}
\usetikzlibrary{matrix,arrows, cd}
\usepackage{breqn}
\usepackage[pdfusetitle]{hyperref}
\hypersetup{pdfauthor={David Burns, Herbert Gangl, Rob de Jeu, Alexander D. Rahm, Dan Yasaki}}
\usepackage{bookmark}

\def\spaceifletter{\futurelet\comingchar\dospaceifletter}                                                                                             
\def\dospaceifletter{\relax\ifmmode\else
  \ifcat A\noexpand\comingchar{} \fi
  \ifcat 0\noexpand\comingchar
  \ifx 0\noexpand\comingchar{} \fi
  \ifx 1\noexpand\comingchar{} \fi\ifx 2\noexpand\comingchar{} \fi
  \ifx 3\noexpand\comingchar{} \fi\ifx 4\noexpand\comingchar{} \fi
  \ifx 5\noexpand\comingchar{} \fi\ifx 6\noexpand\comingchar{} \fi
  \ifx 7\noexpand\comingchar{} \fi\ifx 8\noexpand\comingchar{} \fi
  \ifx 9\noexpand\comingchar{} \fi\fi
  \ifcat $\noexpand\comingchar{} \fi
  \ifcat \noexpand\relax\noexpand\comingchar{} \fi
\fi}

\DeclareMathSymbol{\semicolon}{6}{operators}{`;}

\let\cedille\c
\def\mycedillec{\cedille c}

\def\a{\alpha}
\def\b{\beta}
\def\c{\gamma}
\def\C{\Gamma}
\def\d{\delta}
\def\D{\Delta}

\def\l{\lambda}
\def\s{\sigma}

\def\rtmo{\sqrt{-1}}

\def\phi{\varphi}

\def\ageo{\a_{\mathrm{geo}}}
\def\Bbar #1 {\ol{\pp}(#1)}
\def\bb{\b'}
\def\balg{\b_{\mathrm{alg}}}
\def\bcyc{\b_{\mathrm{cyc}}}
\def\bgeo{\b_{\mathrm{geo}}}
\def\bound #1 {\dd_{2,#1}}

\def\coker{{\rm cok}}
\def\col#1{\begin{matrix} #1 \end{matrix}}
\def\dd{\textup{d}}
\def\ddR #1 #2 {\partial_{#1}^{#2}}

\def\et{\textup{\'et}}
\def\F{F}
\def\Ftwo{\mathbb{F}_2}
\def\Gal{\mathrm{Gal}}
\def\GL{\mathrm{GL}}
\def\id{\mathrm{id}}
\def\im{\textup{im}}
\def\ind{{\mathrm{ind}}}
\def\k #1 {\Q(\sqrt{#1})}
\def\Kind #1 {K_3(#1)^\ind}
\def\Kindtf #1 {K_3(#1)_\tf^\ind}
\def\L{\nu}
\def\LL{\mathcal{L}}
\def\Nm{\textup{Nm}}
\def\ol{\overline}

\def\O{\mathcal{O}}
\def\PGL{\mathrm{PGL}}
\def\PSL{\mathrm{PSL}}
\def\pp{\mathfrak{p}}
\def\pts #1 #2 {C_{#1}(#2)}
\def\SL{\mathrm{SL}}
\def\sustimes #1 #2 {#1 \overset\s\otimes #2}
\def\suswedge #1 {(#1 \otimes #1)_\s}
\def\tbgeo{\tilde\b_{\mathrm{geo}}}
\def\tf{{\mathrm{tf}}}
\def\tpts #1 #2 {\widetilde C_{#1}(#2^2)}

\def\tor{{\rm tor}}
\def\sign{\textup{sgn}}

\def\tL{\L'}
\def\tw{\mathrel{\tilde\wedge}}
\def\twe #1 #2 {(#1) \tw (#2)}
\def\twq #1 #2 {\twt #2/#1\tw #2}
\def\twt{\tilde\wedge^2}

\def\Fil{\textup{Fil}}
\def\CC{\mathbb{C}}
\def\P{\mathbb{P}}
\def\Q{\mathbb{Q}}
\def\Qbar{\ol{\mathbb{Q}}}
\def\R{\mathbb{R}}
\def\Z{\mathbb{Z}}
\newcommand{\ZZ}{\Z}

\def\rightiso{\buildrel{\simeq}\over{\rightarrow}}

\def\ftt{$ \partial $(flat tetrahedra)\spaceifletter}

\newcommand{\Hom}{\operatorname{Hom}}
\newcommand{\Hy}{\mathbb{H}}
\newcommand{\Spec}{{\rm{Spec}}}

\newcommand{\Tor}{{\rm{Tor}}}

\newcommand{\w}{\omega}
\newcommand{\BB}[2]{\overline{B}(#1)_{#2}}
\DeclareMathOperator{\reg}{reg}
\DeclareMathOperator{\vol}{vol}
\DeclareMathOperator{\crr}{cr}
\DeclarePairedDelimiter{\crrterm}{[}{]}
\DeclarePairedDelimiterX\set[1]{\lbrace}{\rbrace}{\def\given{\;\delimsize\vert\;}#1}

\newcommand{\RR}{\mathbb{R}}
\newcommand{\herm}{\mathcal{H}^2}

\newcommand{\QQ}{\mathbb{Q}}
\DeclareMathOperator{\Min}{Min}
\DeclareMathOperator{\Span}{span}

\newcommand{\mat}[1]{\begin{bmatrix} #1 \end{bmatrix}}
\newcommand{\PP}{\P}
\newcommand{\OO}{\O}
\def \CO{\O}

\def\tD{\mathbb{D}}
\def\Di{D}

\def\mm[#1,#2,#3,#4]{$ \begin{psmallmatrix} #1&#2\\#3&#4\end{psmallmatrix} $}

\newtheorem{theorem}[equation]{Theorem}
\newtheorem{conjecture}[equation]{Conjecture}

\newtheorem{proposition}[equation]{Proposition}
\newtheorem{prop}[equation]{Proposition}
\newtheorem{lemma}[equation]{Lemma}
\newtheorem{corollary}[equation]{Corollary}
\theoremstyle{definition}
\newtheorem{definition}[equation]{Definition}
\newtheorem{remark}[equation]{Remark}
\newtheorem{example}[equation]{Example}

\theoremstyle{plain}
\numberwithin{equation}{section}
\numberwithin{figure}{section}

\begin{document}

\title[Hyperbolic tessellations and \texorpdfstring{$K_3$}{K3}]{Hyperbolic tessellations \\
 and generators of \texorpdfstring{\emph{K}$_{\bf 3}$}{K3}\\
  for imaginary quadratic fields}

 \author{David Burns}
 \address{King's College London\\Dept. of Mathematics\\London WC2R 2LS\\United Kingdom}
 \email{david.burns@kcl.ac.uk}
 \urladdr{\url{https://nms.kcl.ac.uk/david.burns/}}

\author{Rob de Jeu}
\address{Faculteit der B\`etawetenschappen\\Afdeling Wiskunde\\Vrije Universiteit Amsterdam\\De Boelelaan 1111\\1081~HV Amsterdam\\The Netherlands}
\email{r.m.h.de.jeu@vu.nl}
\urladdr{\url{http://www.few.vu.nl/~jeu/}}

\author{Herbert Gangl}
\address{Department of Mathematical Sciences, South Road,
  Durham DH1 3LE, United Kingdom}
\email{herbert.gangl@durham.ac.uk}
\urladdr{\url{http://maths.dur.ac.uk/~dma0hg/}}

\author{Alexander~D. Rahm}
\address{Laboratoire de math\'ematiques GAATI \\ Universit\'e de la Polyn\'esie fran\mycedillec aise  \\ BP 6570 --- 98702 Faaa \\ French Polynesia}
\email{rahm@gaati.org}
\urladdr{\url{https://gaati.org/rahm/}}

\author{Dan Yasaki}
\address{Department of Mathematics and Statistics\\
University of North Carolina at Greensboro\\Greensboro, NC 27412\\ USA}
\email{d\_yasaki@uncg.edu}
\urladdr{\url{https://mathstats.uncg.edu/yasaki/}}

\begin{abstract}
We develop methods for constructing explicit generators,
modulo torsion, of the $K_3$-groups of imaginary quadratic number fields. These methods are based on either tessellations of hyperbolic
$3$-space or on direct calculations in suitable pre-Bloch groups, and lead to the very first proven examples of explicit generators, modulo torsion, of any infinite~$K_3$-group of a number field.
As part of this approach, we make several improvements to the theory of Bloch groups for $ K_3 $ of any field, predict the precise power of $2$ that should occur
in the Lichtenbaum conjecture at $ -1 $, and prove that this prediction is valid for all abelian number fields.
\end{abstract}



\maketitle

{\smaller\smaller\smaller\tableofcontents}

\section{Introduction}

\subsection{The general context} \label{general-context} Let $ F $ be a number field with ring of algebraic integers ~$ \O_F $. Then, for each
integer $m \ge 2 1$, the algebraic $K$-group $K_m(\O_F)$
of Quillen  is a fundamental invariant of~$F$, constituting a natural generalization
of the ideal class group of~$\O_F$ if $ m $ is even, and the group of units of $\O_F$ if $ m $ is odd.
By fundamental work of Quillen \cite{qui73b}, and of Borel \cite{Borel74}, this abelian group is known to be finite
if $ m $ is even, and finitely generated if $ m $ is odd.

In addition, as a natural generalization of the analytic class number formula, Lichtenbaum~\cite{lichtenbaum} has conjectured that the leading coefficient $\zeta_F^\ast(1-m)$ in the Taylor expansion at $s=1-m$ of the Dedekind
zeta function $\zeta_F(s)$ of $F$ should satisfy
\begin{equation}\label{lichtenbaum}
\zeta_F^*(1-m) = \pm 2^{n_{m,F}} \frac{|K_{2m-2}(\O_F)|}{|K_{2m-1}(\O_F)_\tor|}  R_m(F).
\end{equation}
Here we write $|X|$ for the cardinality of a finite set $X$, $K_{2m-1}(\O_F)_\tor$ for the torsion subgroup of~$K_{2m-1}(\O_F)$, $R_m(F)$ for the covolume of the image of $K_{2m-1}(\O_F)$ under the Beilinson regulator map, and $n_{m,F}$ for an undetermined integer.

Borel \cite{Borel77} proved that the quotient of $\zeta_F^*(1-m)$ by $R_m(F)$ is
rational (see Theorem~\ref{borel}), but
the identity (\ref{lichtenbaum}) has been proved unconditionally only for $F$ an abelian extension of $\QQ$ (cf. Remark~\ref{kolster remark}).
Moreover, it still is
a difficult problem to explicitly compute, except in special cases,  either~$|K_{2m-2}(\O_F)|$ or $R_m(F)$, or to give
explicit generators of $K_{2m-1}(\O_F)$ modulo torsion.

As our contribution to a solution, in this paper we
clarify the integer $n_{m,F}$, determine the precise
relation between various Bloch groups without ignoring any torsion,
develop techniques for checking divisibility in such groups,
and, in the Bloch group of each imaginary quadratic field~$F$,
algorithmically construct an element that
gives rise to a subgroup of index~$|K_{2}(\O_F)|$ in the quotient group~$K_{3}(\O_F)/K_{3}(\O_F)_{\rm tor}$.
(The number fields $ F $ of lowest degree for which~$K_{3}(\O_F)$ is infinite
are precisely the imaginary quadratic ones.)
These results are of independent interest,
but combining them allows us to give~$ |K_2(\O_F)| $, the value of $R_2(F)$, and
explicit generators of~$ K_3(\O_F) / K_3(\O_F)_\tor $, for
all imaginary quadratic fields~$ F $ of absolute discriminant at most~$1000$.
This data is available online~\cite{database} and gives
$ |K_2(\O_F)| $ for several interesting new cases.
We note it contains the very first proven examples of explicit generators
of a non-trivial $ K_3(\O_F) / K_3(\O_F)_\tor $ for any number
field $ F $, solving a problem open ever since $K_3$-groups were introduced.

Our construction of the element in the Bloch group for an imaginary
quadratic field $ F $ is based on an ideal tessellation of hyperbolic $3$-space
on which $ \GL_2(\O_F) $ acts.
After original work in Dupont-Sah \cite{DuSa82} and
Neumann-Zagier \cite{NZ}, where the gluing condition for
hyperbolic tetrahedra was formulated in an algebraic way, the first
calculations of Bloch elements from hyperbolic tessellations were, to
the best of our knowledge, implicit in
preliminary versions of~\cite{Zag86}, as well as in~\cite{Zag88} and~\cite{Ga}.
(Those tessellations give rise to elements in a Bloch group; see, e.g.,~\cite{mellit}.)

Other explicit constructions, generally based on
a triangulation of a hyperbolic manifold of finite volume and
its homology, or on group homology,
often lead to elements in a Bloch group for $ \CC $ or~$ \Qbar $ instead of
a number field (see, e.g.,~\cite{cis-mol-jon}), 
sometimes even  tensored with $ \Q $
(see, e.g., \cite[Theorem 1.2]{gonvol}, which
also discusses higher dimensional analogues).
Among the earlier results on elements of Bloch groups that are closest to our own are
\cite[Theorem~6.1]{neu-yang99},
which gives an element associated to an oriented hyperbolic 3-manifold
of finite volume,
and~\cite[Corollary~6.2.1]{Bergeron-Falbel-Guilloux},
which obtains such elements from group homology in an essentially geometric way.
By contrast, we argue directly on
the combinatorics of the tessellation under
the action of $ \GL_2(\O_F) $, resulting in a much simpler construction.
In particular, our tetrahedra are not `decorated', we
do not need any hyperbolic 3-manifold, or group homology,
and we construct our element directly for a given imaginary quadratic
field~$ F $ in such a way that it can be explicitly computed algorithmically.

We want to highlight a very interesting by-product of our methods and
calculations. For a field~$ F $, Bloch in the seminal work~\cite{bl00} constructed,
modulo some torsion, a subgroup of $ \Kind F $, the `quotient of indecomposables'
of $ K_3(F) $.
This inspired the paper~\cite{Sus90}, in which Suslin defines,
for infinite~$ F $, a `Bloch group' $ B(F) $ that describes~$ \Kind F $ modulo
some specific torsion.
Although~$ B(F) $ superficially looks very similar to the group
from~\cite{bl00}, the precise relation
between them is mysterious as Bloch's construction is based on
relative $ K $-theory whereas Suslin's is based on group homology.
Still, they are expected to be very closely related.

Instead of the construction of~\cite{bl00} we use
a variation based on a idea in \cite{blo:ltd}, as worked out in \cite{dJ95, BGdJ},
that gives a subgroup of $ \Kind F $ modulo torsion.
(We note in passing that the precise relation between this variation
and the original construction in~\cite{bl00} is not known even though
both use relative $ K $-theory.)
We map $ B(F) $ to the group from~\cite{BGdJ} in Theorem~\ref{psiprop}
but it is not a priori clear if this is compatible with the relations
of both groups with $ \Kind F $.

If~$ F $ is an imaginary quadratic field, then
from the tessellation we obtain
an element in~$ B(F) $ modulo some torsion, which under our map
gives an element in~$ \Kind F $ modulo torsion.
From the precise statement of~\eqref{lichtenbaum} for~$ F $
as obtained in Section~\ref{lich} we can then verify
for many such~$ F $ that our map induces an isomorphism between~$ B(F) $
modulo torsion and the subgroup from~\cite{BGdJ},
and that the latter is the whole of~$ \Kind F $ modulo torsion.
This provides the first concrete evidence that such statements
might hold for all fields.
For more details we refer to Section~\ref{bloch-suslin-link}.

\subsection{The main results}
We  now discuss the main contents of this article in some
more
detail.

In \S\ref{lich} we address the issue of the undetermined exponent
$ n_{m,F} $ in (\ref{lichtenbaum}) by proving that the Tamagawa number conjecture that was formulated
by Bloch and Kato in \cite{bk88} and later
extended by Fontaine and Perrin-Riou in \cite{FPR}, predicts a precise, and
more or less explicit, formula for it. For $m=2$ we can
make this conjectural formula completely explicit by using a result of Levine~\cite{lev89}.
Using results of Huber and Kings \cite{hk}, of Greither and the first author~\cite{bg} and of Flach \cite{fg} relating to the Tamagawa number conjecture,
we can then prove the (unconditional) validity of (\ref{lichtenbaum})
for all $m$ if~$F$ is abelian over $\QQ$, with a precise expression for~$n_{m,F}$.

This result is essential for our subsequent computations but is also of
independent interest. However, the arguments in \S\ref{lich} are technical in nature
and because
these methods are not used elsewhere in the article we invite any reader whose main interest is the determination of explicit generators of $K_3$-groups to
read this section up to the end of \S\ref{sor section} and then pass on to~\S\ref{ktheory}.

In \S\ref{ktheory} we shall introduce,
for any field~$ F $, a pre-Bloch group $\Bbar {F} $
based on (possibly degenerate) configurations of points.
Our approach differs slightly but crucially from that in~\cite[\S3]{gonXpam},
but this ostensibly minor improvement is essential to finding explicit generators of~$K_3$-groups
because we do not ignore any torsion.
We also analyse the corresponding variant of the second exterior power of $ F^\times$
in detail, bearing applications in computer calculations in
mind.

If $ F $ has at least four elements we shall relate $ \Bbar {F} $ to the pre-Bloch group $ \pp(F) $ of Suslin~\cite{Sus90}
(cf.~\cite[VI.5]{WeiKbook}) in a precise way, and use this to determine
the torsion subgroup of the resulting modified Bloch group~$ \BB F {} $ for a number field~$ F $.
For an imaginary quadratic field~$ F $ it turns out that $ \BB F {} $ is torsion-free,
making it much more suitable for computer calculations than~$ B(F) $.

The section actually starts with a review of some earlier results,
including one from~\cite{BGdJ} that enables us for a field $ F $ to construct a
homomorphism $\psi_F$, natural up to a universal choice of sign,
from~$ \BB F {} $ to~$\Kind F $ modulo torsion.
This is the map mentioned at the end of \S\ref{general-context},
for which it is unclear how it fits in with the relation between~$ B(F) $
and $ \Kindtf F $ of \cite{Sus90}, but this way our computations
are compatible with those in~\cite{BGdJ} and shed
light on the relation between the construction of Bloch~\cite{bl00} and
of Suslin~\cite{Sus90}.

In \S\ref{tessellations}-\ref{generators} we specialise to consider the case of an imaginary quadratic field $k$.

In this case we shall, in \S\ref{tessellations},
use the theory of perfect forms to obtain a tessellation of hyperbolic $3$-space~$\Hy^3$
on which $ \PGL_2(\O_k) $ acts, 
and from this construct an explicit well-defined
element~$\bgeo$ of the group~$\BB k {} $.
Humbert's classical formula for $\zeta_k(2)$ in terms of the volume of a fundamental domain for the action of $\PGL_2(\CO_k)$ on~$\Hy^3$ allows us to relate the image under the Beilinson regulator map of $ \psi_k(\bgeo)$ to
$ \zeta_k^\ast(-1)$. From the validity of a precise form of~\eqref{lichtenbaum}
for~$F= k $ and~$m=2$ it then follows that $\psi_k(\bgeo)$ generates a subgroup of $ \Kindtf k $ of index~$ |K_2(\O_k)|$.

The proof that $ \bgeo $ is in $ \BB k {} $ is lengthy and detailed,
since it relies on a precise study of the combinatorics of
the tessellation constructed in \S\ref{tessellations} and, for this reason, it is deferred to \S\ref{proofs}.

Then in \S\ref{generators} we use
results from previous sections to describe two concrete approaches to finding
an explicit generator of $K_3(k)^{\rm ind}_{\rm tf}$ and the order of $K_2(\CO_k)$ for
an imaginary quadratic field~$k$.

The first approach is discussed in \S\ref{first app section} and depends upon dividing $\bgeo$ by $|K_2(\OO_k)|$
directly in $ \BB k {} $ by generating elements in it using a
method involving exceptional $ S $-units (described in~\S\ref{S unit section})
and the defining relations in $ \BB k {} $.
These computations do not use the validity of~\eqref{lichtenbaum}, but they
rely on, and complement, earlier work of Belabas and the third author in \cite{BeGa}
on the orders of such $K_2(\O_k)$.
In particular, they show that the (divisional) bounds on $|K_2(\O_k)|$ obtained in loc.~cit.~ (in those cases where the order could not be precisely established) are sharp.

The second approach, discussed in \S\ref{second app section},
does rely on the known validity of~\eqref{lichtenbaum} for $ F = k  $
and $ m = 2 $ in an essential way.
Combining it with some (in practice sharp) bounds on $|K_2(\OO_k)|$ provided by \cite{BeGa},
we can draw algebraic conclusions from numerical calculations
on elements of $ \BB k {} $ obtained using exceptional $ S $-units,
which leads to the computation of a generator of~$K_3(k)^{\rm ind}_{\rm tf}$ (or of  $\overline{B}(k)$) as well as of $|K_2(\OO_k)|$ in many interesting cases.
As a concrete example, we show that $|K_2(\OO_k)| \, = 233 $ for $ k =  \k -4547 $, thereby verifying a conjecture from ~\cite{BrG}.

Some of the techniques of \S\ref{generators} can be applied to
an arbitrary number field $F$, for which essentially nothing of a
general nature beyond the result of Borel is known.
Doing this can be used to test the validity of Lichtenbaum's conjectural formula (\ref{lichtenbaum})
for $ m = 2 $, but for the sake of brevity we shall not pursue these aspects in the present article.

The article then concludes with two appendices. In Appendix \ref{orders} we shall prove several useful results about finite subgroups of $ \PGL_2(\O_k) $ of an imaginary quadratic field $k$ that are needed in earlier arguments, but for which we could not find a suitable reference. Finally, in Appendix~\ref{huge} we shall give details of the results of applying the geometrical construction of \S\ref{tessellations} and the approach described in \S\ref{first app section} to an imaginary quadratic field~$k$ for which $ |K_2(\O_k)|$ is equal to $22$.

\subsection{Notations and conventions} As a general convention we
let $ F $ denote an arbitrary field (assumed in places to be infinite)
or a number field, and $ k $ an imaginary quadratic field.

For a number field $F$ we write $\O_F$ for its ring of integers, $D_F$ for its discriminant, and~$r_1(F)$ and~$r_2(F)$ for the number of its real and complex places respectively.

For an imaginary quadratic field $k$ we set
\[\omega = \omega_k :=
\begin{cases}
  \sqrt{D_k/4} & \text{if $D_k \equiv 0 \bmod{4}$,}\\
  (1 + \sqrt{D_k})/2 & \text{if $D_k \equiv 1 \bmod{4}$}
\end{cases}
\]
so that $k = \Q(\omega)$ and $\O_k = \Z[\omega]$.

For an abelian group $M$ we write $M_{\rm tor}$ for its torsion subgroup and $M_{\tf}$ for the quotient group~$M/M_{\rm tor}$. The cardinality of a finite set $X$ will be denoted by $|X|$.

\subsection{Acknowledgements} We are grateful to the Irish Research Council for funding a work stay
of the fourth author in Amsterdam, and to the De Brun Centre for
Computational Homological Algebra for funding a work stay of the
second author in Galway.
We also thank Mathematisches Forschungsinstitut Oberwolfach for hosting the last four authors
as part of their Research in Pairs program
for this project.
We would like to thank Philippe Elbaz-Vincent for useful conversations and Nguyen Quang Do for helpful comments.
The fifth author was partially supported by NSA grant H98230-15-1-0228. This manuscript is submitted
for publication with the understanding that the United States government is authorized to produce and distribute reprints.

\section{The conjectures of Lichtenbaum and of Bloch and Kato}\label{lich}

It has long been known that the validity of (\ref{lichtenbaum}) follows from that of the conjecture originally formulated by Bloch and Kato in \cite{bk88} and then reformulated and extended by Fontaine in \cite{F} and by Fontaine and Perrin-Riou in \cite{FPR} (see Remark \ref{combi-remark}(iii) below).
However, for the main purpose of this article, it is essential
to know
not just the validity of (\ref{lichtenbaum}) but also an explicit value of the exponent $n_{m,F}$.
In this section we shall therefore derive an essentially precise formula for~$n_{m,F}$ from the assumed validity of the above conjecture of Bloch and Kato.

For each subring $\Lambda$ of $\mathbb{R}$ and each integer $a$, we write $\Lambda(a)$ for the subset $(2\pi i)^{a}\cdot \Lambda$ of $\CC$.

For a $\ZZ_p$-module $M$ we identify $M_{\tf}$ with its image in $\QQ_p\cdot M := \QQ_p\otimes_{\ZZ_p}M$, and for a homomorphism of $\ZZ_p$-modules
$\theta: M \to N$ we write $\theta_{\tf}$ for the induced homomorphism~$M_{\tf} \to~N_{\tf}$.
We write $ {\rm D}(\ZZ_2)$ for the derived category of $\ZZ_2$-modules and ${\rm D}^{\rm perf}(\ZZ_2)$ for the full triangulated subcategory of ${\rm D}(\ZZ_2)$ comprising complexes that are isomorphic (in ${\rm D}(\ZZ_2)$) to a bounded complex of finitely generated
$ \ZZ_2 $-modules.
(Note that, since the ring $ \ZZ_2 $ is regular, such complexes are precisely those that are quasi-isomorphic to a perfect complex.)

\subsection{Statement of the main result}\label{sor section}

Throughout this section, $F$ denotes a number field.

\subsubsection{}
We first review Borel's Theorem.
For this, we fix an integer $m\ge2$ and recall that Beilinson's regulator map
\[  \mathrm{reg}_m \colon K_{2m-1}(\CC) \to \mathbb{R}(m-1)\]
is compatible with the natural actions of complex conjugation on $K_{2m-1}(\CC)$ and $\mathbb{R}(m-1)$.
For each embedding $ \sigma :F \to \CC $, we consider the composite homomorphism
\[ \mathrm{reg}_{m,\sigma}: K_{2m-1}(\O_F) \overset{\sigma_*}\to K_{2m-1}(\CC) {\buildrel {\mathrm{reg}_m} \over  {\; \longrightarrow  \;}}\R(m-1) \]
where $ \sigma_* $ denotes the induced map on $ K $-groups.
We let $ \zeta_F^*(1-m) $ be the first non-zero coefficient in the Taylor expansion at $s = 1-m$ of the
Dedekind zeta function $\zeta_F(s)$ of $F$,
and set $ d_m(F) $ to be $ r_2(F) $ for even~$m$ and~$ r_1(F) + r_2(F) $ for
odd~$m$.

\begin{theorem}[\bf Borel's theorem] \label{borel} For each integer $m \ge 2$ the following hold.
\begin{itemize}
\item[(i)] The rank of $ K_{2m-2}(\O_F) $ is zero.
\item[(ii)] The rank of $ K_{2m-1}(\O_F) $ is $d_m(F)$.
\item[(iii)] Write $\mathrm{reg}_{m,F}$ for the map $K_{2m-1}(\O_F) \to \prod_{\sigma \colon F \to \CC} \R(m-1) $ given by $(\mathrm{reg}_{m,\sigma})_\sigma$. Then the image of $\mathrm{reg}_{m,F}$ is a lattice in the real vector space
    $ V_{m-1} = \{ (c_\sigma)_\sigma | c_{\overline{\sigma}} = \overline{c_\sigma} \} $,
and its kernel is $K_{2m-1}(\O_F)_{\rm tor} $.
\item[(iv)]
Let $ R_m(F) $ be the covolume of
the image of $\mathrm{reg}_{m,F}$,
with covolumes normalised so that the lattice
$ W_{m-1}^+ = \{ (c_\sigma)_\sigma | c_{\overline{\sigma}} = \overline{c_\sigma} \text{ and } c_\sigma \in \mathbb{Z}(m-1) \} $
has covolume~$ 1 $.
Then $\zeta_F(s)$ vanishes to order $d_m(F)$ at $s=1-m$, and
$\zeta_F^*(1-m) = q_{m,F} \cdot R_m(F)$ for some~$ q_{m,F} $~in~$ \Q^\times $.
\end{itemize}
\end{theorem}

\subsubsection{}For each pair of integers $i$ and $j$ with $j \in \{1,2\}$ and $i \ge j$ and each prime $p$ there exist
natural `Chern class' homomorphisms of $\ZZ_p$-modules
\begin{equation}\label{chern maps} K_{2i-j}(\O_{F})\otimes \ZZ_p \to H^j(\O_{F}[1/p],\ZZ_p(i)). \end{equation}

The first such
 homomorphism
${\rm c}_{F,i,j,p}^{\rm S}$ was constructed using higher Chern class maps by Soul\'e~\cite{soule}
(with
more details for $p=2$ being provided by Weibel \cite{cw0}), and a second~${\rm c}_{F,i,j,p}^{\rm DF}$ was constructed using \'etale $K$-theory by Dwyer-Friedlander \cite{DF}.
By~\cite[Prop.~3]{sou82}, they coincide if~$ p \ne 2 $.
For $p=2$ there is a third,
introduced independently by Kahn \cite{bk} and by Rognes and Weibel \cite{rw}, which will play a key role
for us.
All of
these maps are natural in $\O_F$ and
have finite kernels and cokernels (see Lemma~\ref{second} and
Theorem~\ref{upper bound prop} for ${\rm c}_{F,i,1,2}^{\rm S}$)
and hence induce isomorphisms of the associated $\QQ_p$-vector
spaces.

As we are mostly interested in $p=2$, we usually write ${\rm c}_{F,i,j}^{\rm S}$ and ${\rm c}_{F,i,j}^{\rm DF}$
for ${\rm c}_{F,i,j,2}^{\rm S}$ and~${\rm c}_{F,i,j,2}^{\rm DF}$.
The main result of this section, on the number $q_{m,F}$
in Theorem~\ref{borel}(iv), is then the following.
Here ${\rm det}_{\QQ_2}(\alpha)$ is the determinant of an automorphism~$ \alpha $ of a finite dimensional $\QQ_2$-vector space.

\begin{theorem}\label{main result} Fix an integer $m\ge 2$. Then the Bloch-Kato Conjecture is valid for the motive $h^0({\rm Spec}(F))(1-m)$ if and only if one has
\begin{equation}\label{zeta form} q_{m,F} = (-1)^{s_m(F)}2^{r_2(F) +t_m(F)}\cdot \frac{|K_{2m-2}(\mathcal{O}_F)|}{|K_{2m-1}(\mathcal{O}_F)_{\rm tor}|}.\end{equation}
Here we set
\[ s_m(F) := \begin{cases} [F:\Q]\frac{m}{2} - r_2(F), &\text{ if $ m $ is even,}\\
 [F:\Q]\frac{m-1}{2}, &\text{ if $m$ is odd,}\end{cases}\]
and $t_m(F) := r_1(F)\cdot t^1_m(F) + t^2_m(F)$ with
\[ t^1_m(F) := \begin{cases} -1, &\text{ if $m \equiv 1$ (mod 4)}\\
                           -2, &\text{ if $m \equiv 3$ (mod 4)}\\
                            1, &\text{ otherwise,}\end{cases}\]
and $t^2_m(F)$ the integer which satisfies
\begin{equation*} 2^{t^2_m(F)} := |{\rm cok}({\rm c}^{\rm S}_{F,m,1,{\rm tf}})|\cdot 2^{-a_m(F)} \equiv {\rm det}_{\QQ_2}((\QQ_2\cdot {\rm c}^{\rm S}_{F,m,1}) \circ (\QQ_2\cdot {\rm c}^{\rm DF}_{F,m,1})^{-1}) \,\,\,({\rm mod}\,\,\, \ZZ_2^\times)\end{equation*}
where $a_m(F)$ is $0$ except possibly when both $m \equiv 3$ (mod $4$) and $r_1(F) >0$ in which case it is an integer satisfying $0\le a_m(F) < r_1(F)$. \end{theorem}

\begin{remark}  \label{combi-remark}
(i)
The main result of Burgos Gil's book  \cite{Bur02} implies that the $m$-th Borel regulator of~$F$ is equal to $2^{d_m(F)}\cdot R_m(F)$.
So (\ref{zeta form}) leads directly to a more precise form of the conjectural formula for $\zeta^*_F(1-m)$ in terms of Borel's regulator that is given by Lichtenbaum in \cite{lichtenbaum}.

(ii)
The proof of Lemma \ref{second}(ii) below gives a closed formula for
the integer~$a_m(F)$
in Theorem~\ref{main result} (also see Remark \ref{last remark} in this regard). In addition, in \cite[Th.~4.5]{lev89} Levine shows that~${\rm c}^{\rm S}_{F,2,1}$, and hence also ${\rm c}^{\rm S}_{F,2,1,{\tf}}$, is surjective. It follows that $t^2_2(F)=0$,
so
that the exponent of $2$ in
(\ref{zeta form}) is completely explicit in the case $m=2$.
However, in order to make the kind of numerical computations
we perform in Section~\ref{second app section} for~$ m > 2 $,
 it is important to have an explicit upper bound on $t_m^2(F)$, and such a bound follows directly from Theorem \ref{upper bound prop} below.

(iii)
Up to sign and an unknown
power of $2$, Theorem \ref{main result} is proved by Huber and Kings in \cite[Th.~1.4.1]{hk} under the assumption that ${\rm c}^{\rm S}_{F,m,j,p}$ is bijective for all odd primes $p$, as
conjectured by Quillen and Lichtenbaum.
Suslin has shown the latter conjecture is
implied by the Bloch-Kato conjecture relating Milnor $K$-theory to \'etale cohomology.
Following fundamental work of Voevodsky and Rost, Weibel completed the proof of this
last conjecture \cite{cw}.
So our contribution to Theorem \ref{main result} consists
of specifying the sign (which is easy) and the exponent of~$2$.
\end{remark}

If $F$ is an abelian field, then the Bloch-Kato conjecture for $h^0({\rm Spec}(F))(1-m)$ is known to be valid: up to the $2$-primary part, this was verified indepedendently by Huber and Kings~\cite{hk} and by Greither and the first author \cite{bg},
and the $2$-primary component was subsequently resolved by Flach \cite{fg}. Theorem \ref{main result} thus leads directly to the following result.

\begin{corollary}\label{abel cor} The formula (\ref{zeta form}) is unconditionally valid if $F$ is an abelian field. \end{corollary}

\begin{example}\label{klichtenbaum} For an imaginary quadratic field $k$, Corollary \ref{abel cor} combines with Theorem \ref{borel}(iv), Remark~\ref{combi-remark}(ii) and Example  \ref{torsion comp}
to show that
$ \zeta_k^\ast(-1) = - 12^{-1} |K_2(\O_k)| R_2(k) $.
\end{example}

\begin{remark}\label{kolster remark} Independently of connections to the Bloch-Kato conjecture, the validity of (\ref{lichtenbaum}), but not (\ref{zeta form}), for abelian fields $F$ was first established  by Kolster, Nguyen Quang Do and Fleckinger in \cite{knqdf}
(the main result of loc.~cit.~contains certain erroneous Euler factors but the necessary correction is provided by Benois and Nguyen Quang Do in \cite[\S A.3]{B-NQD02}).
The general approach of \cite{knqdf} also provided motivation for the subsequent work of Huber and Kings in \cite{hk}. \end{remark}

\subsection{The proof of Theorem \ref{main result}: a first reduction}

In the sequel we abbreviate
$r_1(F)$, $r_2(F)$ and $d_m(F)$ to $r_1$, $r_2$ and $d_m$ respectively.
The functional equation of $\zeta_F(s)$ then has the form
\begin{equation}\label{dedekind fe} \zeta_F(1-s) = \frac{2^{r_2}\cdot\pi^{[F:\QQ]/2}}{|D_F|^{1/2}}\left(\frac{|D_F|}{\pi^{[F:\Q]}}\right)^s \left(\frac{\Gamma(s)}{\Gamma(1-s)}\right)^{r_2}\left(\frac{\Gamma(s/2)}{\Gamma((1-s)/2)}\right)^{r_1} \cdot \zeta_F(s),\end{equation}
with $\Gamma(s)$ the Gamma function.
Because $\zeta_F(s)$ converges at $ s = m $ we have $\zeta_F(m) > 0 $.
In addition, the function  $\Gamma(s)$ is analytic and strictly positive for $s>0$, is analytic, non-zero and of sign $(-1)^{\frac{1}{2}-n}$ at each strictly negative half-integer $n$ and has a simple pole at each strictly negative integer $n$ with residue of sign $(-1)^n$.
So from~\eqref{dedekind fe} we find that  $\zeta_F(s)$ vanishes to order~$d_m$  at $s=1-m$
(as stated in Theorem~\ref{borel}(iv)), with leading term of sign equal to~$(-1)^{[F:\Q]\frac{m}{2} - r_2}$ if $ m $ is even and to $(-1)^{[F:\Q]\frac{m-1}{2}}$ if $m$ is odd, as per the explicit formula (\ref{zeta form}).

Therefore,
in view of Remark \ref{combi-remark}(iii), in order to
prove Theorem \ref{main result} it now suffices to show that the $2$-adic component of the Bloch-Kato Conjecture for $h^0({\rm Spec}(F))(1-m)$ is valid if and only if the $2$-adic valuation of the rational number $\zeta^*_F(1-m)/R_m(F)$ is as implied by (\ref{zeta form}).
After several preliminary steps, this will be proved in \S\ref{completion}.

\subsection{The role of Chern class maps}
Regarding $F$ as fixed we set
\[ K_{m,j} := K_{2m-j}(\O_F)\otimes \ZZ_2 \,\,\,\text{ and }\,\,\, H_m^j := H^j(\O_{F}[1/2],\ZZ_2(m))\]
for each strictly positive integer $m$ and each $j \in \{1,2\}$. We also set
\[ Y_m := H^0(G_{\CC/\RR},\prod_{F \to \CC}\ZZ_2(m-1))\]
where $G_{\CC/\RR}$ acts diagonally on the product via its natural action on $\ZZ_2(m-1)$ and via post-composition on the set of embeddings $F \to \CC$.

We recall that in \cite{bk} Kahn uses the Bloch-Lichtenbaum-Friedlander-Suslin-Voevodsky spectral
sequence
to construct, for each pair of integers $i$ and $j$ with $j \in \{1,2\}$ and $i \ge j$,
a homomorphism of $\ZZ_2$-modules ${\rm c}_{i,j}^{\rm K} ={\rm c}_{F,i,j}^{\rm K}$ of the form (\ref{chern maps}) for $ p = 2 $.

We write $R\Gamma_c(\O_{F}[1/2],\Z_2(1-m))$ for the compactly supported \'etale cohomology of $\ZZ_2(1-m)$ on ${\rm Spec}(\O_{F}[1/2])$ (as defined, for example, in \cite[(3)]{bufl96}).
We recall that this complex belongs to the category ${\rm D}^{\rm perf}(\ZZ_2)$, and hence
the same holds for its (shifted) linear dual
\[ C_{m}^\bullet := R\Hom_{\Z_2}(R\Gamma_c(\O_{F}[1/2],\Z_2(1-m)),\Z_2[-2])\,.\]

In the sequel we shall write ${\rm D}(-)$ for the Grothendieck-Knudsen-Mumford determinant functor on ${\rm D}^{\rm perf}(\ZZ_2)$, as constructed in \cite{knumum}. (Note however that, since $\ZZ_2$ is local, one does not lose any significant information by (suppressing gradings and) regarding the values of ${\rm D}(-)$ as free rank one $\ZZ_2$-modules and so this is what we do.)
We shall also write ${\bigwedge}_{\ZZ_2}^aM$ for the (standard) $a$-th exterior power of  
a $ \ZZ_2 $-module $ M $ if $ a $ is a non-negative integer.

\begin{proposition}\label{av description} The map ${\rm c}_{m,1}^{\rm K}$ combines with the Artin-Verdier duality theorem
\cite[Chap.~II, Th.~3.1]{Mil86} to induce a natural identification of $\ZZ_2$-lattices
\begin{equation}\label{used later on2}
{\rm D}(C_m^\bullet) = (2^{r_1})^{t_m^1(F)}\frac{|K_{m,2}|}{|(K_{m,1})_{{\rm tor}}|}\cdot ({\bigwedge}_{\ZZ_2}^{d_m}K_{m,1,{\tf}}) \otimes_{\ZZ_2} \Hom_{\ZZ_2}({\bigwedge}_{\ZZ_2}^{d_m}Y_m,\ZZ_2)\end{equation}
where the integer $t_m^1(F)$ is as defined in Theorem \ref{main result}.
\end{proposition}

\begin{proof} We abbreviate ${\rm c}^{\rm K}_{m,j}$ to ${\rm c}_j$ and set $b_m := r_1$ if $m$ is odd and $b_m := 0$ if $m$ is even.
Then a straightforward computation of determinants shows that
\begin{alignat*}{1}
{\rm D}(C_m^\bullet) & = {\rm D}(H^0(C_m^\bullet)[0])\otimes {\rm D}(H^1(C_m^\bullet)[-1]) = {\rm D}(H^1_{m}[0])\otimes{\rm D}(H^2_{m}[-1])\otimes 2^{-b_m}{\rm D}(Y_{m}[-1])
\\
 & = ({\rm D}(\ker({\rm c}_{1})[0])^{-1}\otimes{\rm D}(K_{m,1}[0])\otimes{\rm D}({\rm cok}({\rm c}_{1})[0]))
 \\
 &\qquad\qquad \otimes ({\rm D}(\ker({\rm c}_{2})[-1])^{-1}\otimes{\rm D}(K_{m,2}[-1])\otimes{\rm D}({\rm cok}({\rm c}_{2})[-1]))\otimes 2^{-b_m}{\rm D}(Y_m[-1])
\\
  & = 2^{-b_m}\frac{|\ker({\rm c}_{1})| \cdot |{\rm cok}({\rm c}_{2})|}{|\ker({\rm c}_{2})|\cdot |{\rm cok}({\rm c}_1)|}\frac{|K_{m,2}|}{|(K_{m,1})_{{\rm tor}}|}\cdot ({\bigwedge}_{\ZZ_2}^{d_m}K_{m,1,{\tf}}) \otimes_{\ZZ_2} \Hom_{\ZZ_2}({\bigwedge}_{\ZZ_2}^{d_m}Y_m,\ZZ_2)
\,.
\end{alignat*}
Here the first equality holds as
$C_m^\bullet$ is acyclic outside degrees zero and one, the second follows from the descriptions of Lemma \ref{construction} below, the third is induced by the tautological exact sequences $0 \to \ker({\rm c}_{j}) \to K_{m,j} \xrightarrow{{\rm c}_{j}} H^j_m \to {\rm cok}({\rm c}_{j})\to 0$ for $j = 1, 2$, and the last
follows from the fact that for any finitely generated $\ZZ_2$-module $M$ and integer $a$ the $ \Z_2 $-lattice ${\rm D}(M[a])$
equals
$(|M_{\rm tor}|)^{-1}\cdot \bigwedge_{\ZZ_2}^{b}M_{\tf}$ with $b = {\rm dim}_{\QQ_2}(\QQ_2\cdot M)$
for $ a $ even and $|M_{\rm tor}|\cdot \Hom_{\ZZ_2}(\bigwedge_{\ZZ_2}^{b}M_{\tf},\ZZ_2)$
for $ a $~odd.
So it
suffices to show the product of $2^{-b_m}$ and $|\ker({\rm c}_{1})|\cdot |{\rm cok}({\rm c}_{2})|/(|\ker({\rm c}_2)|\cdot |{\rm cok}({\rm c}_{1})|)$ is~$(2^{r_1})^{t_m^1(F)}$.

This is an easy computation using that by
\cite[Th.~1]{bk}
there are integers $r_{1,4}$ and $r_{1,5}$~such that
\begin{equation}\label{Kahn isos} \begin{cases} |\ker({\rm c}_{j})|\,\,= |\coker({\rm c}_{j})| = 1, &\text{if $2m-j \equiv 0,1,2,7 \,\, ({\rm mod}\, 8)$}\\
|\ker({\rm c}_{j})|\,\,= 2^{r_1}, \,\, |\coker({\rm c}_{j})|\,\,= 1, &\text{if $2m-j \equiv 3$ \,\, ({\rm mod}\, 8)}\\
|\ker({\rm c}_{j})|\,\,= 2^{r_{1,4}}, \,\, |\coker({\rm c}_{j})|\,\,=1, &\text{if $2m-j \equiv 4$ \,\, ({\rm mod}\, 8)}\\
|\ker({\rm c}_{j})|\,\,=1, \,\, |\coker({\rm c}_{j})|\,\,= 2^{r_{1,5}}, &\text{if $2m-j \equiv 5$ \,\, ({\rm mod}\, 8)}\\
|\ker({\rm c}_{j})|\,\,=1, \,\, |\coker({\rm c}_{j})|\,\,= 2^{r_1}, &\text{if $2m-j \equiv 6$ \,\, ({\rm mod}\, 8)}\end{cases}\end{equation}
and $r_{1,4} = r_{1,5} = 0 $ if $ r_1 = 0 $,
whereas for $ r_1 > 0 $ one has $r_{1,4} \ge 0$, $r_{1,5} >0$ and~$r_{1,4} + r_{1,5} = r_1$. \end{proof}

\begin{remark}  In \cite{rw} Rognes and Weibel use a slightly different approach to Kahn to construct maps of the form ${\rm c}_{i,j}^{\rm K}$.
Their results can be used to give an alternative proof of Proposition \ref{av description}.
\end{remark}

In the sequel we write $\Sigma_\infty$, $\Sigma_\RR$ and $\Sigma_\CC$ for the sets of archimedean, real archimedean and complex archimedean places
of $F$ respectively.

\begin{lemma}\label{construction} The Artin-Verdier duality theorem induces the following identifications.
\begin{itemize}
\item[(i)] $H^0(C_{m}^\bullet) = H^1(\O_{F}[1/2],\ZZ_2(m))$.
\item[(ii)] $H^1(C_{m}^\bullet)_{\rm tor} =  H^2(\O_{F}[1/2],\ZZ_2(m)).$
\item[(iii)] $H^1(C_{m}^\bullet)_{\tf}$ is the submodule %
 \[\qquad\quad \bigoplus_{w\in \Sigma_\RR}2\cdot H^0(G_{\CC/\RR},\Z_2(m-1))\oplus \bigoplus_{w\in \Sigma_\CC}H^0(G_{\CC/\RR},\Z_2(m-1)\cdot \sigma_w \oplus \Z_2(m-1)\cdot \overline{\sigma_w})\]
 of
 \[\qquad\qquad Y_m = \bigoplus_{w\in \Sigma_\RR}H^0(G_{\CC/\RR},\Z_2(m-1))\oplus \bigoplus_{w\in \Sigma_\CC}H^0(G_{\CC/\RR},\Z_2(m-1)\cdot \sigma_w \oplus \Z_2(m-1)\cdot \overline{\sigma_w}),\]
where for each place $w \in \Sigma_\CC$ we choose a corresponding embedding $\sigma_w: F \to \CC$ and write~$\overline{\sigma_w}$ for its complex
  conjugate.
\end{itemize}
\end{lemma}

\begin{proof} For each $w$ in $\Sigma_\infty$ we write $R\Gamma_{\rm Tate}(F_w,\ZZ_2(1-m))$ for the standard
complex that computes Tate cohomology for $ F_w $ with coefficients $\ZZ_2(1-m)$,
and $R\Gamma_\Delta(F_w,\ZZ_2(1-m))$ for the mapping fibre of the natural morphism
$R\Gamma(F_w,\ZZ_2(1-m)) \to R\Gamma_{\rm Tate}(F_w,\ZZ_2(1-m))$.
We recall (see, e.g., \cite[Prop.~4.1]{bufl96}), that Artin-Verdier duality
gives an exact triangle in ${\rm D}^{\rm perf}(\ZZ_2)$ of the form %
\begin{equation}\label{exact tri} C_m^\bullet \to \bigoplus_{w\in \Sigma_\infty}R\Hom_{\Z_2}(R\Gamma_\Delta(F_w,\ZZ_2(1-m)),\ZZ_2[-1]) \to R\Gamma(\O_{F}[1/2],\ZZ_2(m))[2] \to C_m^\bullet[1]. \end{equation}
Explicit computation shows $R\Hom_{\Z_2}(R\Gamma_\Delta(F_w,\ZZ_2(1-m)),\ZZ_2[-1])$ is represented by the complex
\[  \begin{cases} \ZZ_2(m-1)[-1], &\text{if $w \in \Sigma_\CC$}\\
(\ZZ_2(m-1) \xrightarrow{\delta_m^0} \ZZ_2(m-1) \xrightarrow{\delta_m^1} \ZZ_2(m-1) \xrightarrow{\delta_m^0} \cdots )[-1], &\text{if $w\in \Sigma_\RR$,}\end{cases}\]
with $\delta_m^i $ equal to multiplication by $ 1+(-1)^{i+m}$ for $i = 0,1$.
From this and the fact that $H^2(\O_{F}[1/2],\ZZ_2(m))$ is finite,
the long exact cohomology sequence of (\ref{exact tri}) directly
implies claims~(i) and~(ii).
It also gives an exact sequence of $\ZZ_2$-modules
as in the top row of
\begin{equation*}
\begin{tikzcd}
0  \ar{r} &  H^1(C_{m}^\bullet)_\tf  \ar{r}  & \bigoplus_{w\in \Sigma_\infty} H^0(F_w,\Z_2(m-1)) \ar{r} \ar{d}{\theta} & H^3(\O_{F}[1/2],\Z_2(m)) \ar{d}{\cong}
\\
 & & \bigoplus_{w\in \Sigma_\infty} H^{3}(F_w,\Z_2(m)) \ar[equal]{r}
 & \bigoplus_{w\in \Sigma_\infty} H^3(F_w,\Z_2(m))
\end{tikzcd}
\end{equation*}
where the right hand vertical map is the isomorphism induced
by~\cite[I.4.20(b)]{Mil86}.
We define $\theta$ so that the square commutes, which
respects the direct sum decompositions of its
source and target (see the proof of \cite[Lem.~18]{bufl96}).
Since $H^{3}(F_w,\Z_2(m))$ is isomorphic to $\Z_2/2\Z_2$ for $w\in \Sigma_\RR$ and $m$ odd,
and vanishes in all
other cases, this diagram implies the description of $H^1(C_{m}^\bullet)_{\tf}$ in claim~(iii). \end{proof}

\begin{remark} The proof of Lemma \ref{construction} also shows that
\cite[Rem.~after Prop.~1.2.10]{hk} has to be modified. In
terms of the notation of loc.~cit., for the given statement to
hold one must replace~${\rm det} (T_2(r)^+)$
by $|\hat H^0(\RR,T_2(r))| \cdot \, {\rm det}(T_2(r)^+)$
instead of the asserted ${\rm det}(\hat H^0(\RR,T_2(r)))$.
\end{remark}

\subsection{Completion of the proof of Theorem \ref{main result}}\label{completion}
We write $\mathcal{E}_\RR$ for the set of embeddings $F\to \CC$ with
image in $\RR$ and set $\mathcal{E}_\CC := \Hom(F,\CC)\setminus \mathcal{E}_\RR$. For each integer $a$ we set $W_{a,\RR} := \prod_{\mathcal{E}_\RR}(2\pi i)^{a}\ZZ$ and $W_{a,\CC} := \prod_{\mathcal{E}_\CC}(2\pi i)^{a}\ZZ$ and endow the direct sum $W_a := W_{a,\RR}\oplus W_{a,\CC}$ with the diagonal action of $G_{\CC/\RR}$ that uses its natural action on $(2\pi i)^{m-1}\ZZ$ and post-composition on the embeddings~$F \to \CC$.

We write $\tau$ for the non-trivial element of $G_{\CC/\RR}$ and for each $G_{\CC/\RR}$-module $M$ we use $M^{\pm}$ to denote the submodule comprising the
elements upon which $\tau$ acts as multiplication by $\pm 1$.

Then the perfect pairing
\begin{equation*} (\QQ\otimes_\ZZ W_a) \times (\QQ\otimes_\ZZ W_{1-a}) \to (2\pi i)\QQ\end{equation*}
that sends each element $((c_\sigma),(c'_{\sigma}))$ 
to $\sum_{\sigma}c_\sigma c'_{\sigma}$ restricts to induce an identification
 \[ \Hom_\ZZ(W^+_a,(2\pi i)\ZZ) = W_{1-a,\RR}^- \oplus (W_{1-a,\CC}/(1+\tau)W_{1-a,\CC})\]
and hence also
\[ \Hom_\ZZ(W^+_a,\ZZ) = W_{-a,\RR}^+\oplus (W_{-a,\CC}/(1-\tau)W_{-a,\CC}).\]
In particular, after identifying $\QQ\otimes (W_{-a,\CC}/(1-\tau)W_{-a,\CC})$ and $\QQ\otimes W_{-a,\CC}^+$ in the natural way, one obtains an isomorphism
\begin{equation}\label{rational iso} \QQ \otimes W_{-a}^+ \cong \Hom_\QQ(\QQ\otimes W_a^+,\QQ)\end{equation}
that identifies $W_{-a}^+$ with a sublattice of $\Hom_\ZZ(W_a^+,\ZZ)$ in such a way that
\begin{equation}\label{index} |\Hom_\ZZ(W_a^+,\ZZ)/W_{-a}^+| = 2^{r_2}.\end{equation}
We then let $ \beta_m $ be defined as the composition
\begin{equation*} \beta_m: \RR\otimes K_{2m-1}(\O_F) \to \RR\otimes W_{m-1}^+ \cong \Hom_\RR(\RR\otimes W_{1-m}^+,\RR) \end{equation*}
with the first map $ \reg_{m,F} $ as in Theorem~\ref{borel}(iii)
and the second induced by (\ref{rational iso}) with $a =1-m$.
The map $\beta_m$ is bijective and we write $\beta_{m,*}$ for its induced isomorphism of $\RR$-vector spaces
\[ \RR\otimes ({\bigwedge}_{\ZZ}^{d_m} K_{2m-1}(\O_F)) \otimes \Hom_{\ZZ}({\bigwedge}_{\ZZ}^{d_m}\Hom_\ZZ(W_{1-m}^+,\ZZ),\ZZ) \to \RR\,.\]
Making explicit the formulation of \cite[Conj.~1]{bf} (which originates with Fontaine and Perrin-Riou \cite[Prop.~III.3.2.5]{FPR}) and the construction of \cite[Lem.~18]{bf}, one finds that the Bloch-Kato Conjecture for $h^0({\rm Spec}(F))(1-m)$ uses the
map $\beta_m$ rather than $ \reg_{m,F} $.
In addition, if one fixes a topological generator $\eta$ of $\ZZ_2(m-1)$,  then
mapping $\eta$ to the element of $\Hom_\QQ((2\pi i)^{1-m}\QQ,\QQ)$ that sends $(2 \pi i)^{1-m}$ to $1$
identifies $Y_m$ in Lemma~\ref{construction} with $\ZZ_2\otimes \Hom_\ZZ(W_{1-m}^+,\ZZ)$.

Given these observations, the discussion of \cite[\S1.2]{hk} shows the Bloch-Kato Conjecture asserts the
following:
if one fixes an identification of $\CC$ with $\CC_2$, then $\zeta_F^*(1-m)$ is a generator over $\ZZ_2$ of the image of ${\rm D}(C_m^\bullet)$ under the composite isomorphism
\[ \CC_2\cdot {\rm D}(C_m^\bullet) \cong \CC_2\cdot (({\bigwedge}_{\ZZ_2}^{d_m} K_{m,1}) \otimes_{\ZZ_2} \Hom_{\ZZ_2}({\bigwedge}_{\ZZ_2}^{d_m}\cdot Y_m,\ZZ_2)) \cong \CC_2.\]
Here the first map is constructed as~\eqref{used later on2} was, but
using ${\rm c}^{\rm S}_{F,m,1}$ instead of ${\rm c}^{\rm K}_{F,m,1}$
and coefficients~$ \CC_2 $ as opposed to $ \ZZ_2 $,
and the second isomorphism is
$\CC_2\otimes_\RR \beta_{m,*}$.
(Note that $ c_1^S $ has finite kernel and cokernel by
Lemma~\ref{second} below.)

Combining the above and Proposition \ref{av description} one
then finds that the Bloch-Kato Conjecture predicts
the image under $\CC_2\otimes_\RR \beta_{m,*}$ of the lattice
on the right hand side of~\eqref{used later on2}
to be equal to $\ZZ_2 \cdot \zeta_F^*(1-m)\cdot {\rm det}_{\QQ_2}((\QQ_2\cdot {\rm c}^{\rm S}_{F,m,1})\circ (\QQ_2\cdot {\rm c}^{\rm K}_{F,m,1})^{-1})^{-1}$.
Moreover, since $R_m(F)$ is defined with respect to the lattice $W_{m-1}^+$ rather than $\Hom_\ZZ(W_{1-m}^+,\Z)$, the formula (\ref{index}) implies that
\[ (\CC_2\otimes_\RR \beta_{m,*})(({\bigwedge}_{\ZZ_2}^{d_m}K_{m,1,{\tf}}) \otimes_{\ZZ_2} \Hom_{\ZZ_2}({\bigwedge}_{\ZZ_2}^{d_m}Y_m,\ZZ_2)) = 2^{r_2}\cdot R_m(F)\cdot \ZZ_2 \subset \CC_2.\]
Now to deduce (\ref{zeta form}), by direct substitution, we only
need to note that Lemma \ref{second} below implies
\[ 2^{a_m(F)}|{\rm cok}({\rm c}^{\rm S}_{F,m,1,{\tf}})|^{-1}\cdot\,{\rm det}_{\QQ_2}((\QQ_2\cdot {\rm c}^{\rm S}_{F,m,1})\circ (\QQ_2\cdot {\rm c}^{\rm K}_{F,m,1})^{-1})\in \ZZ_2^\times\]
for an integer $ a_m(F) $ as in the statement of Theorem~\ref{main result}.

The remainder of Theorem~\ref{main result} is now an immediate
consequence of Lemma~\ref{second}(ii) below.

\begin{lemma}\label{second} \
\begin{itemize}
\item[(i)] One has ${\rm det}_{\QQ_2}((\QQ_2\cdot {\rm c}^{\rm DF}_{F,m,1})\circ (\QQ_2\cdot {\rm c}^{\rm K}_{F,m,1})^{-1})\in \ZZ_2^\times$.

\item[(ii)]
The kernel and cokernel of $ {\rm c}^{\rm S}_{F,m,1} $ are finite.
If $2^{a_m(F)} = |{\rm cok}({\rm c}^{\rm K}_{F,m,1,{\tf}})|$
for $a_m(F) \ge 0$, then
\[ 2^{a_m(F)}|{\rm cok}({\rm c}^{\rm S}_{F,m,1,{\tf}})|^{-1}\cdot {\rm det}_{\QQ_2}((\QQ_2\cdot {\rm c}^{\rm S}_{F,m,1})\circ (\QQ_2\cdot {\rm c}^{\rm DF}_{F,m,1})^{-1})\in \ZZ_2^\times\,.\]
Moreover, $a_m(F) = 0$ except possibly when both $m \equiv 3$ (mod $4$) and $r_1 >0$, in which case one has $a_m(F) = r_{1,5}-1$.
\end{itemize}
\end{lemma}

\begin{proof}
Set $F' := F(\sqrt{-1})$ and $\Delta := G_{F'/F}$. With $\theta_{F'}$ denoting either ${\rm c}^{\rm S}_{F',i,1}$, $c^{\rm DF}_{F',i,1}$ or ${\rm c}^{\rm K}_{F',i,1}$ we write $\theta_F$ for the corresponding map
for~$F$. Then there is a commutative diagram
\begin{equation*}
\begin{tikzcd}
 K_{2m-1}(\O_{F'})\otimes \ZZ_2 \ar{r}{\theta_{F'}} & H^1(\O_{F'}[1/2],\ZZ_2(m))
 \\
K_{2m-1}(\O_{F})\otimes \ZZ_2 \ar{r}{\theta_F} \ar{u} & H^1(\O_{F}[1/2],\ZZ_2(m))\ar{u}
\end{tikzcd}
\end{equation*}
with the vertical maps the pullbacks.
Since $\theta_{F'}$ is natural it is $\Delta$-equivariant, hence
we may identify~$\QQ_2\cdot\theta_{F}$ with
$ H^0(\Delta,\QQ_2\cdot\theta_{F'}) \colon H^0(\Delta, K_{2m-1}(\O_{F'}) \otimes \QQ_2) \to H^0(\Delta, H^1(\O_{F'}[1/2],\QQ_2(m)))$
by using the projection formula.

By \cite[Th.~8.7 and Rem.~8.8]{DF} we have that
${\rm c}^{\rm DF}_{F',m,1}$ is surjective, and as it is a map between finitely generated $\ZZ_2$-modules of the same rank,
the induced map ${\rm c}^{\rm DF}_{F',m,1,{\tf}}$ is bijective.
This also holds for~${\rm c}_{F',m,1}^{\rm K}$ since $r_1(F') = 0$
in~\eqref{Kahn isos}, so that
$ \phi := {\rm c}^{\rm DF}_{F',m,1,{\tf}} \,\circ \, ({\rm c}^{\rm K}_{F',m,1,{\tf}})^{-1}$
is a $ \Delta $-equivariant automorphism of the $\Z_2[\Delta]$-lattice $ H^j(\O_{F'}[1/2],\ZZ_2(m))_{\tf}$.
Therefore
\begin{align*} {\rm det}_{\QQ_2}((\QQ_2\cdot {\rm c}^{\rm DF}_{F,m,1})\circ (\Q_2\cdot {\rm c}^{\rm K}_{F,m,1})^{-1})
 & = {\rm det}_{\QQ_2}(H^0(\Delta,\QQ_2\cdot {\rm c}^{\rm DF}_{F',m,1})\circ H^0(\Delta,\QQ_2\cdot {\rm c}^{\rm K}_{F',m,1})^{-1})\\
 & = {\rm det}_{\QQ_2}(H^0(\Delta,\QQ_2\cdot \varphi))\end{align*}
is in $\Z_2^\times$, proving claim~(i).

We shall prove in Theorem~\ref{upper bound prop} below that~$ {\rm c}^{\rm S}_{F,m,1,\tf} $ has finite
cokernel.
This implies $ {\rm c}^{\rm S}_{F,m,1} $ has finite kernel and
cokernel because its source and target are finitely generated $ \ZZ_2 $-modules
of the same rank.
For the remainder of the claims in~(ii), by~(i) it 
suffices to prove those with ${\rm c}^{\rm DF}_{F,m,1}$ replaced by ${\rm c}^{\rm K}_{F,m,1}$.
By (\ref{Kahn isos}) one also knows ${\rm c}^{\rm K}_{F,m,1}$, and hence also ${\rm c}^{\rm K}_{F,m,1,{\tf}}$, is surjective except possibly if $m \equiv 3$ (mod $4$) and $r_1 > 0$. In the latter case, the $\ZZ_2$-module $K_{m,1}$ is torsion-free (by
\cite[Th.~0.6]{rw}) and, since $r_1 >0$ and $m$ is odd, it is straightforward to check that $| H^1_{m,{\rm tor}}| = 2$ (see, for example, \cite[Props.~1.8 and 1.9(b)]{rw}). In this case therefore, the computation (\ref{Kahn isos}) implies that
$|{\rm cok}({\rm c}^{\rm K}_{F,m,1,{\tf}})| = 2^{-1}\cdot |{\rm cok}({\rm c}^{\rm K}_{F,m,1})| = 2^{r_{1,5}-1}$.
The remaining part of claim (ii) now follows immediately by a
computation with determinants using 
any fixed $\ZZ_2$-bases of $(K_{2m-1}(\O_{F})\otimes \ZZ_2)_{\tf}$ and
$ H^j(\O_{F}[1/2],\ZZ_2(m))_{\tf} $.
\end{proof}

\begin{remark}\label{last remark} In \cite[just after Th.~1]{bk} Kahn
asks if for $m \equiv 3$ (mod $4$) and $r_1>0$ one
has $r_{1,5} = 1$ in~\eqref{Kahn isos}
(so that $ a_m(F) = 0 $ for all $ m $).
He points out that it amounts to asking whether, in this case, the image of $H^1(\O_F[1/2],\ZZ_2(m))$ in $H^1(\O_F[1/2],\ZZ/2) \subset F^\times/(F^\times)^2$ is contained in the subgroup generated by the classes of $- 1$ and the totally
positive elements of~$F^\times$.
\end{remark}

\subsection{An upper bound for \texorpdfstring{$t_m^2(F)$}{tm2(F)}}

We now fix an integer $m\ge2$.
For a number field $E$, we let~$E' := E(\sqrt{-1})$, and write ${\rm c}_{E}$
for Soul\'e's $2$-adic Chern class map
\[ {\rm c}^{\rm S}_{E,m,1,2}: K_{2m-1}(\O_{E})\otimes \ZZ_2 \to H^1(\O_{E}[1/2],\ZZ_2(m))\,.\]

In the proof of Lemma~\ref{second} we used that~${\rm cok}({\rm c}_{F,{\tf}}) $
is finite.
Although this is commonly believed to be true, we were not able to locate a proof in the literature (cf.~Remark \ref{hw remark} below).
In addition, and as already discussed in Remark~\ref{combi-remark}(ii), for
numerical computations
one must have a computable upper bound for its order~$ 2^{t_m^2(F)} $.
For this, we prove the following result.

\begin{theorem}\label{upper bound prop} 
$|{\rm c}_{F,{\tf}}|$ is finite and divides
$  [F':F]^{d_m(F)}((m-1)!)^{d_m(F')}|K_{2m-2}(\O_{F'})| $.
\end{theorem}

As preparation for its proof, we first consider universal norm subgroups in \'etale cohomology. To do this we write $F'_\infty$ for the cyclotomic $\ZZ_2$-extension of $F'$, and
we let $F_n'$
be the unique subfield of $F'_\infty$ with $[F_n' : F'] = 2^n$ for $ n \ge 0 $, .
We also set $\Gamma := G_{F'_\infty/F'}$ and write $\Lambda$ for the Iwasawa algebra $\ZZ_2[[\Gamma]]$. For each  $\Lambda$-module $N$ and integer $a$ we write $N(a)$ for the $\Lambda$-module $N\otimes_{\ZZ_2}\ZZ_2(a)$ upon which $\Gamma$ acts diagonally.
For a finite extension $E$ of $F'$ we set $\O'_{E} := \O_{E}[1/2]$. 

We then define the `universal norm' subgroup $H_\infty^1(\O'_{F'},\ZZ_2(m))$ of $H^1(\O'_{F'},\ZZ_2(m))$ to be the image of the natural projection map $\varprojlim_nH^1(\O'_{F'_n},\ZZ_2(m)) \to H^1(\O'_{F'},\ZZ_2(m))$ where the limit is taken with respect to the natural corestriction maps.

\begin{prop}\label{univ norm prop} The index of $H_\infty^1(\O'_{F'},\ZZ_2(m))$ in $H^1(\O'_{F'},\ZZ_2(m))$ divides $|K_{2m-2}(\O_{F'})|$.
\end{prop}

\begin{proof} We write $C_\infty^\bullet$ for the object of the derived category of perfect complexes of $\Lambda$-modules that is obtained as the inverse limit of the complexes $R\Gamma(\O'_{F'_n},\ZZ_2(m))$ with respect to the natural projection morphisms
\[ R\Gamma(\O'_{F'_{n+1}},\ZZ_2(m)) \to \ZZ_2[\Gamma_{n}]\otimes_{\ZZ_2[\Gamma_{n+1}]}^{\mathbb{L}}R\Gamma(\O'_{F'_{n+1}},\ZZ_2(m)) \cong R\Gamma(\O'_{F'_{n}},\ZZ_2(m)).\]
We recall that there is a natural isomorphism $\ZZ_2\otimes^{\mathbb{L}}_{\Lambda}C^\bullet_\infty \cong R\Gamma(\O'_{F'},\ZZ_2(m))$ in
$D^{\rm perf}(\ZZ_2)$ and that this induces a natural short exact sequence of $\ZZ_2$-modules
\begin{equation}\label{descent seq} 0 \to \ZZ_2\otimes_{\ZZ_2[[\Gamma]]}H^1(C^\bullet_\infty) \xrightarrow{\pi} H^1(\O'_{F'},\ZZ_2(m)) \to H^0(\Gamma,H^2(C^\bullet_\infty)) \to 0.\end{equation}

In addition, in each degree $i$ one has $H^i(C_\infty^\bullet) = \varprojlim_nH^i(\O'_{F'_n},\ZZ_2(m))$, where the limits are taken with respect to the natural corestriction maps.
Therefore ${\rm im}(\pi) = H_\infty^1(\O'_{F'},\ZZ_2(m))$ and the $\Lambda$-module $H^2(C^\bullet_\infty)$ is isomorphic to $(\varprojlim_nH^2(\O'_{F'_n},\ZZ_2(1)))(m-1)$.

Now for each $n\ge 0$, class field theory identifies $H^2(\O'_{F'_n},\ZZ_2(1))_{\rm tor}$ with the ideal class group ${\rm Pic}(\O'_{F'_n})$ of $\O'_{F'_n}$ and $H^2(\O'_{F'_n},\ZZ_2(1))_\tf$ with a submodule of the free $\ZZ_2$-module on the set of places of $F_n'$ that are either archimedean or $2$-adic.
Hence, upon passing to the limit over $n$ and then taking $\Gamma$-invariants, we
obtain an exact sequence of $\ZZ_2$-modules
\begin{multline*} 0 \to H^0(\Gamma,X_\infty'(m-1)) \to H^0(\Gamma, H^2(C^\bullet_\infty)) \to
\bigoplus_{v\in \Sigma_2(F')\cup \Sigma_\infty(F')} H^0(\Gamma, \ZZ_2[[\Gamma/\Gamma_v]](m-1))
\end{multline*}
where $X_\infty'$ is the Galois group of the maximal unramified pro-$2$ extension of $F_\infty'$
in which all $2$-adic places split completely,  $\Sigma_2$ the set of $2$-adic places of $F'$,
and $\Gamma_v$ the decomposition subgroup of $v$ in $\Gamma$. Since $m \ge 2 $ it is also clear that each $H^0(\Gamma, \ZZ_2[[\Gamma/\Gamma_v]](m-1))$ vanishes.
Therefore $H^0(\Gamma,X_\infty'(m-1))= H^0(\Gamma, H^2(C^\bullet_\infty))$,
and by~\eqref{descent seq} it suffices to show that $H^0(\Gamma,X_\infty'(m-1))$ is finite
and of order dividing $|K_{2m-2}(\O_{F'})|$.

Next we recall that, by a standard `Herbrand Quotient' argument in Iwasawa theory (see, for example, \cite[Exer.~13.12]{wash}), if a finitely generated $\Lambda$-module $N$ is such that $H_0(\Gamma,N)$ is finite, then $H^0(\Gamma, N)$ is both finite and of order at most $|H_0(\Gamma,N)|$. In addition, since $F'$ is totally imaginary, the argument
in \cite[\S6, Lem.~1]{ps} (see also the discussion
in \cite[just before Lem.~1.2]{nqd}) shows that  $H_0(\Gamma,X_\infty'(m-1))$ is naturally isomorphic to the `\'etale wild kernel'
\[ {\rm WK}^\et_{2m-2}(F') := \ker(H^2(\O_{F'}',\ZZ_2(m)) \to \bigoplus_{w\in \Sigma_2(F')\cup \Sigma_\infty(F')}H^2(F'_w,\ZZ_2(m)))
\]
of $F'$, where the arrow denotes the natural diagonal localisation map.
Hence, to deduce the claimed result, we need only recall that, as $F'$ is totally imaginary, the group $H^2(\O_{F'}',\ZZ_2(m))$ is naturally isomorphic to the finite group $K_{2m-2}(\O_{F'})\otimes_\ZZ\ZZ_2$
by~\eqref{Kahn isos}).
\end{proof}

Turning now to the proof of Theorem~\ref{upper bound prop}, we consider for each $n$ the following diagram
\begin{equation*}
\begin{tikzcd}
 K_{2m-1}(\O'_{F'_n},\ZZ/2^n) \ar{rr}{m\cdot {\rm c}_{F'_n,2^n}} \ar{d} & & H^1(\O_{F'_n}',(\ZZ/2^{n})(m)) \ar{r}{\iota_n} \ar{d} & H^1(F'_n,(\ZZ/2^{n})(m)) \ar{d}
\\
K_{2m-1}(\O'_{F'},\ZZ/2^n) \ar{rr}{m\cdot {\rm c}_{F',2^n}} & & H^1(\O_{F'}',(\ZZ/2^{n})(m)) \ar{r}{\iota} & H^1(F',(\ZZ/2^{n})(m))
\,.
\end{tikzcd}
\end{equation*}
Here we write ${\rm c}_{E,2^n}$ for the Chern class maps $K_{2m-1}(\O'_{E},\ZZ/2^n) \to H^1(\O_{E}',(\ZZ/2^{n})(m))$ of Soul\'e, as  discussed by Weibel in \cite{cw0}, the arrows $\iota_n$ and $\iota$ are the natural inflation maps, the left hand vertical arrow is the natural transfer map, and the remaining  vertical arrows are the natural corestrictions.
The results of \cite[Prop.~2.1.1 and 4.4]{cw0}
imply that the outer rectangle of this diagram commutes.
Therefore the first square also commutes
as the maps $\iota_n$ and $\iota$ are injective.

Since the first square is compatible with change of $n$ in the natural way we may then pass to the inverse limit over $n$ to obtain a commutative diagram
\begin{equation*}
\begin{tikzcd}
\varprojlim_n K_{2m-1}(\O'_{F'_n},\ZZ/2^n) \ar{rr}{(m\cdot {\rm c}_{F'_n,2^n})_n} \ar{d} & & \varprojlim_n H^1(\O_{F'_n}',\ZZ_2(m)) \ar{d}
\\
K_{2m-1}(\O'_{F'})\otimes_\ZZ \ZZ_2 \ar{rr}{m\cdot {\rm c}_{F'}} & & H^1(\O_{F'}',\ZZ_2(m))
\,.
\end{tikzcd}
\end{equation*}
Here we use the fact that, since $K_{2m-1}(\O'_{F'})$ is finitely generated, $\varprojlim_n K_{2m-1}(\O'_{F'},\ZZ/2^n)$ identifies with $K_{2m-1}(\O'_{F'})\otimes_\ZZ \ZZ_2$ in such a way that the limit $(m\cdot {\rm c}_{F',2^n})_n$ is equal to $m\cdot {\rm c}_{F'}$.

Now the image of the right hand vertical arrow in this diagram is $H_\infty^1(\O'_{F'},\ZZ_2(m))$ and, as each $F_n'$ contains all roots of unity of order $2^n$, from \cite[Cor.~5.6]{cw0} one knows that the exponent of ${\rm cok}((m\cdot {\rm c}_{F'_n,2^n})_n)$ divides $m!$. From the commutativity of the above diagram we
then deduce $\im(m\cdot {\rm c}_{F'})$ contains $m!\cdot H_\infty^1(\O'_{F'},\ZZ_2(m))$,
so that $\im({\rm c}_{F',{\tf}})$ contains $(m-1)!\cdot H_\infty^1(\O'_{F'},\ZZ_2(m))_{\tf}$. This inclusion implies that
$ {\rm cok}({\rm c}_{F',{\tf}}) $ is finite of order dividing
\[ |(H^1(\O'_{F'},\ZZ_2(m))_{\tf}/H_\infty^1(\O'_{F'},\ZZ_2(m))_{\tf})|\cdot |(H_\infty^1(\O'_{F'},\ZZ_2(m))_{\tf}/(m-1)!\cdot H_\infty^1(\O'_{F'},\ZZ_2(m))_{\tf})|\]
and hence also $|(H^1(\O'_{F'},\ZZ_2(m))/H_\infty^1(\O'_{F'},\ZZ_2(m)))|\cdot ((m-1)!)^{d_m(F')}$.
Proposition \ref{univ norm prop} now implies that
Theorem~\ref{upper bound prop} is true if $F = F'$.

To deduce the general case of Theorem \ref{upper bound prop} we assume $F \ne F'$, write $\tau$ for the unique non-trivial element of $\Delta$ and note \cite[Prop. 4.4]{cw0} implies there is a commutative diagram
\begin{equation*}
\begin{tikzcd}
K_{2m-1}(\O'_{F'})\otimes_\ZZ\ZZ_2 \ar{rr}{{\rm c}_{F',{\tf}}} \ar{d}{T^1_\Delta } & & H^1(\O_{F'}',\ZZ_2(m))_{\tf} \ar{d}{T^2_\Delta }
\\
K_{2m-1}(\O'_{F})\otimes_\ZZ\ZZ_2 \ar{rr}{{\rm c}_{F,{\tf}}} & &  H^1(\O_{F}',\ZZ_2(m))_{\tf}\\
\end{tikzcd}
\end{equation*}
where the maps $T_\Delta^i$ are induced by the respective actions of $1+ \tau\in \ZZ_2[\Delta]$.
It follows that the index of ${\rm c}_{F,{\tf}}(\im(T^1_\Delta))$ in $\im(T^2_\Delta)$ divides $|{\rm cok}({\rm c}_{F',{\tf}})|$.
From the projection formula we see that $ \im(T^2_\Delta) $
contains~$ 2 H^1(\O_{F}',\ZZ_2(m))_{\tf} $, so that its index
in~$ H^1(\O_{F}',\ZZ_2(m))_{\tf} $ divides~$ 2^{d_m(F)} $ as
the latter is a free $ \ZZ_2 $-module of rank~$ d_m(F) $.
The statement of Theorem \ref{upper bound prop} for~$ F $ now follows.

\begin{remark}\label{hw remark} The argument of Huber and Wildeshaus in \cite[Th.~B.4.8 and Lem.~B.4.7]{hw} aims to show, amongst other things, that  ${\rm cok}({\rm c}^{\rm S}_{F,m,1,2})$ is finite. However, this argument uses in a key way results of Dwyer and Friedlander from \cite[Th.~8.7 and Rem.~8.8]{DF} that relate to ${\rm c}^{\rm DF}_{F,m,1,2}$ rather than ${\rm c}^{\rm S}_{F,m,1,2}$. To complete this argument one would thus need to investigate the relation between ${\rm cok}({\rm c}^{\rm DF}_{F,m,1,2})$ and ${\rm cok}({\rm c}^{\rm S}_{F,m,1,2})$. \end{remark}

\section{\texorpdfstring{$K$}{K}-theory, wedge complexes, and configurations of points} \label{ktheory}

Let $\F$ be an infinite field. Then it is well-known by work of Bloch and (subsequently) Suslin that $K_3(\F)$ is closely related to the Bloch group $B(\F)$ (as defined in \S\ref{connection to B} below). However, the group $B(\F)$ often contains non-trivial elements of finite order and so can be difficult for the purposes of explicit computation. With this in mind, in this section we shall introduce,
for any field~$ F $, a slight variant $\overline{B}(\F)$ of $B(\F)$ over which we have better control.

We shall also construct a natural (but unique up to a universal choice of sign)
homomorphism~$\psi_\F$ from $\overline{B}(\F)$ to $K_3(\F)_\tf^\ind$ (see Theorem~\ref{psiprop}) and are motivated to conjecture, on the basis of extensive computational evidence, that~$\psi_\F$ is bijective if $\F$ is a number field (see Conjecture~\ref{psi is iso}). We note, in particular, that these observations provide the first concrete evidence to suggest both that the groups defined by Suslin in terms of group homology and by Bloch in terms of relative $K$-theory should be related in a very natural way, and also that Bloch's group should account for all of $K_3^\ind(F)$, at least modulo torsion (cf. Remark \ref{conj rem}).

If $\F$ is imaginary quadratic, then the groups $\overline{B}(\F)$ and $K_3(\F)_\tf^\ind$ are both isomorphic to~$\ZZ$  (see Corollary \ref{Bcor} for $\overline{B}(\F)$) and we shall later use $\psi_\F$ to reduce the problem of finding a generating element of $K_3(\F)_\tf^\ind$ to computational issues in $\overline{B}(\F)$.

\subsection{Towards explicit versions of \texorpdfstring{$K_3(\F)$}{K3(F)} and \texorpdfstring{$\reg_2$}{reg2}}
\smallskip
In this section we review some earlier results that we shall need.

\subsubsection{}We first recall some basic facts concerning the \texorpdfstring{$K_3$}{K3}-group of a general field $ \F $.
For this, we write $ \Kind {\F} $ for the quotient of $ K_3(\F) $
by the image of the Milnor $K$-group~$ K_3^M(\F) $ of~$ \F $.

We recall that if $ F $ is a number field,
then the abelian group $ K_3^M(F)$ has exponent~$ 1 $ or~ $2$ and order $2^{r_1(F)}  $ (cf. \cite[p.146]{Wei04}), so that $ \Kindtf F $ identifies with $ K_3(F)_\tf$, hence is a free abelian group of rank~$ r_2(F) $ as a consequence of Theorem~\ref{borel}(ii).

We further recall that for any field $\F $ the torsion subgroup of $ \Kind {\F} $ is explicitly described by Levine in \cite[Cor.~4.6]{lev89}.

\begin{example}\label{torsion comp}
For an imaginary quadratic field $k$ in $\Qbar$, the latter result gives isomorphisms
\[ \Kind k _\tor \simeq
H^0(\Gal(\Qbar/k), \Q(2)/\Z(2) ) =  H^0(\Gal(\Qbar/\Q), \Q(2)/\Z(2) ) = \Z/24\Z \]
where the first equality is valid because complex conjugation acts trivially on~$ \Q(2)/\Z(2) $ and the second follows by explicit computation.
For any such field $k$ the abelian group $ \Kind k = K_3(k)$ is therefore isomorphic to a direct product of the form $\Z \times \Z/24\Z$.
\end{example}

\subsubsection{}\label{bloch-wigner}
 For an arbitrary field $ \F $ we set
\begin{equation*}
\twt \F^\times \coloneqq
\frac{\F^\times \otimes_\Z \F^\times}{\langle (-x) \otimes x \text{ with }  x \text{ in } \F^\times \rangle}
\,,
\end{equation*}
which is a quotient of the usual exterior power $ \F^\times \otimes_\Z \F^\times / \langle x \otimes y + y \otimes x \text{ with }  x,u \text{ in } \F^\times \rangle $.
We write $ a \tw b  $ for the class of $ a \otimes b $ in $ \twt \F^\times $,
and note that
$ a \tw b + b \tw a $ is trivial.

We next let $ \Z[\F^\flat]$  be the free abelian group on $ \F^\flat := \F \setminus \{0,1\} $,
and define the homomorphism
\begin{equation} \label{trivialnu}
 \delta_{2,\F} \colon \Z[\F^\flat] \to  \twt \F^\times
\end{equation}
by sending $[x]$ to $ (1-x) \tw x$ for each $x$ in $\F^\flat$.

We write $ \Di \colon \CC^\flat \to \R $ for the
Bloch-Wigner dilogarithm.
The value of this function at~$ z $
is defined by Bloch in \cite{bl00} by integrating 
$ \log|w|\cdot \dd {\rm arg}(1-w) - \log|1-w|\cdot \dd {\rm arg}(w) $
along any path from a point $ z_0$ in~$\R^\flat $ to $ z $.
We recall that, by differentiating, one easily shows 
the identities
\begin{equation} \label{dilogeqns}
\begin{gathered}
 \Di(z) + \Di(z^{-1}) = 0,
\qquad
 \Di(z) + \Di(1-z) = 0,
\qquad
 \Di(z) + \Di(\ol{z}) = 0,
\\
 \Di(x)-\Di(y)+\Di\Bigl(\frac{y}{x}\Bigr)-\Di\Bigl(\frac{1-y}{1-x}\Bigr)+\Di\Bigl(\frac{1-y^{-1}}{1-x^{-1}}\Bigr) = 0,
\end{gathered}
\end{equation}
for $ x $, $ y $, and $ z $ in $ \CC^\flat $ with~$ x \ne y $.

We note, in particular, that the third identity here implies that the map $ i \Di $ from $\CC^\flat$ to $\R(1)$ is equivariant with respect to the natural action of complex conjugation.

\smallskip

We quote a result connecting these notions to $ \Kindtf F $. It
underlies our construction of elements in $ \Kindtf F $
for a number field $ F $ (in particular,
in Section~\ref{tessellations}, for $ F $ imaginary quadratic).

\begin{theorem}[{\cite[Th.~4.1]{BGdJ}}] \label{regthm}
With the above notation, the following claims are valid.

\begin{itemize}
\item[(i)] There exists a homomorphism
\begin{equation*}
 \phi_\F \colon \ker(\delta_{2,\F}) \to \Kindtf {\F}
\,,
\end{equation*}
that is natural up to sign, and, after fixing a choice of sign, functorial in $ \F $.

\item[(ii)]
If $ \F $ is a number field, then the cokernel of $ \phi_F $ is finite.

\item[(iii)] There exists a universal choice of sign such that if $
F $ is any number field, and $ \sigma\colon F\to \CC$ is any embedding, then
the composition
\begin{equation*}
\begin{tikzcd}
\reg_\sigma \colon
 \ker(\delta_{2,F}) \ar{r}{\phi_F} &
 \Kindtf F = K_3(F)_\tf \ar{r}{\sigma_*} &
 K_3(\CC)_\tf \ar{r}{\reg_2} &
 \R(1)
\end{tikzcd}
\end{equation*}
is induced by sending each element $ [x] $ for $x$ in $F^\flat$ to $ i \Di(\sigma(x)) $.
\end{itemize}
\end{theorem}

\subsection{Analysing our wedge product.} In this section we obtain explicit information on the structure of $ \twt \F^\times $ for a general field $\F$.
With an eye towards implementation for the purposes of numerical calculations, we pay special attention to the case that~$ \F $ is a number field.

\subsubsection{}We first consider the abstract structure of $ \twt \F^\times $.

\begin{proposition} \label{twdescription}
For a field $ \F $,  we have a filtration $ \{0\} = \Fil_0 \subseteq \Fil_1 \subseteq \Fil_2 \subseteq \Fil_3 = \twt \F^\times $,
with~$ \Fil_1 $ the image of $ \F_\tor^\times \otimes \F_\tor^\times $ and
$ \Fil_2  $ the image of $ \F_\tor^\times \otimes \F^\times $.
Then
\begin{equation*}
 \Fil_2 =  \frac{\F_\tor^\times \otimes \F^\times }{\langle (-x) \otimes x \text{ with $ x $ in $ \F_\tor^\times $} \rangle }
\end{equation*}
and there are natural isomorphisms
\begin{alignat*}{1}
 \Fil_1/\Fil_0 & = \twt \F_\tor^\times = \frac{\F_\tor^\times \otimes_\Z \F_\tor^\times}{\langle (-x) \otimes x \text{ with } x \text{ in } \F_\tor^\times \rangle}
\\
 \Fil_2/\Fil_1 & \simeq \F_\tor^\times \otimes \F_\tf^\times
\\
 \Fil_3/\Fil_2 & \simeq \frac{\F_\tf^\times \otimes_\Z \F_\tf^\times}{\langle x \otimes x \text{ with } x \text{ in } \F_\tf^\times \rangle}
\,,
\end{alignat*}
with the last two induced by the quotient maps
$ \F_\tor^\times \otimes \F^\times \to \F_\tor^\times \otimes \F_\tf^\times $
and
$ \F^\times \otimes \F^\times \to \F_\tf^\times \otimes \F_\tf^\times $.
\end{proposition}

\begin{proof}
By taking filtered direct limits, it suffices to prove those
statements with $ \F^\times $ replaced with a finitely generated subgroup $ A $
of $ \F^\times $ that contains~$ -1 $.
We can then obtain a splitting $ A \simeq A_\tor \oplus A_\tf $
and find that the quotient for $ A $ is isomorphic to
\begin{equation*}
 \frac
   {A_\tor \otimes A_\tor \oplus A_\tor \otimes A_\tf \oplus A_\tf \otimes A_\tor \oplus A_\tf \otimes A_\tf}
   {\langle ( (-u) \otimes u , (-u) \otimes c, c \otimes u, c \otimes c ) \text{ with $ u $ in $ A_\tor $ and $ c $ in $ A_\tf $}\rangle}
\,.
\end{equation*}
Our claims follow for~$ A $ if we prove that the intersection of
\begin{equation*}
   A_\tor \otimes A_\tor \oplus A_\tor \otimes A_\tf \oplus A_\tf \otimes A_\tor \oplus 0
\end{equation*}
with the group in the denominator equals
\begin{equation*}
 \langle  ( (-u) \otimes u , u \otimes c, c \otimes u, 0 ) \text{ with $ u $ in $ A_\tor $ and $ c $ in $ A_\tf $} \rangle
\end{equation*}
as the latter is the product
\begin{equation*}
 \langle  (-u) \otimes u \text{ with $ u $ in $ A_\tor $} \rangle \times \langle ( v \otimes c , c \otimes v) \text{ with $ v $ in $ A_\tor $ and $ c $ in $ A_\tf $} \rangle \times \{0\}
 \,.
\end{equation*}
From the identity
\[  (-uc) \otimes uc - (-c) \otimes c  =  (-u) \otimes u + u \otimes c + c \otimes u \]
in $ A \times A $ it is clear that this intersection contains the given subgroup.
In order to show equality holds, 
choose a basis $ b_1,\dots,b_s $ of $ A_\tf $ and assume
that, for some integers $ m_i $, the last position in
\begin{equation} \label{intsubelt}
   \sum_i m_i ( (-u_i) \otimes u_i , (-u_i) \otimes c_i, c_i \otimes u_i, c_i \otimes c_i )
\end{equation}
is trivial.
If $ b_j $ has coefficient~$ a_{i,j} $ in $ c_i $, then $ \sum_i m_i a_{i,j}^2 = 0 $
for each~$ j $.
So each $ \sum_i m_i a_{i,j} $ is even,
\begin{equation*}
 \sum_i m_i (-1) \otimes c_i =  \sum_ j \sum_i m_i a_{i,j} (-1) \otimes b_j
\end{equation*}
is trivial, and in the second position of the element in~\eqref{intsubelt} we can replace each
$ -u_i $ with~$ u_i $.
\end{proof}

\begin{remark} \label{twotorremark}
Clearly, $ \Fil_1 $ is trivial if $ \F $ has characteristic~2.
It is also trivial
if the characteristic is not equal to~2 but $ \F^\times $ contains an element of order~$4$:
if $ u $ in $ \F^\times $ has order $ 2m $ with $ m $ even, then $ u \tw u = (-1) \tw u = m ( u \tw u ) $ in $ \twt \F^\times $,
and $ \gcd(m-1,2m) = 1 $.
Finally, if $ \F $ has characteristic not equal to~2, and $ \F^\times $
does not contain an element of order~$4$, then by decomposing $ \F_\tor^\times $
into its primary components, one sees that $ \Fil_1 $ is cyclic of order~$2$, generated by $ \twe -1 -1 $.
\end{remark}

\begin{corollary} \label{twcomputecor}
Let $ F $ be a number field, $ n $ the order of~$ F_\tor^\times $, and $ c_1,c_2,\dots $ in $ F^\times $ such that they give a
basis of $ F_\tf^\times $. Let $ m = 1 $ and $ u=1 $ if $ n $ is divisible by~$4$,
and $ m = 2 $ and $ u=-1 $ otherwise. Then the map
\begin{alignat*}{1}
\Z/m\Z \times \oplus_i \Z/n\Z \times \oplus_{i<j} \Z & \to \twt F^\times
\\
 (a, (b_i)_i, (b_{i,j})_{i,j} ) & \mapsto
  u^a \tw u
+ \sum_i  u^{b_i} \tw c_i
+ \sum_{i<j}  c_i^{b_{i,j}} \tw c_j
\end{alignat*}
is an isomorphism.
\end{corollary}

\begin{proof}
The domain has a filtration $ \Fil_l' $ for $ l=0 $, 1, 2 and 3
by taking the last $ 3-l $ positions to be trivial,
and under the homomorphism we map $ \Fil_l' $ to $ \Fil_l $ as in Proposition~\ref{twdescription}
so it induces a homomorphism $ \Fil_l'/\Fil_{l-1}' \to \Fil_l/\Fil_{l-1} $
for $ l=1 $, 2~and~3.
For $ l = 1 $ this is an isomorphism by Remark~\ref{twotorremark},
and for $ l = 2 $ and $ l = 3 $ by Proposition~\ref{twdescription}.
We now apply the five-lemma.
\end{proof}

\begin{remark} \label{twcomputerem}\
\noindent{}(i) One can get finitely many of the $ c_i $ in the corollary by taking
a basis of the free part of the $ S $-units
for a finite set $ S $ of primes of the ring of integers $ \O $
of~$ F $. If one extends~$ S $ to $ S' $, then one can add more $ c_j $ in
order to obtain a similar basis for the $ S' $-units.

\smallskip

\noindent{}(ii) For a generator $ u $ of $ F_\tor^\times $ of order~$ 2l $
and finitely many of the $ c_i $ in the corollary,
which together generate a subgroup $ A $ of $ F^\times $, the isomorphism of the corollary
becomes explicit on the image of $ \twt A $ by writing its elements
in terms of the generators, and using
$ c_i \tw c_j + c_j \tw c_i = 0  $ if $ i \ne j $,
$ c_i \tw c_i = (-1) \tw c_i $,
$ u \tw c_i + c_i \tw u = 0  $,
as well as that $ u \tw u $ equals $ \twe -1 -1 $ for~$ l $ odd and
is trivial for~$ l $ even.
\end{remark}

\subsubsection{}Any field extension $\F \to \F'$ induces a homomorphism from $ \twt \F^\times $ to $ \twt (\F')^\times $.
We determine its kernel for $\F = \QQ$ and $\F'$  imaginary quadratic.
This will be important for Theorem~\ref{mthm}.

\begin{lemma} \label{kertwmap}
Let $ d $ be a positive square-free integer, and let $ k = \k -d $.
If $ d \ne 1 $ then the kernel of the map $ \twt {\Q^\times} \to \twt k^\times $ has
order~$2$, with $ \twe -1 -d $ as non-trivial element.
If $ d = 1 $ then the kernel is non-cyclic of order~$4$, and is generated by
$ \twe -1 -1 $ and $ (-1) \tw 2 $.
\end{lemma}

\begin{proof}
The given elements are clearly in the kernel
(for $ \k -1 $ use $ 2 = \rtmo (1-\rtmo)^2 $
and that $ \rtmo \tw \rtmo $ is trivial by Remark~\ref{twotorremark}).
In addition, by Proposition~\ref{twdescription}
(or Corollary~\ref{twcomputecor} with the prime numbers as $ c_i $)
the elements generate a subgroup of the stated order.
It is therefore enough to check the size of the kernel.

Clearly, from the description in Proposition~\ref{twdescription}
the kernel is contained in $ \Fil_2 $ on $ \twt {\Q} $.
Because the map on the $ \Fil_1$-pieces is surjective by Remark~\ref{twotorremark},
with kernel of order~$2$ if $ d = 1 $ and trivial otherwise, we
only have to show that the kernel for $ \Fil_2/\Fil_1 $ has order~$2$.

For this we use the description of $ \Fil_2/\Fil_1 $ in Proposition~\ref{twdescription}. If $ |k_\tor^\times| = 2m $
then the kernel corresponds to the kernel of the map $ \Q_\tf^\times/2 \to k_\tf^\times/2m $
given by raising to the $ m $th power, as $ -1 $ is the $ m $th
power of a generator of $ k_\tor^\times $. We solve $ a^m = u \a^{2m} $
with $ a $ in $ \Q^\times $, $ u $ in $ k_\tor^\times $, and $ \a $ in~$ k^\times $, or, equivalently, $ a = v \a^2 $ for some $ v $ in
$ k_\tor^\times $ as $ u $ is an $ m $th power in~$ k^\times $.

After some calculation, we find that, for $ d \ne 1 $, $ a $ is of the form $ \pm b^2 $ or $ \pm d\cdot b^2 $
for some $ b $ in~$ \Q^\times $, and for $ d = 1 $ that
it is of the form $ \pm b^2 $ or $ \pm 2 b^2 $.
In either case, this leads to two elements in $ \Q_\tf^\times/2 $ that
are in the kernel, as required.
\end{proof}

\subsection{Configurations of points, and a modified Bloch group.}

In order to be able to apply Theorem~\ref{regthm} in our
geometric construction
of elements in the indecomposable $K_3$-groups of imaginary quadratic fields
in \S\ref{tessellations}, it is convenient to make technical modifications of well-known constructions of Suslin \cite{Sus90} and of Goncharov~\cite[p.~73]{gonXpam}.
As a result, we shall be able to be more precise about torsion
in the resulting Bloch groups and some of the homomorphisms involved.

However, in order to be able to take finite non-trivial torsion in stabilisers
of points in $ \P_\F^1 $ into account, and to be able to work with
groups like $ \PGL_2(\F) $ instead of $ \GL_2(\F) $ whenever necessary, we are forced to work in somewhat greater generality.

\subsubsection{}\label{cr3 sec def} Let $ \F $ be a field, and fix two subgroups $ \L \subseteq \tL $ of $ \F^\times $.
(Typically, we have in mind $ \L = \{1\} $ or $ \{\pm1\} $,
and $ \tL $ the torsion subgroup of the units of the ring of algebraic
integers in a number field.)
Let $ \D = \GL_2(\F)/\L $.
Let $ \LL $ be the set of orbits for the action of $ \tL $
on $ \F^2 \setminus\{(0,0)\} $ given by scalar multiplication,
which has a natural map to $ \P_\F^1 $.
The extreme cases $ \tL = \F^\times $ and $ \tL = \set{1} $ give $ \LL = \P_\F^1 $
and $ \LL = \F^2\setminus\set{(0,0)} $ respectively.
For $ n \ge 0 $ we let $ \pts n {\LL} $ be the free abelian group
with as generators $ (n+1) $-tuples $ (l_0,\dots,l_n) $ of elements
in $ \LL $ such that if $ l_{i_1} $ and $ l_{i_2} $ have the same image
in $ \P_\F^1 $ then $ l_{i_1} = l_{i_2} $
(see Remark~\ref{degenrem} for an explanation of this condition).
We shall call such a tuple $ (l_0,\dots,l_n) $ with all $ l_i $ distinct in $ \LL $
(or equivalently, in $ \P_\F^1 $) non-degenerate, and we shall
call it degenerate otherwise.
Then $ \D $ acts on $ \pts n {\LL} $ as $ \L \subseteq \tL $, and
with the usual boundary map
$ \dd \colon \pts n {\LL} \to \pts n-1 {\LL} $ for
$ n\ge 1 $ given by
\[  \dd (l_0,\dots,l_n) = \sum_{i=0}^n (-1)^i (l_0,\dots,\widehat{l_i},\dots,l_n), \]
where $\widehat{l_i}$ indicates that the term $l_i$ is omitted, we get a complex
\begin{equation} \label{eq:configuration-complex}
\begin{tikzcd}
 \cdots \ar{r}{\dd} &  \pts 4 {\LL} \ar{r}{\dd} & \pts 3 {\LL}
 \ar{r}{\dd} & \pts 2 {\LL} \ar{r}{\dd} & \pts 1 {\LL}
 \ar{r}{\dd} &  \pts 0 {\LL}
  \end{tikzcd}
\end{equation}
of $ \Z[\D] $-modules.

For three non-zero points $ p_0 $, $ p_1 $ and $ p_2 $ in $ \F^2 $
with distinct images in $ \P_\F^1 $,
we define $ \crr_2(p_0,p_1,p_2) $ in $ \twt \F^\times $
by the rules:
\begin{itemize}
\item[$\bullet$]
$ \crr_2(gp_0,gp_1,gp_2) =  \crr_2(p_0,p_1,p_2) $ for every $ g $ in $ \GL_2(\F) $;

\item[$\bullet$]
$ \crr_2((1,0),(0,1),(a,b)) = a \tw b $.
\footnote{Goncharov in \cite[\S3]{gonXpam} maps this to $ \twe -1 -1 + b \tw a $.}
\end{itemize}
From
$ \crr_2((0,1), (1,0), (a,b) ) = b \tw a $
and
$ \crr_2((1,0), (a,b), (0,1) ) = (-ab^{-1}) \tw b^{-1} = b \tw a $
we see that $ \crr_2 $ is alternating. It is also clear that if
we scale one of the $ p_i $ by $ \l $ in~$ \tL $ then~$ \crr_2(p_0,p_1,p_2) $
changes by a term $ {\l} \tw c $ with $ c $ in~$ \F^\times $.  Let
\begin{equation} \label{messdef}
 \twq {\tL} {\F^\times} 
=
 \frac{\twt \F^\times}{\langle\twe {\l} c \text{ with }\l\text{ in }\tL\text{ and }c\text{ in }\F^\times\rangle}
\,.
\end{equation}
We then define a homomorphism
\[ f_{2, \F} \colon \pts 2 {\LL} \to \twq {\tL} {\F^\times} \]
by letting it be trivial on a degenerate generator $ (l_0,l_1,l_2) $,
and by mapping a non-degenerate generator $ (l_0,l_1,l_2) $
to $ \crr_2(p_0,p_1,p_2) $ with $ p_i $ a point in $ l_i $.
(We suppress $ \tL $ from the notation.)

We next define a homomorphism
\[ f_{3,\F} \colon \pts 3 {\LL} \to \Z[\F^\flat]\]
as follows. On a degenerate generator $ (l_0,l_1,l_2,l_3) $ we let $ f_{3,\F} $
be trivial, and we let it map a non-degenerate generator $ (l_0,l_1,l_2,l_3) $
to $ [\crr_3(\ol{l_0},\ol{l_1},\ol{l_2},\ol{l_3})] $, the generator
for the cross-ratio $ \crr_3 $ of the images of the points in $ \P_\F^1 $.
Recall that $ \crr_3 $ is defined by rules similar to those for~$ \crr_2 $:
\begin{enumerate}
\item[$\bullet$]
$ \crr_3(g\ol{l_0},g\ol{l_1},g\ol{l_2},g\ol{l_3}) = \crr_3(\ol{l_0},\ol{l_1},\ol{l_2},\ol{l_3}) $
for every $ g $ in $ \GL_2(\F) $;
\item[$\bullet$]
$ \crr_3([1,0],[0,1],[1,1],[x,1]) = x $ for $ x $ in~$ \F^\flat $.
\end{enumerate}

\begin{remark} \label{classicalcr}
From the $ \GL_2(\F) $-equivariance of $ \crr_3 $ one sees by a direct calculation
that, for~$ l_0,l_1,l_2,l_3 $ different non-zero
points in $ \F^2 $,
\[ \crr_3(\ol{l_0}, \ol{l_1}, \ol{l_2}, \ol{l_3})  = \frac{\det(\mat{l_1~l_3})  \det(\mat{l_2~ l_4})}{\det(\mat{l_1~ l_4}) \det(\mat{l_2~ l_3})}.\]
As is well known, from this, or by a direct calculation, we see that permuting
the four points can give the following related possibilities
for a cross-ratio:
$ x $, $ 1-x^{-1} $, $ (1-x)^{-1} $
for even permutations, and
$ 1-x $, $ x^{-1} $, $ (1-x^{-1})^{-1} $
for odd ones, with the subgroup $ V_4 $ of $ S_4 $ acting~trivially.
\end{remark}

\subsubsection{}In the next result we consider the homomorphism
\[  \d_{2,\F}^{\tL} \colon \Z[\F^\flat] \to \twq {\tL} {\F^\times} \]
that sends each element $[x]$ for $x$ in $\F^\flat$ to the class
of $(1-x) \tw x$.  If $ \tL $ is trivial then this is still the
map~$ \d_{2,\F} $ of~\eqref{trivialnu}.

\begin{lemma}\label{commutative lemma} The following diagram commutes.
\begin{equation*}
\begin{tikzcd}
  \pts 3 {\LL} \ar{r}{\dd} \ar{d}{f_{3,\F}} &   \pts 2 {\LL} \ar{d}{f_{2,\F}}
\\
 \Z[\F^\flat] \ar{r}{\d_{2,\F}^{\tL}} & \twq {\tL} {\F^\times} 
\end{tikzcd}
\end{equation*}
\end{lemma}

\begin{proof} It suffices to check commutativity for each generating element $ (l_0,l_1,l_2,l_3) $ of $\pts 3 {\LL} $.

For $ (l_0,l_1,l_2,l_3) $ non-degenerate this follows by an explicit computation. Specifically, by using the $ \GL_2(\F) $-invariance of both $ f_{3,\F} $ and $ f_{2,\F} $
and the $ \GL_2(\F) $-equivariance of $ \dd $ one can assume that $ l_0,l_1,l_2,l_3 $ are
the classes of $ (a,0) $, $ (0,b) $, $ (1,1) $ and $ (xc,c) $
in $ \LL $
for some and $ a $, $ b $ and~$ c $ in $ \F^\times $ and $ x $ in $ \F^\flat $,
which results in $ [x] $ in $ \Z[\F^\flat] $ under $ f_{3,\F} $
and the class of $ (1-x) \tw x $ under~$ f_{2,\F} \circ \dd $.

For a degenerate tuple $ (l_0,l_1,l_2,l_3) $ the commutativity
is obvious if $ \{l_0,l_1,l_2,l_3\} $
has at most two elements as then $ f_{2,\F} $ is trivial on every
term in $ \dd (l_0,l_1,l_2,l_3) $.

If $ \{l_0,l_1,l_2,l_3\} $ consists of
$ A $, $ B $ and $ C $ with $ A $ occurring twice among
$ l_0 $, $ l_1 $, $ l_2 $ and $ l_3 $, then up to permuting $ B $ and
$ C $ the
possibilities for $ (l_0,l_1,l_2,l_3) $ are
$ (A,A,B,C) $, $ (A,B,A,C) $, $ (A,B,C,A) $, $ (B,A,A,C) $,
$ (B,A,C,A) $ and $ (B,C,A,A) $.
After cancellation of identical terms with opposite signs in
$ \dd (l_0,l_1,l_2,l_3) $,
we see that commutativity follows because $ f_{2,\F} $ is alternating.
\end{proof}

\begin{remark} \label{degenrem} The argument used to prove Lemma \ref{commutative lemma} provides the motivation for considering only tuples
$ (l_0,\dots,l_n) $ of elements in $ \LL $ such that if $ l_{i_1} $ and $ l_{i_2} $ have the same image
in $ \P_\F^1 $ then~$ l_{i_1} = l_{i_2} $.
It seems reasonable to define $ f_{2,\F} $ and $ f_{3,\F} $ to be trivial
on tuples for which some points have the same image in $ \P_\F^1 $.
Starting with such a tuple $ (A,A',B,C) $ where $ A $ and $ A' $
have the same image, but $ A $, $ B $ and $ C $ have different
images, we require $ f_{2,\F} $ takes the same value
on~$ (A,B,C) $ and $ (A',B,C) $ and so must limit the amount
of scaling between $ A $ and $ A' $ to~$ \tL $.
\end{remark}

\subsubsection{}We now set
\begin{equation*}
 \Bbar {\F} :=  \frac{\Z[\F^\flat]}{(f_{3,\F}\circ \dd ) (\pts 4 {\LL} ) }
\,.
\end{equation*}
 Then the diagram in Lemma \ref{commutative lemma} induces a commutative diagram
\begin{equation} \label{newCD}
\begin{tikzcd}
\cdots \ar{r}{\dd} &  \pts 4 {\LL} \ar{d} \ar{r}{\dd} &  \pts 3 {\LL}
\ar{d}{f_{3,\F}}\ar{r}{\dd} & \pts 2 {\LL} \ar{d}{f_{2,\F}}
\ar{r}{\dd} & \pts 1 {\LL} \ar{d} \ar{r}{\dd} & \pts 0 {\LL} \\
& 0  \ar{r}  &  \Bbar {\F} \ar{r}{\ddR {2,\F} {\tL} } & \twq {\tL} {\F^\times} \ar{r} & 0
\end{tikzcd}
\end{equation}
in which $ \ddR {2,\F} {\tL} $ denotes the map induced by~$ \d_{2,\F}^{\tL} $.
(If $ \tL $ is trivial then we shall use the notation~$ \ddR {2,L} {} $
for the induced map.)
We observe that we could take $ \GL_2(\F) $-coinvariants in
the top row because of the properties of $ f_{3,\F} $ and $ f_{2,\F} $.
In particular, $ f_{3,\F} $ induces a homomorphism~$ H_3(\pts {\bullet} {\LL} _{\GL_2(\F)}) \to \BB {\F} {\tL} $, where we set
\begin{equation*}
\BB {\F} {\tL} := \ker(\ddR {2,\F} {\tL} )
\,.
\end{equation*}
We shall denote this latter group more simply as $ \BB {\F} {} $ if $ \tL $ is trivial.

The following result provides an explicit and very useful description of the relations in $\Bbar {\F} $.

\begin{lemma}\label{explicit relations} The subgroup $(f_{3,\F}\circ \dd )(\pts 4 {\LL} )$ of $\Z[\F^\flat]$ is generated by all elements of the form
\begin{equation} \label{5-term}
 [x]-[y]+[y/x]-[(1-y)/(1-x)]+[(1-y^{-1})/(1-x^{-1})]
\end{equation}
for $ x \ne y $ in $ \F^\flat $, and
\begin{equation} \label{2-2-term}
 [x]+[x^{-1}] \text{ and } [y]+[1-y]
\end{equation}
for $x$ and $y$ in $\F^\flat$.
\end{lemma} 

\begin{proof} We note first that for each non-degenerate generator
$ (l_0,\dots,l_4) $ of $\pts 4 {\LL} $ one has
\[ (f_{3,\F}\circ \dd)((l_0,\dots,l_4)) = \sum_{i=1}^5 (-1)^i \crr_3(\ol{l_0},\dots,\widehat{\ol{l_i}},\dots,\ol{l_4}),\]
where $ \ol{l_0},\dots,\ol{l_4} $ are distinct points in  $ \P_\F^1 $ and $\widehat{\ol{l_i}}$ indicates that the term $\ol{l_i}$ is omitted.

In view of the invariance of $ \crr_3 $ under the action of $ \GL_2(\F) $,
and the fact that for $ \crr_3 $ we can use points in $ \P_\F^1 $,
we may assume the points are $ (1,0) $, $ (0,1) $, $ (1,1) $, $ (x,1) $
and $ (y,1) $ for $ x \ne y $ in~$ \F^\flat $. This then under $ f_{3,\F} \circ \dd $ yields the element (\ref{5-term}).

Let now $ (l_0,\dots,l_4) $ be a degenerate generator. Then its image under $ f_{3,\F} \circ \dd $
is trivial if~$ \{l_0,\dots,l_4\} $ has at most three elements since then all the terms in $ \dd (l_0,\dots,l_4) $ are degenerate.
On the other hand, if $ \{l_0,\dots,l_4\} $ has four elements, then after cancelling possible identical
terms in $ \dd (l_0,\dots,l_4) $ and applying $ \crr_3 $ to the
result we see that it is of the form
\[  [\crr_3(\ol{m_1},\dots,\ol{m_4})] - \sign(\sigma) [\crr_3(\ol{m_{\sigma(1)}},\dots,\ol{m_{\sigma(4)}})]\]
for a permutation $ \sigma $ in $ S_4 $
with sign $ \sign(\sigma) $, and four distinct points
$ \ol{m_i} $ in $ \P_\F^1 $. The subgroup generated by these images coincides with the subgroup generated by the terms (\ref{2-2-term}).
 (This shows, in particular, that the map $ f_{3,\F} $ is alternating.)
\end{proof}

\begin{remark} \label{factorremark}
(i)
If $ \tL $ is finite of order~$ a $, then multiplying an element in $ \BB {\F} {\tL} $ in the bottom row of~\eqref{newCD} by $ a $ gives an element in  $ \BB {\F} {} $.

\noindent(ii)
If~$ \s \colon F \to \CC $ is an embedding of a number field~$  $, then 
the map $ \Z[F^\flat] \to \R(1) $ in Theorem~\ref{regthm}(iii),
which maps $ [x] $ to $ i \Di(\sigma(x)) $, by~\eqref{dilogeqns} induces a map $ \tD_\s \colon \Bbar F \to \R(1) $.
\end{remark}

\subsubsection{}\label{connection to B}
We next show that if $ |F| \ge 4 $ then $ \BB {\F} {} $ as
defined above is naturally isomorphic to a quotient of the `Bloch group' $ B(\F) $ that is defined
and studied by Suslin in \cite{Sus90} if $ F $ is infinite
and treated in \cite[Chap. VI, \S5]{WeiKbook} for $ |F| \ge 4 $.
This result motivates us to regard~$ \BB {\F} {} $
as a modified Bloch group (which explains our choice of notation).
In fact, we shall establish the precise relation between our groups~$ \BB {\F} {} $ and $ \Bbar {\F} $ and the corresponding groups~$ B(\F) $ and~$ \pp(\F) $.

Following those two sources, for a field $ F $ with $ |\F| \ge 4 $ we set
its pre-Bloch group to be
\begin{equation*}
 \pp(\F) = \frac{\Z[\F^\flat]}{ \langle  [x]-[y]+\bigl[\frac{y}{x}\bigr]+\bigl[\frac{1-x}{1-y}\bigr] - \bigl[\frac{1-x^{-1}}{1-y^{-1}}\bigr]
 \text{ with } x, y \text{ in } \F^\flat, x \ne y
\rangle
 }.
\end{equation*}
We then define its Bloch group $ B(\F) $ to be the kernel of the homomorphism
\begin{equation} \label{gon boundary}
\begin{aligned} 
 \pp(\F) & \to \suswedge {\F^\times} 
\\
[x] & \mapsto \sustimes x (1-x)
\,,
\end{aligned}
\end{equation}
where we set
\[ \suswedge {\F^\times} :=  \frac{\F^\times\otimes \F^\times}{\langle x \otimes y + y
\otimes x \text{ with } x, y \text{ in } \F^\times \rangle}\]
and write $ \sustimes a b $ for the class in the quotient of an element $ a \otimes b $.

We further recall the existence of an exact sequence
\begin{equation} \label{tork3}
 0 \to \Tor(\F^\times, \F^\times)^\sim \to \Kind {\F} \to B(\F) \to 0,
\end{equation}
natural in $ \F $, where $ \Tor(\F^\times, \F^\times)^\sim $ denotes the unique non-trivial
extension of $ \Tor(\F^\times, \F^\times) $ by~$ \Z/2\Z $
for $ \F $ of characteristic different from~$2$ and $ \Tor(\F^\times, \F^\times) $ otherwise,
and~$K_3(\F)^{\rm ind}$ is the cokernel of the natural homomorphism from the Milnor $K$-group $K_3^M(\F)$ to~$K_3(\F)$.
We also recall that the element $ c_\F = [x]+[1-x] $ of $ B(\F) $ is independent of
$ x $ in~$ \F^\flat $ and has order dividing~6 by \cite[Lem.~1.3, 1.5]{Sus90}
or \cite[VI.5.4]{WeiKbook}, and that
$ c_\Q $ in $ B(\Q) $ has order~6~\cite[Prop.~1.1]{Sus90}.

\subsubsection{}Lemma \ref{explicit relations} implies that $ \Bbar {\F} $ is obtained by quotienting out~$ \pp(\F) $ by the subgroup generated
 by all elements of the form $ [x] + [x^{-1}] $ with $ x $ in $ \F^\flat $ and~$ [y]+[1-y] $ with $ y $ in $ \F^\flat $.

Since the latter  elements generate the same group as does the element $c_\F$ defined above we are motivated to consider the following short exact sequence of complexes (with vertical differentials)

\begin{equation*}
\begin{tikzcd}
0 \ar{r} & \langle [x]+ [x^{-1}] \text{ with } x \text{ in } \F^\flat \rangle \ar{r} \ar{d}{f}
& \pp(\F)/\langle c_\F \rangle \ar{r} \ar{d} & \Bbar {\F} \ar{r} \ar{d}{\ddR {2,\F} {} } & 0
\\
0 \ar{r} & \langle \sustimes x (-x) \text{ with } x \text{ in } \F^\flat \rangle \ar{r} & \suswedge {\F^\times} \ar{r} & \twt \F^\times \ar{r} & 0.
\end{tikzcd}
\end{equation*}

\begin{theorem} \label{ppmodcBbar} If $|\F| \ge 4$,
then the homomorphism $ f $ in the above
diagram is bijective, and, in particular, the diagram induces an isomorphism $ B(\F)/\langle c_\F \rangle \to \BB {\F} {} $.
\end{theorem}

\begin{proof}
A calculation shows that $ f $ maps the class of $ [x]+[x^{-1}] $ to $ \sustimes x (-x) $,
so $ f $ is surjective.

For proving that $ f $ is injective, recall that by \cite[Lem.~1.2]{Sus90} or \cite[VI.5.4]{WeiKbook}, the map
$ \F^\times \to \pp(\F) $ sending $ x $
to $ [x] + [x^{-1}] $ if $ x \ne 1 $ and $ 1 $ to $ 0 $, is a
homomorphism with $ (\F^\times)^2 $ in its kernel.
We shall consider its composition with the quotient map to $\pp(\F) / \langle c_\F \rangle $,
giving a surjective homomorphism
\begin{equation*}
 g \colon \F^\times \to \langle [x]+ [x^{-1}] \text{ with } x \text{ in } \F^\flat \rangle
 \,,
\end{equation*}
with the target in~$  \pp(\F) / \langle c_\F \rangle $.
If $ -1 $ is a square then we already know that $ [-1] + [-1] = 0 $
in~$ \pp(\F) $. If $ -1 $ is not a square then $ 2 \ne 0 $, so
$ 2 [-1] = 2 c_\F - 2 [2] = 2 c_\F + 2 [\frac12] = 3 c_\F $ in~$ \pp(\F) $
(cf.~\cite[Lem.~1.4]{Sus90} or~\cite[VI.5.4]{WeiKbook}). In either case, we have that
$ \{\pm1\} \cdot (\F^\times)^2 \subseteq \ker(g) $, and that
$ \im(g) $ is the subgroup generated by the classes
of $ [x] + [x^{-1}] $ with $ x $ in $ \F^\flat $.
We also want to consider $ \ker(f \circ g) $. For this,
we fix a basis $ \mathcal{B} $ of $ \F^\times/(\F^\times)^2 $ as $ \Ftwo $-vector space, making sure to include
$ -1 $ in $ \mathcal{B} $ if~$ -1 $ is not a square in $ \F^\times $.
For $ b $ in $ \mathcal{B} $, the homomorphism~$ \F^\times / (\F^\times)^2 \to \Ftwo \cdot b \simeq \Ftwo $
obtained from the projection onto $ \Ftwo \cdot b $ can be applied
twice in the tensor product in order to give a composite homomorphism~$ \F^\times \otimes \F^\times \to \F^\times/ (\F^\times)^2 \otimes \F^\times/ (\F^\times)^2 \to \Ftwo \otimes \Ftwo \simeq \Ftwo $.
This induces a homomorphism
$ h_b \colon \suswedge {\F^\times} \to \Ftwo $, mapping
$ \sustimes x y $ to the product of the coefficients of $ b $
in the classes of $ x $ and~$ y $ in $ \F^\times/(\F^\times)^2 $.
If $ x $ in $ \F^\times $ is in $ \ker(f \circ g) $, then
$ h_b(\sustimes x (-x)) = 0 $ for all~$ b $.
If $ -1 $ is a square, this means $ x $ is a square. If $ -1 $
is not a square, then $ x $ or $ -x $ must be a square.
In either case, it follows that $ \ker(f \circ g) \subseteq \{\pm1\} \cdot (\F^\times)^2 $.
Because $ \{\pm1\} \cdot (\F^\times)^2 \subseteq \ker(g) $ and $ g $
is surjective, it follows that $ f $ is injective.

Now that $ f $ is an isomorphism, the snake lemma implies the
isomorphism in the theorem.
\end{proof}

\begin{remark}
Note that in the proof above, it also follows that $ \ker(g) = \{\pm1\} (\F^\times)^2 $.
So $ g $ induces an isomorphism from $  \F^\times / \{\pm1\} \cdot (\F^\times)^2 $
to the subgroup of $ \pp(\F) / \langle c_F \rangle $ generated by the classes of $ [x] + [x^{-1}] $ with
$ x $ in $ \F^{\flat}$, given by mapping the class of $ x $ to the
class of $ [x] + [x^{-1}] $.
And~$ f \circ g $ induces an isomorphism from $  \F^\times / \{\pm1\} \cdot (\F^\times)^2 $
to the subgroup of $ \suswedge {\F^\times} $ generated by the $ \sustimes x (-x) $,
mapping the class of $ x $ to the class of $ \sustimes x (-x) $.
\end{remark}

\subsubsection{}We can now state the main result of this section.
It concerns the map $ \phi_\F $ in Theorem~\ref{regthm}.

\begin{theorem} \label{psiprop} For any field $ F $, the map $ \phi_\F $ induces a homomorphism $ \psi_\F \colon \BB {\F} {} \to \Kindtf {\F} $.
\end{theorem}

\begin{proof} In view of Lemma \ref{explicit relations}, it
is enough suffices to show that $ \phi_\F $
is trivial on all elements of the form~\eqref{5-term} and~\eqref{2-2-term}.

If $ \F $ is a number field $ F $ then this follows from Theorem~\ref{regthm}(iii), 
\eqref{dilogeqns}, and Borel's theorem, Theorem~\ref{borel}, by letting $ \sigma $ run through all
embeddings of $ F $ into~$ \CC $.

In order to see that it holds for all fields $ \F $ as in Theorem~\ref{regthm}, we can tensor with $ \Q $, in which case the construction underlying
the construction of the map~$ \phi_\F $ in Theorem~\ref{regthm} is the simplest case of the constructions that are made by the second author in~\cite{dJ95}.

One can then verify that the elements in~\eqref{2-2-term} and~\eqref{5-term}
are trivial by working over $ \Z[x,x^{-1},(1-x)^{-1}] $ or
$ \Z[x,y,(1-x)^{-1},(1-y)^{-1},(x-y)^{-1}] $ as the base schemes,
along the lines of the
proofs of \cite[Prop.~6.1]{dJ95} and \cite[Lem.~5.2]{dJ00a}.
We leave the precise details of this argument to an interested reader.
\end{proof}

\begin{remark} In this remark we let $ \tL $ be trivial and explain the advantages
of the definitions that we have adopted in comparison to those
used by Goncharov in \cite{gonXpam}.

For this,
we recall that in~(3.8) of loc.~cit.~a key role is played by the
map in~\eqref{gon boundary}
that sends a generator $ [x] $  to $ \sustimes x (1-x) $. Our group~$ \twt \F^\times $ is a quotient of $ \suswedge {F^\times} $ (cf.~the
diagram just before~Theorem~\ref{ppmodcBbar}), and
$ \d_{2,\F}^{} $ maps $ [x] $ to the inverse of the image of  $ \sustimes x (1-x) $ in~$ \twt \F^\times $.

Now the map that Goncharov constructs from
non-degenerate triples of non-zero points in~$ \F^2 $ to
the right hand side of (\ref{gon boundary}) is not itself $\GL_2(\F) $-equivariant since letting a matrix with determinant $ c $ act changes the result by the class of $ c \otimes (-c) $.
In addition, the calculation with the points
$ (a,0) $, $ (0,b) $, $ (1,1) $ and $ (xc,c) $
in the proof of Lemma~\ref{commutative lemma} would similarly
result in the element
$ \sustimes x (1-x) + \sustimes c (-c) $ which is not what one wants.

Whilst these problems could be simply resolved by multiplying any of the relevant maps by a factor of two, this would in the end lead
either to a smaller subgroup of $ \Kindtf {\F} $ if we multiply~$ f_{3,\F} $ by~2, or a (new) Bloch group that is too large (if
we multiply Goncharov's boundary map by~2; cf.~\cite{neu-yang99} and many other papers).
It is therefore better
to avoid the problem by replacing the right hand side of (\ref{gon boundary}) as the target of the boundary map by its quotient~$ \twt \F^\times $.

But the elements of the form $ [x] + [x^{-1}] $ that are in the kernel of $ \d_{2,\F}^{} $
could
then result in a potentially large and undesired subgroup in
the kernel of the boundary map, even modulo the 5-term relations~\eqref{5-term}
(see Theorem~\ref{ppmodcBbar} and its proof). To avoid this,
we have also imposed
the relations~\eqref{2-2-term} when defining~$ \Bbar F $ by working with degenerate configurations.
\end{remark}

\subsection{Torsion elements in Bloch groups}
In this section we study the torsion subgroup of the modified Bloch group $ \BB F {} $ of a number field $ F $ by means of a comparison with the Bloch group $ B(F) $ defined by Suslin (and recalled in \S\ref{connection to B}).

In this way, we find that $ \BB F {} $ is torsion free if
$ F $ is equal to either $ \Q $, or to an imaginary quadratic number field or is generated over $\Q$ by a root of unity (of any given order).
In the case of imaginary quadratic fields, this fact will then play an important role in \S~\ref{tessellations}.

\subsubsection{}For the sake of simplicity, we formulate and prove the next result only for number fields.

In its statement, if $ p $ is a prime number, then we denote
the $ p $-primary torsion subgroup of a finitely generated abelian group by
means of the subscript $ p $.
Because all torsion groups here are finite and cyclic, this determines
their structures.

\begin{proposition} \label{torsionprop}
Let $ F \subset \Qbar $ be a number field.
For a prime $ p $, let $ p^s $ be the number of $ p $-power roots of unity in $ F $,
and let $ r $ be the largest integer such that the maximal
totally real subfield~$ \Q(\mu_{p^r})^+ $ of $ \Q(\mu_{p^r}) $
is contained in $ F $.
Then the orders of the $ p $-power torsion subgroups in the various groups are
as follows.

\begin{center}
\def\vph{\vphantom{$ b^{b^b}$}}
\begin{tabular}{l|c|c|c|c|c}
prime & $ | \Tor(F^\times, F^\times)_p^\sim| $ & $ |\Kind F _p| $ & $ |B(F)_p| $  &  $ |\BB F {p} |$ & condition
\\
\cline{1-6}
\multirow{2}{*}{$ p \ge 5 $} & \multirow{2}{*}{$ p^s $}     & \multirow{2}{*}{$ p^r $}     & $ p^r $     & $ p^r $     & $ \zeta_p \notin F $\vph
\\
\cline{4-6}
                             &                              &                              &    $1$      &    $1$      & $ \zeta_p \in F $\vph
\\
\cline{1-6}
\multirow{2}{*}{$ p = 3 $}   & \multirow{2}{*}{$ 3^s $}     & \multirow{2}{*}{$ 3^r $}     & $ 3^r $     & $ 3^{r-1} $ & $ \zeta_3 \notin F $\vph
\\
\cline{4-6}
                             &                              &                              &    $1$      &    $1$      & $ \zeta_3 \in F $\vph
\\
\cline{1-6}
\multirow{2}{*}{$ p = 2 $}   & \multirow{2}{*}{$ 2^{s+1} $} & \multirow{2}{*}{$ 2^{r+1} $} & $ 2^{r-1} $ & $ 2^{r-2} $ & $ \zeta_4 \notin F $\vph
\\
\cline{4-6}
                             &                              &                              &    $1$      &    $1$      & $ \zeta_4 \in F $\vph
\\
\cline{1-6}
\end{tabular}
\end{center}

\medskip

\noindent
\textup{(}Note $ r \ge 2 $ and $ s \ge 1 $ if $ p =2 $, and $ r \ge 1 $ if $ p = 3 $.\textup{)}
\end{proposition}

\begin{proof}
We compute $ |\Kind F _p |  $ (which is faster than
using \cite[Chap. IV, Prop.~2.2 and~2.3]{WeiKbook}).

Let~$ A \subseteq \Z_p^\times $ be the image of $ \Gal(\Qbar/F) $ in $ \Gal( \Q(\mu_{p^\infty})/\Q) \simeq \Z_p^\times $.
 Then there are identifications
\[  \Kind F _p \simeq H^0(\Gal(\Qbar/F), \Q_p(2)/\Z_p(2)) \simeq \cap _{a \in A} \ker( \Q_p/\Z_p \overset{a^2-1}\to \Q_p/\Z_p),\]
where the first follows from \cite[Cor.~4.6]{lev89} and the second is clear.

We assume for the moment that $  p \ne 2 $. Then $ r = 0 $ is equivalent with $ A \not \subseteq \{\pm1\} \cdot (1+p \Z_p) $,
so some $ a^2 -1  $ is in $ \Z_p^\times $ and the resulting kernel is trivial.
For $ r \ge 1 $, we have $ A \subseteq \{\pm 1\} \cdot (1+ p^r \Z_p) $
but $ A \not \subseteq \{\pm 1\} \cdot (1+ p^{r+1} \Z_p) $.
Then $ 1 + p^r \Z_p \subseteq A $ because $ p $ is odd, hence $ A^2 = 1 + p^r \Z_p $ and the statement is clear.

To deal with the case $ p = 2 $ we note that $ A $ is the image of $ \Gal(\Qbar/F) $
in $ \Gal( \Q(\mu_{2^\infty})/\Q) \simeq \Z_2^\times = \{\pm1\} \cdot (1+4 \Z_2) $.
Then $ A \subseteq \{\pm 1\} \cdot (1+ 2^r \Z_2) $
but $ A \not \subseteq \{\pm 1\} \cdot (1+ 2^{r+1} \Z_2) $, for
some $ r \ge 2 $. In this case, $ A $ contains an element of
$ \{\pm1\} \cdot (1 + 2^r \Z_2^\times) $, and $ A^2 = 1 + 2^{r+1} \Z_2 $, from which
the statement follows.

We always have $ r \ge s $.
If $ \zeta_{2p} $ is in $ F $ then $ r = s $ because $ \Q(\zeta_{p^r})^+(\zeta_{2p}) = \Q(\zeta_{p^r}) $.
If $ \zeta_{2p} $ is not in $ F $, then $ s = 0 $ for $ p \ne 2 $,
and $ s = 1 $ for $ p = 2 $.
The entries for $ |B(F)_p| $ are now immediate from~\eqref{tork3}.
From this, we recover that $ B(\Q) $ has order~6, so is generated
by $ c_\Q $. Then we can compare the sequences~\eqref{tork3} for
the field $ \Q $ and for~$ F $. Using that $ \Kind F _\tor $
is cyclic, and that $ c_\Q $ maps to $ c_F $ under the injection
$ \Kind {\Q} _\tor \to \Kind F _\tor $,
it follows that $ 3 c_F $ has order~2 if and only
if $ \zeta_4 \notin F $, and that $ 2 c_F $ has order~3 if and
only if $ \zeta_3 \notin F $: in order to have those orders,
we must have $ | \Tor(F^\times, F^\times)_p^\sim| = | \Tor(\Q^\times, \Q^\times)_p^\sim|  $
for $ p = 2 $ or $ p = 3 $ respectively.
(Cf.\ \cite[Lem.~1.5]{Sus90}.)
This gives the entries for $ | \BB F {p} | = | ( B(F) / \langle c_F  \rangle )_p | $.
\end{proof}

If $ F $ is a number field, and $ p $ a prime number, then
combining Theorem~\ref{ppmodcBbar} and Proposition~\ref{torsionprop}
gives a precise statement when $ \BB F {} $ has no $ p $-torsion.

\begin{theorem} \label{Bthm}
Let $ F \subset \Qbar $ be a number field. Then, for a given prime number $p$, the group~$ \BB F {} $ has no $ p $-torsion if and only if the following condition is satisfied:

\begin{itemize}
\item[-] if $ p \ge 5 $, then either $ \Q(\zeta_p)^+ \not \subseteq F $ or $ \Q(\zeta_p) \subseteq F $;

\item[-] if $ p = 3 $, then either $ \Q(\zeta_9)^+ \not \subseteq F $ or $ \Q(\zeta_3) \subseteq F $;

\item[-] if $ p = 2 $, then either $ \Q(\sqrt2) = \Q(\zeta_8)^+ \not \subseteq F $ or $ \Q(\zeta_4) \subseteq F $.
\end{itemize}

\end{theorem}

\begin{proof}
Using Theorem~\ref{ppmodcBbar}, it is clear when $ \BB F {} $ has non-trivial $ p $-torsion.
With notation as in Proposition~\ref{torsionprop}, this
is the case if and only if the following condition is satisfied:
\begin{itemize}
\item[-] if $ p \ge 5 $, then $ r \ge 1 $ and $ \zeta_p $ is not in $ F $;

\item[-] if $ p = 3 $, then $ r \ge 2 $ and $ \zeta_3 $ is not in $ F $;

\item[-] if $ p = 2 $, then $ r \ge 3 $ and $ \zeta_4 $ is not in $ F $.
\end{itemize}
Upon negating this statement for each prime $ p $ one obtains the claimed result.
\end{proof}

\begin{corollary} \label{Bcor}
We have that

\begin{itemize}
\item[(i)] $ \BB {\Q} {} $  is trivial;

\item[(ii)] $ \BB F {} \simeq \Z^{[F:\Q]/2} $ if $ F = \Q(\zeta_N) $ with $ N \ge 3 $;

\item[(iii)] $ \BB F {} \simeq \Z $ if $ F $ is an imaginary quadratic field.
\end{itemize}

In addition, for each of these fields the composition of the maps~$ \psi_F \colon \BB F {} \to \Kindtf F $
and $ \Kindtf F \to \prod_{\s} \R(1) $, where $ \s $ runs through
the places of $ F $, is injective.
\end{corollary}

\begin{proof}
We first let $ F $ be any number field. Then $ \BB F {} $ is a finitely generated abelian group of the same rank as $ K_3(F) $
by~\eqref{tork3} and Theorem~\ref{ppmodcBbar},
hence, by Theorem~\ref{regthm}(ii), the kernel of~$ \psi_F $ is the torsion subgroup of $ \BB F {} $.
Because of the behaviour of the regulator with respect to complex conjugation, in Theorem~\ref{borel}(iii)
we only have to consider all places of $ F $, not all embeddings into~$ \CC $.

It therefore suffices to check that for the number fields listed in
the corollary, $ \BB F {} $ has no $ p $-torsion
for every prime number~$ p $. But this follows from Theorem~\ref{Bthm}.
\end{proof}

\subsubsection{}
 We conclude this subsection with a result on $ \Bbar {\Q} $
that we shall use in Sections~\ref{tessellations} and~\ref{proofs}.
 
\begin{proposition} \label{BbarQprop}
\mbox{}
\begin{itemize}
\item[(i)] The torsion subgroup of $ \Bbar {\Q} $ has order~$2$ and is generated
by~$ [2] $.

\item[(ii)] If $ k $ is an imaginary quadratic number field, then the
natural map $ \Bbar {\Q} \to \Bbar k $
and its composition with $ \ddR {2,k} {} $ are injective when
$ k \ne \k -1 $, but for $ k = \k -1 $ both kernels are generated by~$ [2] $.
\end{itemize}
\end{proposition}

\begin{proof} To prove claim (i) we note Corollary~\ref{Bcor} implies
that $ \Bbar {\Q} $ injects into $ \twt \Q^\times $. We also know
from Proposition~\ref{twdescription} (or Corollary~\ref{twcomputecor})
that the natural map $ \langle -1 \rangle \otimes \Q^\times \to \twt \Q^\times $
gives an isomorphism with the torsion subgroup of the latter.

So we want to compute the kernel of the natural map
$ \langle -1 \rangle \otimes \Q^\times \to K_2(\Q) $. Using the
tame symbol, and the fact that $ \{-1,-1\} $ is non-trivial
in $ K_2(\Q) $, one sees that this kernel is cyclic of order~$2$,
generated by $ (-1) \tw 2 = \ddR {2,\Q} {} ([2]) $.
And $ 0 = [\frac12] + [1-\frac12] = 2 [\frac12] = -2 [2]$.

Turning to claim (ii), we note that the kernel of the composition
by Corollary~\ref{Bcor}(i) under~$ \ddR {2,\Q} {} $
must inject into the kernel of $ \twt \Q^\times \to \twt k^\times $,
which we computed in Lemma~\ref{kertwmap}.
In particular, as this kernel is a torsion group, we see from Corollary~\ref{Bcor}(ii) that the kernel
of the composition and that of $ \Bbar {\Q} \to \Bbar k $ coincide
as the torsion of $ \Bbar k $ injects into $ \twt k^\times $ under~$ \ddR {2,k} {} $.

In addition,  by~claim (i), those kernels are either trivial or generated
by $ [2] $.  They contain $ [2] $ if and only
if $ (-1) \tw 2 $ is in the kernel of $ \twt \Q^\times \to \twt k^\times $,
which by Lemma~\ref{kertwmap} holds if and only if $ k = \k -1 $.

This completes the proof.
\end{proof}

\begin{remark} Note that $ [2] = 0 $ in $ \Bbar {\k -1 } $
follows explicitly from~\eqref{5-term} with $ x = \rtmo $, $ y = -\rtmo $,
which gives $ [-1] = 0 $, as $ [2] = -[-1] $ by~\eqref{2-2-term}.
\end{remark}

\subsection{A conjectural link between the groups of Bloch and Suslin} \label{bloch-suslin-link}
If $ \F $ is an infinite field, then~\eqref{tork3} gives an isomorphism $ \Kindtf {\F} \rightiso B(\F)_\tf $ and Theorem~\ref{ppmodcBbar} gives an isomorphism $ B(\F)_\tf \rightiso \BB {\F} {\tf} $.
(We ignore the case of finite fields as then all these groups are trivial.)

By Theorem~\ref{psiprop}, one also knows that the homomorphism $ \phi_\F $ in Theorem~\ref{regthm} induces a homomorphism of the form
$ \psi_\F \colon \BB {\F} {} \to \Kindtf {\F} $.
This in turn induces a homomorphism
\[ \psi_{\F,\tf} \colon \BB {\F} {\tf} \to K_3(\F)^{\rm ind}_{\rm tf}, \]
thereby allowing us to form the composition
\begin{equation}\label{unknown}
 \Kindtf {\F} \rightiso B(\F)_\tf \rightiso \BB {\F} {\tf} \buildrel{\psi_{\F,\tf}}\over{\longrightarrow} \im(\psi_{\F,\tf}) \subseteq
 \Kindtf {\F} \,.
\end{equation}
Both $ B(F) $ and $ \im(\psi_{\F,\tf}) $ are described as the kernel of a map from an
abelian group that is
generated by elements of the form $ [x] $ for $ x $ in $ \F^\flat $, sending each $ [x] $ to the class of $ (1-x) \otimes x $ in either $ \F^\times \otimes \F^\times $ or a variant like $ \twt \F^\times $.
As mentioned in Section~\ref{general-context}, these groups are
widely expected to be closely related even though there is no
obvious map between them, the former being
 constructed using group homology of $ \GL_ 2(\F) $ and the latter using
relative $ K $-theory.

Our approach provides the first concrete evidence (in situations in which the groups are non-trivial) to suggest
that~$ B(F) $ and $ \im(\psi_{F,\tf}) $ should be related in a very natural way, and that
the latter is all of $ \Kindtf F $.
To be specific, if $ \F $
is a number field,
then $ \psi_{F,\tf} $ is injective by the
proof of Theorem~\ref{Bthm}, so by
Proposition~\ref{regthm} the composite map (\ref{unknown}) is
an injection of a finitely generated free abelian group into itself.
One can therefore determine the (finite) index of~$ \im(\psi_{F,\tf}) $  
in~$ \Kindtf F $ by comparing the results of the regulator map on $ \im(\psi_{F,\tf}) $ and on~$ \Kindtf F $.
Extensive evidence that we have obtained by computer calculations in the case that $F$ is imaginary quadratic (cf. \S\ref{generators})
motivates us to formulate the following conjecture.

\begin{conjecture}\label{psi is iso}
If $F$ is a number field, then $\psi_{F,\tf}$ is an isomorphism. \end{conjecture}

\begin{remark}\label{conj rem} As mentioned above, the map $ \psi_{F,\tf} $ is injective and so the main
point of Conjecture \ref{psi is iso} is that $ \im(\psi_{F,\tf}) = \Kindtf {F} $.
However, for a general number field~$ F $,
this equality would not itself resolve the problem of finding an explicit description of the
resulting composite isomorphism in~\eqref{unknown}.
Of course, for $ F $ imaginary quadratic
all groups occurring in the composite are isomorphic to~$ \Z $,
and so the conjecture would imply that (\ref{unknown}) is multiplication by~$\pm 1$. (Recall
that $ \psi_{F,\tf} $ is itself natural up to a universal choice of sign
since this is true for $ \phi_F $.)

It would also seem reasonable to hope that~\eqref{unknown} has a very simple description
for any infinite field~$ \F $, such as, perhaps, being given by multiplication by
some integer that is independent of~$ \F $. Assuming this to be the case, our numerical calculations would imply that this integer is $ \pm 1 $. If true, this would in turn imply that the isomorphism $\Kindtf {\F} \to B(\F)_\tf$ constructed by Suslin in \cite{Sus90} could be given a more direct, and more directly $K$-theoretical, description, at least up to sign, as the inverse of the composite isomorphism  $ B(\F)_\tf \to \BB {\F} {\tf} \to \Kindtf {\F} $ where the first map is induced by  Theorem~\ref{ppmodcBbar} and the second is~$\psi_{F,\tf}$.
\end{remark}

\section{A geometric construction of elements in the modified Bloch group} \label{tessellations}

Let $ k $ be an imaginary quadratic number field with ring of
algebraic integers $ \O $.

In this section, we shall use a geometric construction, the Voronoi theory of
Hermitian forms, in order to
construct a non-trivial element $ \bgeo $ in $ \BB k {} \simeq \Z $.

To do this we shall invoke a tessellation of hyperbolic 3-space for~$ k $,
based on perfect forms, in order to construct an element of the kernel of the homomorphism $\dd \colon \pts 3 {\LL} \to \pts 2 {\LL} $
 that occurs in ~\eqref{eq:configuration-complex} with $ F = k $ and $ \tL = \{1\} $.
By applying $ f_{3,k} $ to this element we shall then obtain $ \bgeo $ by using the commutativity
of the diagram ~\eqref{newCD}.

Furthermore, we are able to explicitly determine the image of this element under the regulator map
and compare it to the special value $\zeta_k'(-1)$ by using a celebrated formula of Humbert.
This will in particular show that the element $ \psi_k(\bgeo) $
of $\Kindtf k $ that is constructed in this geometric fashion
generates a subgroup of index $|K_2 (\O)|$ (cf.~Corollary \ref{gammageozeta}(i)).

\subsection{Voronoi theory of Hermitian forms}

Our main tool is the polyhedral reduction theory for $\GL_2(\OO)$
developed by Ash \cite[Chap.~II]{AMRT_SmoothCompactification} and
Koecher \cite{koecher}, generalizing work of Voronoi \cite{VoronoiI}
on polyhedral reduction domains arising from the theory of perfect
forms (see~\cite[\S 3]{Yasaki} and~\cite[\S 2, \S 6]{imquad-coh}
for a description of the algorithms involved).
We recall some details here to set notation.

We fix a complex embedding $k \hookrightarrow \CC$ and
identify $k$ with its image. We extend this identification to
vectors and matrices as well.  We use $\overline{\cdot}$ to denote
complex conjugation on $\CC$, the non-trivial Galois automorphism on
$k$.  Let $V = \herm(\CC)$ be the $4$-dimensional real vector space of
$2 \times 2$ complex Hermitian matrices with complex coefficients.
Let $C \subset V$ denote the codimension~0 open cone of positive
definite matrices.  Using the chosen complex embedding of $k$, we can
view $\herm(k)$, the $2 \times 2$ Hermitian matrices with coefficients
in $k$, as a subset of $V$.  Define a map~$q \colon \OO^2 \setminus
\set{0} \to \herm(k)$ by $q(x) = x \overline{x}^t$.  For each $x \in
\OO^2$, we have that $q(x)$ is on the boundary of~$C$.  Let $C^*$
denote the union of $C$ and the image of~$q$.

The group $\GL_2(\CC)$ acts on $V$
by $g \cdot A = g A \overline{g}^t$.   The image of $C$ in the
quotient of $V$ by positive homotheties can be identified with
hyperbolic $3$-space $\Hy$.
The image of $q$ in this quotient is identified with $ \P_k^1 $,
the set of cusps.
The action induces an action of $\GL_2(k)$ on $\Hy$ and the cusps of
$\Hy$ that is compatible with other models of $\Hy$
(see \cite[Chap.~1]{EGMbook} for descriptions of other models). 
We let  $\Hy^* = \Hy \cup \PP^1_k$.

Each $A \in V$ defines a Hermitian form
$A[x] = \overline{x}^tAx$, for $x \in \CC^2$.  Using the chosen
complex embedding of $k$, we can view $\OO^2$ as a subset of $\CC^2$.
\begin{definition}
  For $A \in C$, we define the \emph{minimum of $A$} as
\[\min(A) \coloneqq  \inf_{x \in \OO^2 \setminus \set{0}} A[x].\]
Note that $\min(A)>0$ since $A$ is positive definite.
A vector $v \in \OO^2$ is called a
\emph{minimal vector of $A$} if $A[v] = \min(A)$.
We let $\Min(A)$ denote the set of minimal vectors of $A$.
\end{definition}

These notions depend on the fixed choice of the imaginary quadratic
field $k$.  Since $ \OO^2 $ is discrete in the topology of $ \CC^2 $,
a compact set $ \{z \mid  A[z]\le\text{bound}\} $ in $ \CC^2 $
gives a finite set in $ \OO^2 $.
Thus $ \Min(A) $ is finite.

\begin{definition}
We say a Hermitian form $A \in C$ is a \emph{perfect Hermitian
form over $k$} if
\[\Span_\RR\set{q(v) \mid v \in \Min(A)} = V.\]
\end{definition}

By a \emph{polyhedral cone} in $V$ we mean a subset $\sigma$ of the form
\[\sigma = \set*{\sum_{i = 1}^n \lambda_i q(v_i) \mid \lambda_i \geq
  0},\]
where $v_1, \dots, v_n$ are non-zero vectors in $\OO^2$.  A set of polyhedral
cones $S$ forms a \emph{fan} if the following two conditions hold.
Note that a face here can be of codimension higher than~$1$.
\begin{enumerate}
\item If $\sigma$ is in $S$ and $\tau$ is a face of $ \sigma $, then
  $\tau$ is in $S$.
\item If $\sigma$ and $\sigma'$ are in $S$, then $\sigma \cap \sigma'$ is a
  common face of $\sigma$ and $\sigma'$.
\end{enumerate}
The reduction theory of Koecher \cite{koecher} applied in this setting
gives the following theorem.

\begin{theorem}
  There is a fan $\tilde{\Sigma}$ in $V$ with $\GL_2(\OO)$-action such
  that the following hold.
  \begin{itemize}
  \item[(i)] There are only finitely many $\GL_2(\OO)$-orbits in
    $\tilde{\Sigma}$.
  \item[(ii)] Every $y \in C$ is contained in the interior of a unique cone
    in $\tilde{\Sigma}$.
  \item[(iii)] Any cone $\sigma \in \tilde{\Sigma}$ with non-trivial
    intersection with $C$ has finite stabiliser in
    $\GL_2(\OO)$. 
  \item[(iv)] The $4$-dimensional cones in $\tilde{\Sigma}$ are in bijection with the perfect forms over $k$.
\end{itemize}
\end{theorem}

The bijection in claim (iv) of this result is explicit and allows one to compute the structure of~$\tilde{\Sigma}$ by using a modification of Voronoi's algorithm
\cite[\S 2, \S 6]{imquad-coh}.  Specifically, $\sigma$  is a
$4$-dimensional cone in $\tilde{\Sigma}$ if and only if
there exists a perfect Hermitian form $A$ such that
\[\sigma = \set*{\sum_{v \in \Min(A)} \lambda_v q(v) \given \lambda_v
  \geq 0}.\]
Modulo positive homotheties, the fan $\tilde{\Sigma}$ descends to a
$\GL_2(\OO)$-tessellation of $\Hy$ by ideal polytopes.
The output of the computation described above is a collection of
finite sets $\Sigma^*_n$, $n = 1, 2, 3$, of
representatives of the $ \GL_2(\O) $-orbits of the  $n$-dimensional cells in $\Hy^*$
that meet~$\Hy$. The cells in each $ \Sigma_n^* $
have vertices described explicitly by finite sets of non-zero vectors in~$\OO^2$.

\subsection{Bloch elements from ideal tessellations of hyperbolic space}
The collection of 3-cells
\[\Sigma^*_3 = \set{P_1, P_2, \dots, P_m}\]
above gives
rise, after choosing a triangulation of each, to an element in  $ \BB k {} $, as follows.

We first establish a useful interpretation of a classical formula of Humbert in this setting.
For the sake of brevity, we shall write $ \C $ for $\PGL_2(\OO)$.

\begin{lemma} \label{FDvol}  Let $\C_{P_i}$ denote the stabiliser in $ \C $ of $P_i$. Then one has
\[ \sum_{i = 1}^m\frac{1}{|\C_{P_i}|}\vol(P_i) = - \pi\cdot \zeta_k'(-1)\,.\]
\end{lemma}

\begin{proof} One has
\begin{equation}\label{humbert} \sum_{i = 1}^m\frac{1}{|\C_{P_i}|}\vol(P_i) = \vol(\PGL_2(\OO)\backslash \Hy) = \frac{1}{8\pi^2}|D_k|^{\frac{3}{2}}\cdot \zeta_k(2).\end{equation}
Here the first equality is clear and the second is a celebrated result of Humbert (see ~\cite{borel}, where the formula is given for general number fields).
The claimed formula now follows since an analysis of the functional equation (\ref{dedekind fe}) shows that the final term in (\ref{humbert}) is equal to $-\pi \cdot\zeta_k'(-1)$.
\end{proof}

We next subdivide each polytope $P_i$ into ideal tetrahedra $T_{i,j}$ with
positive volume without introducing any new vertices,
\begin{equation} \label{subdiv}
P_i = T_{i,1} \cup T_{i,2} \cup \dots \cup T_{i,{n_i}}
\,.
\end{equation}
Here we assume that the subdivision is such that the
faces of the tetrahedra that lie in the interior of the $ P_i $ match.
An ideal tetrahedron $T$ with vertices $v_1, v_2, v_3, v_4$
has volume 
\[ \vol(T) =  \Di(\crr_3(v_1, v_2, v_3, v_4)).\]
Here $D$ denotes the Bloch-Wigner dilogarithm defined in \S\ref{bloch-wigner} above
and $\crr_3$ denotes the cross-ratio discussed in \S\ref{cr3 sec def}.
The ordering of vertices is chosen so that the right hand side is positive.

To ease the notation, we let $ r_{i,j}$ denote a resulting cross-ratio for $T_{i,j}$. We note that, whilst there is some ambiguity in choosing the order the four vertices of $T_{i,j}$ when defining this cross-ratio, the transformation rules in Remark~\ref{classicalcr} combine with
the relations in~\eqref{2-2-term} to imply that the induced element $ [ r_{i,j} ] $ of $ \Bbar k $ is
indeed independent of that choice.

We can now formulate the main result of this section (the proof of which will be given in \S\ref{proofs}).

By Corollary~\ref{cor:lcm} we know that each $ |\C_{P_i}|$ divides~24,
so the coefficients in the next theorem are integers.
We also note that, by Proposition~\ref{BbarQprop}, the map $ \Bbar {\Q} \to \Bbar k $
is injective unless $ k = \k -1 $, that the map $ 2\cdot \Bbar {\Q} \to \Bbar k $
is always injective, and that $ 2\cdot \Bbar {\Q} $ is torsion-free.
Moreover, the composition of the map $ \Bbar {\Q} \to \Bbar k $ with $ \ddR {2,k} {} \colon \Bbar k \to  \twt k^\times $
is injective if $ k \ne \k -1 $, and if $ k = \k -1 $ then this
composition has the same kernel as the map $ \Bbar {\Q} \to \Bbar k $.
Therefore the image of $ \Bbar {\Q} $ in $ \Bbar k $ always injects
into $ \twt k^\times $ under~$ \ddR {2,k} {} $.

\begin{theorem}\label{mthm} Let $ k $ be an imaginary quadratic number field, with the
polytopes $ P_i $ and cross-ratios $ [r_{i,j}] $ chosen as above. Then the following hold.
\begin{itemize}
\item[(i)] There exists a unique element $ \b_\Q $ in the image of $ \Bbar {\Q} $ in $ \Bbar k $,
such that the element
\begin{equation*}
    \bgeo = \b_\Q + \sum_{i = 1}^m  \frac{24}{|\C_{P_i}|}  \sum_{j = 1}^{n_i} [r_{i,j}]
\end{equation*}
belongs to $ \BB k {} $.
If $ k \ne \k -2  $ then $ \b_\Q $ belongs to the image of $ 2\cdot \Bbar {\Q} $. In all cases the element $ \bgeo $ is independent of the choice of representatives in $ \Sigma_3^* $, and the resulting
subdivision~\eqref{subdiv} into tetrahedra.

\item[(ii)]  If no stabiliser of an element in $ \Sigma_3^* $ or $ \Sigma_2^* $
has order divisible by~$4$, then there is a unique $ \tilde \b_\Q $
in $ \Bbar {\Q} $, which lies in $ 2\cdot \Bbar {\Q} $, such that the element
\begin{equation*}
\tbgeo = \tilde \b_\Q + \sum_{i = 1}^m  \frac{12}{|\C_{P_i}|}  \sum_{j = 1}^{n_i} [r_{i,j}]
\end{equation*}
belongs to $ \BB k {} $. 
Moreover, one has $ 2 \cdot\tbgeo = \bgeo $ and $ 2\cdot \tilde \b_\Q = \b_\Q $.
\end{itemize}
\end{theorem}

\begin{remark}
The situation for $ k $ equal to either $ \k -1 $ or $\k -2 $ is more complicated because
the order of the stabiliser of the (in both cases unique) element of $ \Sigma_3^* $
has order~24. For $ k = \k -2 $ it can be subdivided in several
different ways, resulting in the exception in Theorem~\ref{mthm}(i).
In fact, the subdivision in this case determines whether $ \b_\Q $ either belongs, or does not belong, to $ 2\cdot \Bbar {\Q} $, and both cases occur; see the argument in \S\ref{proofs} below for more details.
\end{remark}

\begin{remark}\label{mpropnew} It is sometimes computationally convenient to avoid explicitly computing the element  $ \b_\Q $ in Theorem \ref{mthm}(i). In this regard it is useful to note that the injectivity in~Corollary~\ref{Bcor} combines with Theorem~\ref{regthm}(iii) and the equality $ \Di(\ol{z}) = - \Di(z) $ in~\eqref{dilogeqns} to imply that
\begin{equation*}
 2\cdot \bgeo = \sum_{i = 1}^m \frac{24}{|\C_{P_i}|}  \sum_{j = 1}^{n_i} ( [r_{i,j}] - [\ol{r_{i,j}}] ).
\end{equation*}
\end{remark}

\subsection{Regulator maps and \texorpdfstring{$K$}{K}-theory}
As we fixed an injection of $ k $ into $ \CC $, by the behaviour
of the regulator map $ \reg_2 $ with respect to complex conjugation (see \eqref{dilogeqns}),
we can compute regulators by considering only the
composition  \[ K_3(k) \to K_3(\CC) \buildrel{\reg_2}\over\longrightarrow \R(1) .\]
By slight abuse of notation, we shall denote this composition
by~$ \reg_2 $ as well.

\begin{corollary} \label{gammageozeta} Assume the notation and hypotheses of Theorem~\ref{mthm}. Then the following hold.
\begin{itemize}
\item[(i)] The element $\psi_k(\bgeo) $ satisfies
\[
\frac{\reg_2(\psi_k(\bgeo))}{2 \pi i} = -12 \cdot\zeta_k'(-1)
\,
\]
and generates a subgroup of the infinite cyclic group $ \Kindtf k $ of index $ |K_2(\OO)|$.

\item[(ii)] If no stabiliser of an element in $ \Sigma_3^* $ or $ \Sigma_2^* $
has order divisible by~$4$, then $|K_2(\O)|$ is even, and $\psi_k(\tbgeo) $ generates a subgroup 
of  $ \Kindtf k $ of index $ |K_2(\OO)|/2$.
\end{itemize}
\end{corollary}

\begin{proof} Before proving claim (i) we note that for each polytope $P_i$ in $\Sigma^*_3$ one has
\begin{equation}\label{summary} \rtmo\cdot\vol(P_i) = \sum_{j=1}^{n_i}\tD([r_{i,j}]),\end{equation}
where $\tD$ is the homomorphism $\Bbar k \to \R(1)$ that is defined in Remark~\ref{factorremark}(ii) with respect to a fixed embedding $k\to \CC$.
This is true because $\vol(P_i) = \sum_{j=1}^{n_i}\vol(T_{i,j})$ whilst for each $i$ and $j$ one has $\rtmo\cdot\vol(T_{i,j}) = \rtmo\cdot D(r_{i,j}) = \tD([r_{i,j}])$.

Turning now to the proof of claim (i), we observe that, since the element $ \b_\Q $ that occurs in the definition of $\bgeo$ lies in the image of the map $ \Bbar {\Q} \to \Bbar k $, it also lies in the kernel of the composite homomorphism $\reg_2\circ \psi_k$. One therefore computes that
\begin{align*}
\reg_2(\psi_k(\bgeo)) & = \sum_{i = 1}^m  \frac{24}{|\C_{P_i}|}  \sum_{j = 1}^{n_i} \reg_2(\psi_k([r_{i,j}]))\\
                                    & = 24\cdot\sum_{i = 1}^m  \frac{1}{|\C_{P_i}|}  \sum_{j = 1}^{n_i} \tD([r_{i,j}])\\
                                    & = 24\rtmo\cdot \sum_{i = 1}^m  \frac{1}{|\C_{P_i}|}\vol(P_i)\\
                                    & = -24\pi\rtmo\cdot \zeta_k'(-1)
\,,
\end{align*}
where the second equality follows from Theorem \ref{regthm}(iii) and Remark~\ref{factorremark}(ii), the third from~(\ref{summary}) and the last  from Lemma \ref{FDvol}.
This proves the first assertion of claim (i),
and then the final assertion of claim (i) follows directly from
Example~\eqref{klichtenbaum}.

For claim (ii) we note that, under the stated conditions, Theorem \ref{mthm}(ii) implies~$ \bgeo = 2\cdot \tbgeo $,
so this follows from the final assertion of claim~(i).
\end{proof}

\begin{remark}
(i)
The condition in Theorem \ref{mthm}(ii) and Corollary~\ref{gammageozeta}(ii)
holds for many fields~$ \k -d $.
Ordered by~$ d $ the first five are
$ \k -15 $, $ \k -30 $, $ \k -35 $, $ \k -39 $ and~$ \k -42 $.
For the first, third and fifth of those~$ \psi_k(\tbgeo) $ generates $ \Kindtf k $
as~$ |K_2(\O)| = 2 $.

\noindent{}(ii)
The example discussed in~\S\ref{minus5example} below shows that one cannot ignore
the condition on the  stabilisers of the elements of~$ \Sigma_2^* $ in Theorem \ref{mthm}(ii)
Corollary~\ref{gammageozeta}(ii). Specifically, in this case the stabilisers of the elements of~$ \Sigma_3^* $ have order~2
or~3, one element of $ \Sigma_2^* $ has stabiliser of order~4, but $ \bgeo $ generates~$ \BB k {} $
so cannot be divided by~2.
\end{remark}

\subsection{A cyclotomic description of \texorpdfstring{$\bgeo$}{beta geo}}

Let $ k $ be an imaginary quadratic field of conductor~$N$.
Fixing an injection $k \to \CC $, the image is in the cyclotomic field $ F = \Q(\zeta_N) $ for~$\zeta_N := e^{2 \pi i/N}$.

The following result shows that the image under the induced map $ \BB k {} $  $\longrightarrow$ $\BB F {} $ of the element $ \bgeo $ constructed in Theorem \ref{mthm}(i) has a simple description in terms of elements constructed directly from roots of unity. (This result is, however, of very limited practical use  since it is generally much more difficult to compute explicitly in $ \Bbar F $ rather than in $ \Bbar k $).

\begin{proposition} The image of $ \bgeo$ in $ \BB F {} $ is equal to~$ N \sum_{\ol{a} \in \Gal(F/k) } [\zeta_N^a] $.
\end{proposition}

\begin{proof} At the outset we note that there is a commutative diagram
\begin{equation*}
\begin{tikzcd}
\Kindtf k \ar{r}{\simeq}\ar{d} & \BB k {} \ar{r}{\psi_k} \ar{d} & \Kindtf k \ar{d} \ar{r}{\reg_2} & \R(1) \ar[equal]{d}
\\
\Kindtf F \ar{r}{\simeq} & \BB F {} \ar{r}{\psi_F} & \Kindtf F  \ar{r}{\reg_2} & \R(1)
\,,
\end{tikzcd}
\end{equation*}
where the isomorphisms are obtained from~\eqref{tork3}, as well
as Theorems~\ref{ppmodcBbar} and~\ref{Bthm}, and $ \reg_2 $
is the regulator map corresponding to our chosen embeddings of
$ k $ and $ F $ into $ \CC $.

We further recall (from, for example, \cite[Prop.~5.13]{Sri96} with $ X = Y' = \Spec(F) $ and $ Y = \Spec(k) $) that the composition $ K_3(F) \to K_3(k) \to K_3(F) $ of the norm and pullback is given by the trace, and that the same is also true for the induced maps on $ \Kindtf F = K_3(F)_\tf $ and $ \Kindtf k = K_3(k)_\tf $.
By applying this fact to the element of $ \Kindtf F $
corresponding to $  N [\zeta_N] $ in $ \BB F {} $,
we deduce from the left hand square in the above diagram that there exists an element~$ \bcyc $ in $ \BB k {} $
that maps to $ N \sum_{\ol{a} \in \Gal(F/k) } [\zeta_N^a] $ in~$ \BB F {} $.

We now identify $ \Gal(F/k) $ with a subgroup of index 2 of $ (\Z/N\Z)^\times $,
which is the kernel of a primitive character $ \chi : (\Z/N\Z)^\times \to \CC^\times $
of order~2, corresponding to $ k $ (so $ \chi(\ol{-1}) = -1 $).
Then from the above diagram and Theorem~\ref{regthm}(iii), one finds that
$  (2 \pi i)^{-1}\cdot\reg_2(\psi_k(\bcyc)) $ is equal to
\begin{alignat*}{1}
\frac N{4\pi i} \sum_{\ol a \in \Gal(F/k)} \sum_{n \ge 1} \frac{ \zeta_N^{na} - \zeta_N^{-na} }{ n^2}
& =
 \frac N{4 \pi i} \sum_{\ol a \in \Gal(F/\Q)} \sum_{n \ge 1} \chi(\ol{a}) \frac{ \zeta_N^{na} }{ n^2}
\\
& =
\frac {N^{3/2}}{4 \pi } L(\Q, \chi, 2)
\\
& =
- 12  \zeta_k'(-1)
\end{alignat*}
as $ \zeta_k(s) = \zeta_\Q(s) L(\Q,\chi,s) $,
 with the Gau\ss\ sum $ \sum_{\ol a \in \Gal(F/\Q) } \chi(\ol a) \zeta_N^a = i \sqrt{N} $
(see \cite[\S58]{Hecke}).
(Cf.~the more general (and involved) calculation of \cite[p. 421]{Zag91},
or the calculation in the proof of \cite[Th.~3.1]{BGdJ} with $ r=-1 $, $ \ell = 1 $ and $ \O = \Z $.)

According to Theorem~\ref{mthm}, one has $(2\pi i)^{-1}\cdot \reg_2(\psi_k(\bgeo)) = -12 \cdot\zeta_k'(-1) $ and so, by the injectivity of $ \reg_2 $ on $ \Kindtf k $ (cf. Corollary~\ref{Bcor}) one has  $\psi_k(\bgeo) = \psi_k(\bcyc)$.
It then follows that $ \bgeo = \bcyc $ because~$ \psi_k $ is injective.
\end{proof}

\section{The proof of Theorem \ref{mthm}} \label{proofs}

Throughout this section we fix an imaginary quadratic field $k$ with ring of integers $\CO$, as in
\S\ref{tessellations}. In \S\ref{theproof} we also use the embedding of $ k $ into $ \CC $
chosen there.

\subsection{A preliminary result concerning orbits} We start by proving a technical result that will play an important role in later arguments.

We set $ V = k^2 \setminus\{(0,0)\} $ and let $ \C$ denote either $\SL_2(\O) $ or $ \GL_2(\O) $.

\begin{lemma}\label{first tech}\
\begin{itemize}
\item[(i)] For $ v $ in $ V $, $ \O^\times $ acts on the orbit $ \C v $, and
the natural map $ V \to \P_k^1 $ induces an injection of
$ \C v / \O^\times $ into $ \P_k^1 $, compatible with the action of~$ \C $.

\item[(ii)] For $ v_1 $ and $ v_2 $ in $ V $, the images of $ \C v_1 / \O^\times $
and $ \C v_2 / \O^\times $ are either disjoint or coincide.
\end{itemize}
\end{lemma}

\begin{proof} That $ \O^\times $ acts on $ \C v $ is clear if $ k \ne \k -1 $ or
$ \k -3 $ as $ \O^\times = \{\pm1\} $ and $ \C $ contains $ \pm \id_2 $.
For the two remaining cases, $ \O $ is Euclidean, and based on
iterated division with remainder in $ \O $ it is easy to find $ g $ in
$ \SL_2(\O) $ with $ g v = \begin{psmallmatrix} c\\0 \end{psmallmatrix} $
for some $ c $ in $ k^\times $, so if $ u $ is in $ \O^\times $, then $ uv $
is in the orbit of $ v $ because
$ g^{-1} \begin{psmallmatrix} u&0\\0& u^{-1} \end{psmallmatrix} g $
maps $ v $ to $ uv $.
Alternatively, this follows immediately from \cite[Chap.~7, Lem.~2.1]{EGMbook}
because if $ v = (\a, \b) $, and $ u $ is in $ \O^\times $, then $ (\a, \b) = (u \a, u \b) $.

Now assume that $ g_1 v = c g_2 v $ with $ c $ in $ k^\times $ and
the $ g_i $ in $ \C $. Then $ v $ is an eigenvector with eigenvalue~$ c $
of the element $ g_2^{-1} g_1 $ in $ \C $, which has determinant
in~$ \O^\times $.
Hence $ c $ is in $ \O^\times $, and $ g_1 v $ and $ g_2 v $ give the same element
in $ \C v / \O^\times $. So we do get the claimed
injection, and it is clearly compatible with the action of~$ \C $.

For the last part, suppose $ g_1 v_1 = c g_2 v_2 $ for some $ c $
in $ k^\times $, $ v_i $ in $ V $, and $ g_i $ in $ \C $.
Then~$ \C v_1 = c \C v_2 $ and the result is clear.
\end{proof}

\begin{proposition} \label{liftvectors}
If $ h $ is the class number of $ k $, then we can find $ v_1,\dots,v_h $
in $ V $ such that $ \P_k^1 $ is the disjoint union of the images
of the $ \C v_i / \O^\times $.
In particular, every element
in $ \P_k^1 $ lifts uniquely to some~$ \C v_i / \O^\times $ and this
lifting is compatible with the action of~$ \C $.
\end{proposition}

\begin{proof} By \cite[Chap.~7, Lem.~2.1]{EGMbook} we may identify $ \C \backslash V $ with
the set of fractional ideals of $\CO$ and hence $~\C \backslash V / k^\times = \C \backslash \P_k^1 $
with the ideal class group of $ k $.
We can then apply Lemma \ref{first tech}.
\end{proof}

\subsection{The proof of Theorem~\ref{mthm}} \label{theproof}

\subsubsection{}We first establish some convenient notation and conventions.

For a $2$-cell, or more generally, any flat polytope with vertices $ v_1,\dots,v_n $ in that order along
its boundary, we indicate an orientation by $ [v_1,\dots,v_n] $
up to cyclic rotation. The inverse orientation corresponds to
reversing the order of the vertices.
If we want to denote the face with either orientation, we write $ (v_1,\dots,v_n) $.
In particular, an orientated triangle is the same as a
3-tuple $ [v_1,v_2,v_3] $ of vertices up to the action (with
sign) of $ S_3 $.
Similarly, an orientated tetrahedron is the same as a 4-tuple
$ [v_1,v_2,v_3,v_4] $ up to the action (with sign) of $ S_4 $.
Recall that we defined maps $ f_{3,k} $ and $ f_{2,k} $ just
after~\eqref{messdef}. As mentioned above,
the map $ f_{3,k} $ is compatible with the action of $ S_4 $
by Remark~\ref{classicalcr} and~\eqref{2-2-term}.
By the properties of $ \crr_2 $ mentioned just before~\eqref{messdef},
the map $ f_{2,k} $ is also compatible with the action of $ S_3 $ on orientated triangles if
we lift them to elements of $ \pts 3 {\LL} $ for some suitable~$ \LL  $.

\subsubsection{}

By our discussion before the statement of the theorem,
the uniqueness of elements $ \b_\Q $ and $ \tilde \b_\Q $ with the stated properties is clear.
It is also clear that for any element $ \b_\Q $ in $ \Bbar k $  the explicit sum $ \tbgeo $ belongs to $ \pp(k) $.
In addition, the uniqueness of $ \b_\Q $ combines with the explicit expression for $\tbgeo$ to imply that $ 2 \tilde \b_\Q = \b_\Q $,
hence that $ 2 \tbgeo = \bgeo $.

The fact that $ \bgeo $ is independent of the subdivision~\eqref{subdiv},
and of the choice of representatives in $ \Sigma_3^* $ also follows directly from the equality in Corollary \ref{gammageozeta}(i) and the injectivity assertions in Corollary \ref{Bcor},
once it is known that $ \bgeo $ is in $ \BB k {} $.

Hence to prove Theorem \ref{mthm}, it suffices to prove the existence of
$ \b_\Q $ and $ \tilde \b_\Q $ in the stated groups such that the sums $ \bgeo $ and $ \tbgeo $ belong to $ \BB k {} $, and to do this we shall use the tessellation.

This argument is given in the next subsection. The basic idea is that, for $ k$ different from $\k -1 $ and $ \k -2 $, the sum $12\cdot |\C_{P_i}|^{-1}\cdot  \sum_{j = 1}^{n_i} [r_{i,j}] $ has integer coefficients and belongs to the kernel of $ \ddR {2,k} {} $ since the faces
of the polytopes $ P_i $ with those multiplicities can be matched
under the action of~$ \C $. This argument uses that $ f_{2,k} $
is invariant under the action of~$ \GL_2(k) $ and behaves compatibly with respect to permutations, just as $ f_{3,k} $.

The precise argument is slightly more complicated because the subdivision~\eqref{subdiv} induces triangulations
of the faces of the $ P_i $ which may not correspond, necessitating
the introduction of `flat tetrahedra', which give rise to the
term~$ \b_\Q $. Also, the
faces themselves may have orientation reversing elements in their
stabilisers.
But the resulting matching of faces does imply that the
explicit sum~$ \bgeo $ lies in the kernel $ \BB k {} $ of $ \ddR {2,k} {} $.

For the special cases $ k = \k -1 $, $k = \k -2 $, and $k= \k -3 $, we have to compute more explicitly
for the single polytope involved in each case.

\subsubsection{}We note first that each polytope $ P $ in the tessellation of $ \Hy $
comes with an orientation corresponding to it having
positive volume. For a face (2-cell) $ F $ in the tessellation,
we fix an orientation, and consider the group
$ \oplus_F \Z [F] $, where we identify $ [F^\dagger] $ with $ - [F]$ if
$ F^\dagger $ denotes $ F $ with the opposite orientation.

To $ P $ we associate its boundary $ \partial P $ in
this group, where each face has the induced orientation.
As the action of $ \C $ on $ \Hy $ preserves the orientation,
it commutes with the boundary map.

We now need to do some counting. For a face $ [F] $, we let $ \C_F $
denote the stabiliser of the (non-oriented) face~$ F $, and $ \C_F^+ $
the subgroup that preserves the orientation $ [F] $.
 We note that the index of $\C_F^+$ in $\C_F$ is either $1$ or $2$.

Let $ P $ and $ P'$ be the polytopes in the tessellation that
have $ F $ in their boundaries. If $ g $ is in~$ \C_F $ then $ gP = P $ or $ P' $,
and $ gP = P $ precisely when $ g $ is in $ \C_F^+ $. Therefore
$ \C_F^+ = \C_F \cap \C_P $.

It is convenient to distinguish between the following two cases for the $ \C $-orbits of $ F $.

\begin{itemize}
\item
$ \C_F = \C_F^+ $.
If $ P $ and $ P' $ in $ \Sigma_3^* $ are such that their
boundaries each contain an element in the $ \C $-orbit of $ [F] $,
then $ P $ and $ P' $ are in the same $ \C $-orbit, hence are the same.
Therefore there is exactly one $ P $ in $ \Sigma_3^* $ that contains faces in the $ \C $-orbit
of $ [F] $.  If two faces of $ P $ are in the $ \C $-orbit
of $ [F] $, then they are transformed into each other already
by $ \C_P $. Hence the number of elements in the $ \C $-orbit
of $ F $ in $ \partial P $ is $ [\C_P : \C_F] = [\C_P : \C_F^+] $.
If $ P' $ is the element in $ \Sigma_3^* $ that has an element
in the $ \C $-orbit of $ [F^\dagger] $ in its boundary (with $ P = P' $
and $ P \ne P' $ both possible), then there are $ [\C_{P'} : \C_F] = [\C_{P'} : \C_F^+] $
elements in the $ \C $-orbit of $ [F^\dagger] $ in the boundary of $ P' $.

\item
$ \C_F \ne \C_F^+ $. Note that in this case $ [\C_F : \C_F^+] = 2 $.
Here $ [F] $ and $ [F^\dagger] $ are in the same $ \C $-orbit
and as above one sees that there is only one element $ P $ of
$ \Sigma_3^* $ that has elements in this $ \C $-orbit in its
boundary.  Any two such elements can be transformed into each
other using elements of $ \C_P $, so there are $ [\C_P : \C_F^+] $
of those in the boundary of $ P $.
\end{itemize}

\subsubsection{}In this subsection we prove Theorem \ref{mthm}
for $ k $ not equal to $ \k -1 $,  $\k -2 $ or $ \k -3 $.

In this case Corollary~\ref{cor:lcm} implies that the order of each group $\C_{P_i}$ divides~12, and so both the formal sum of elements in $ \Sigma_3^* $ given by
\begin{equation*}
 \pi_P = \sum_{i = 1}^m  \frac{12}{|\C_{P_i}|}  [P_i]
 \,,
\end{equation*}
 and the formal sum of of tetrahedra resulting from the subdivision~\eqref{subdiv}
\begin{equation*}
 \pi_T = \sum_{i = 1}^m  \frac{12}{|\C_{P_i}|}  \sum_{j = 1}^{n_i} [T_{i,j}]
 \,,
\end{equation*}
have integral coefficients.

We extend the boundary
map $ \partial $ to such formal sums, where for $ \pi_T $ the
boundary is a formal sum of ideal triangles, contained in the
original faces of the polytopes $ P_i $.
(Note that the subdivision~\eqref{subdiv} may introduce `internal
faces' inside each polytopes, but by construction the parts of
the boundaries of the tetrahedra here cancel exactly.
This also holds after lifting all vertices to $ \O^2 $ as there
is no group action involved in order to match them. So we may,
and shall, ignore those internal faces.)

The subdivision~\eqref{subdiv} induces a triangulation of each
face of each $ P_i $ in $ \Sigma_3^* $. We let $ [\D_F] $
denote the induced triangulation.
But if $ F $ is a face of such a $ P_i $ with $ [F] \ne [F^\dagger] $, then the
induced triangulations $ [\D_F] $ and $ [\D_{F^\dagger}] $ (which may come
from a different element in $ \Sigma_3^* $) may not match.
Similarly, if $ gF $ and $ F $ are both faces of $ P_i $ with $ g $ in $ \C_{P_i} $, then
$ g [\D_F] $ and $ [\D_{gF}] $ may not match.

A typical example of non-matching triangulations is that
of a `square' $ F = [v_1,v_2,v_3,v_4] $ that is cut into
two triangles using either diagonal, resulting in the triangulations
 $ [v_1,v_2,v_4] + [v_2,v_3,v_4] $ and
 $ [v_1,v_2,v_3] + [v_1,v_3,v_4] $.
 But the boundary of
the orientated tetrahedron $ [v_1,v_2,v_3,v_4] $ gives exactly the former
minus the latter.
Using induction on the number of vertices of a face $ F $
it is easily seen that any two triangulations of $ [F] $ with
the same orientation differ by the boundary of a formal sum of 
tetrahedra contained in $ F $. Such tetrahedra have
no volume, and the cross-ratio of its four cusps is in $ \Q^\flat $.
We refer to them as `flat tetrahedra', and if $ [\D_F^1] $ and
$ [\D_F^2] $ are two triangulations (with the same orientation)
 of an orientated face~$ [F] $, we shall write
\[ `[\D_F^1] \equiv [\D_F^2] \,\,\text{ mod. } \partial(\text{flat tetrahedra})\text{'}\]
if $ [\D_F^1] - [\D_F^2] $ is the boundary of a formal sum of such
flat tetrahedra.

In particular, if $ [\D_F] $ and $ [\D_{F^\dagger}] $  are any triangulations
of the faces $ [F] $ and $ [F^\dagger] $ (so with opposite orientation),
then $ [\D_F] + [\D_{F^\dagger}] \equiv 0 $ mod. \ftt.

We extend the boundary map to the free abelian group $ \oplus_{P \in \Sigma_3^*} \Z[P] $
by linearity. For a given~$ [F] $, in $ \partial (\pi_P) $ we find
find $12\cdot |\C_{F}^+|^{-1}$ copies of $ [F] $ (up to the action
of $ \C $).
If $ [F] \ne [F^\dagger] $ (i.e., if $ \C_F = \C_F^+ $) then
this equals the number of copies of $ [F^\dagger] $, and we
combine~$ [F] $ and~$ [F^\dagger] $.

We consider four cases, based on the exponents
of~2 in $|\C_F|$ and $|\C_F^+|$.
Note that $ \C_F^+ $ is cyclic, so by Lemma~\ref{finite-order-lemma}
and our assumptions on $ k $ we can write its order as $ 2^s m $
with $ m = 1 $ or~3, and~$ s = 0 $ or 1, with the case $ m= 3 $
and $ s = 1 $ not occuring. Then $ |\C_F| = 2^t |\C_F^+| $ with
$ t = 0 $ or~1.

\begin{itemize}
\item[(1)] $ s = t = 0 $. Here $ F^\dagger $ is not in the same $ \C $-orbit as $ F $, and in
$ \partial(\pi_P) $ the contribution of their $ \C $-orbits is
$ \frac{12}m [F] + \frac{12}m [F^\dagger] $, modulo the action of $ \C $.
Then for $ \partial(\pi_T) $ they contribute $ \frac{12}m [\D_F] + \frac{12}m [\D_{F}^\dagger] $
mod. \ftt and modulo the action of~$ \C $.

\item[(2)] $ s = 0 $ and $ t = 1 $. Here $ F $ and $ F^\dagger $ are in the same~$ \C $-orbit, and in the boundary of~$ \pi_P $, up to the action of $ \C $, we have
$ \frac{12}m [F] = \frac6m [F] + \frac6m [F^\dagger] $. Then in $ \partial(\pi_T) $ we obtain
$ \frac6m [\D_F] + \frac6m [\D_F^\dagger] $ mod.~\ftt and modulo the
action of~$ \C $.

\item[(3)] $s = 1$ and $  t = 0 $. This is similar to case~(1), but now $ [F] $ and $ [F^\dagger] $ both occur with
coefficient~ $ 6 $ in~$ \partial(\pi_P) $ because $ m = 1 $.
In~$  \partial(\pi_T) $ we obtain
$ 6 [\D_F] + 6 [\D_F^\dagger] $ mod. \ftt, and modulo the
action of~$ \C $.

\item[(4)] $ s = t = 1 $. This is similar to case~(2), but now in~$  \partial(\pi_P) $ we find
$ 6 [F] = 3 [F] + 3 [F^\dagger] $, again because $ m = 1 $, hence in
$ \partial(\pi_T) $ this gives
$ 3 [\D_F] + 3 [\D_F^\dagger] $ mod. \ftt and modulo the
action of~$ \C $.
\end{itemize}

We see that there exists some $ \a $, a formal sum of flat tetrahedra,
such that $ \pi_T + \a $ has boundary, up to the action of $ \C $,
a formal sum with terms $ [t] + [t^\dagger] $ with $ t $ an ideal triangle.
Lifting all cusps to $ \LL = \C v_1/\O^\times \coprod \dots \coprod \C v_h/\O^\times $
as in Proposition~\ref{liftvectors}, and applying $ f_{3,k} $ as
in~\eqref{newCD}, we see from the $ \C $-equivariance of $ f_{2,k} $, and the fact that this
map is alternating, so kills elements of the form $ [t] + [t^\dagger] $,
that $ \sum_{i=1}^m 12\cdot |\C_{P_i}|^{-1}\cdot\sum_{j=1}^{n_i} [r_{i,j}] + \b' $
is in the kernel of $ \ddR {2,k} {\O^\times} $, where $ \b' $ is the
image of~$ \a $.

Note that $ \b' $ lies in the image of $ \Bbar {\Q} $
in $ \Bbar k $. Multiplying by $ |\O^\times| = 2 $ and setting $ \b_\Q := 2 \b' $ we complete the proof of Theorem \ref{mthm}(i) in this case.

The proof of Theorem \ref{mthm}(ii) is done in the same way, starting with
$ \sum_{i=1}^m 6\cdot |\C_{P_i}|^{-1}\cdot \sum_{j=1}^{n_i} [r_{i,j}] $ (which has integer coefficients under the stated assumptions).
In this case the coefficients in the above cases (1), (2) and~(3) are divided by~2, and case~(4) is ruled out by the assumptions.

\subsubsection{} \label{specialfields}
We now consider the special fields $ \k -1 $, $ \k -2 $, and $ \k -3 $.

In each of these cases either $|\C_{P_i}|$ does not divide~12 or $ |\O^\times| $ is
larger than~2.
However, one also knows that $ \Sigma_3^* $ has only one element and its stabiliser
has order 12 or~24 and so the result of Theorem \ref{mthm}(ii) does not apply. It is therefore enough to prove Theorem \ref{mthm}(i) for these fields.

If $k = \k -1 $, then $ \Sigma_3^* $ is an octahedron, with stabiliser isomorphic to~$ S_4 $.
Using that an ideal tetrahedron with positive volume in this octahedron
must contain exactly two antipodal points, it is easy to see that the
subdivision is unique up to the action of the stabiliser.
Hence the resulting element under $ f_{3,k} $ is well-defined.
Computing it explicitly as $ \sum_{j=1}^4 [r_{1,j}]  = 4 [\w] $
one finds that it is in the kernel of $ \ddR {2,k} {} $ as $ \w^2 = -1 $,
so we can simply take $ \b_\Q = 0 $.

If $k= \k -3 $, then the polytope is a tetrahedron, with stabiliser
isomorphic to~$ A_4 $. Computing its image $ [r_{1,1}] $ under
$ f_{3,k} $ explicitly one finds
$ [\w] $ with $ \w^2 = \w-1 $,
and $ \ddR {2,k} {} ([r_{1,1}]) = \w \tw (1-\w) = \twe -1 -1 $ in $ \twt k^\times $, which has order~2 by Remark~\ref{twotorremark}.
So we can again take~$ \b_\Q = 0 $.

If $k = \k -2 $, then the polytope is a rectified cube (i.e., a cuboctahedron), with stabiliser isomorphic to~$ S_4 $,
so it has six 4-gons and eight triangles as faces.
By the commutativity of~\eqref{newCD}, for $ \tL=\{1\} $, we can compute
$ \ddR {2,k} {} (\sum_{j}[r_{1,j}]) $ by choosing lifts of all
vertices involved, and applying~$ f_{2,k} $ to each of the lifted
triangles (with correct orientation) of the induced triangulation
of the faces of $ P_1 $.
(This provides an alternative approach for $ \k -1 $ and $ \k -3 $ as well.)
Note that any triangulation of the faces occurs for some subdivision:
fix a vertex $ V $ and use the cones on all the triangles that
do not have~$ V $ as a vertex.

So we must consider all triangulations.
Giving the six 4-gons one of the two possible triangulations at random resulted in
$ \twe -1 -1 = (-1) \tw 2 =  \ddR {2,k} {} ([2]) $ in $ \twt k^\times $.
The other triangulations we obtain from this one by triangulating
one or more of the 4-gons differently.
For each 4-gon, this adds
the image under $ \ddR {2,k} {} $ of the cross-ratio
of the corresponding flat tetrahedron.
As the 4-gons are equivalent under $ \C $,
one easily computes this equals $  \ddR {2,k} {} ([2]) $.
So by Proposition~\ref{BbarQprop}
we find that $ \ddR {2,k} {} (\sum_j [r_{1,j}])$ equals either~$ 0 $
or~$\ddR {2,k} {} ([2]) $, and both occur.
By Corollary~\ref{Bcor}, we must take $ \b_\Q = 0 $ or~$ [2] $,
and by Proposition~\ref{BbarQprop} the latter is not in $ 2 \Bbar {\Q} $.

This completes the proof of Theorem \ref{mthm}.

\subsubsection{}We make several observations concerning the above argument.

\begin{remark} \label{tworemark}
Let $ k $ not be equal to $ \k -1 $, $ \k -2 $ or $ \k -3 $.
 One can try to find a better element than $ \bgeo $ by going
through the calculations in the proof of Theorem~\ref{mthm}
after replacing $ \pi_P $ by an element of the form $\sum_{i=1}^m M\cdot |\C_{P_i}|^{-1} [P_i]$ for some positive integer $ M $ that is divisible by
the orders of the stabilisers $ \Gamma_{P_i} $. If we start with
$ M $ equal to the least common multiple of the orders $|\C_{P_i}|$, then we may
have to multiply this element by 2 perhaps twice in the proof in order to ensure that the resulting element in $ \Bbar k $ belongs to $ \BB k {} $:
\begin{enumerate}[label=(\arabic*)]
\item\label{it:tworemark-1}
in order to ensure that the boundary $ \partial $ of the resulting
analogue of~$  \pi_T $ is trivial up to the action of~$ \C $, which is not
automatic if some $ P_i $ has a face with reversible orientation
under~$ \C $ and  $M\cdot |\C_{P_i}|^{-1}$ is odd;

\item\label{it:tworemark-2}
in order to ensure that $ \sum_{i=1}^m M\cdot |\C_{P_i}|^{-1}\cdot \sum_{j=1}^{n_i} [r_{i,j}] + \b' $
is in the kernel of $ \ddR {2,k} {} $ and not just $ \ddR {2,k} {\O^\times} $, where $ M $
results from~\ref{it:tworemark-1}, and $ \b' $ (coming from flat tetrahedra)
is in $ \Bbar {\Q} $, which we view as inside $ \Bbar k $ by Proposition~\ref{BbarQprop}.
\end{enumerate}
Note that in the second part here we use $ |\O^\times| = 2 $,
which excludes $ k = \k -1 $ or $ \k -3 $.

For our $ k $, with $ \Nm \colon k \to \Q $ the norm,
the Hermitian form
$ (x,y) \mapsto \Nm(x)+\Nm(y) + \Nm(x-y) $ on $ \CC^2 $
has minimal vectors $ \{\pm(1,0), \pm(0,1), \pm (1,1)\} $.
By \cite[Th.~2.7]{imquad-coh}, this means that the triangle
with vertices $0$, $1$ and $ \infty $ is a $2$-cell of the tessellation.
The element $ \begin{psmallmatrix} 0 & -1\\ 1 & -1 \end{psmallmatrix} $
in $ \C $ of order~3  stabilises this triangle while preserving its orientation. Therefore
the 3-cells that share this triangle as faces have stabilisers
with orders divisible by~3, and $ M $ is divisible by~3.

Also, the elements $ \begin{psmallmatrix*}[r] -1 & 0\\ 0 & 1 \end{psmallmatrix*} $
and $ \begin{psmallmatrix*} 0 & 1\\ 1 & 0 \end{psmallmatrix*} $
have order~2, and generate a subgroup of $ \C $ of order~$4$.
The first has as axis of rotation the $1$-cell connecting $0$ and
$ \infty $, so the axis of rotation of the second, which meets
this 1-cell, must meet either a $3$-cell, or a $2$-cell with vertices
$0$, $ \infty $, and purely imaginary numbers. In the first
case we start with $ M $ divisible by~$6$. In the second case,~\ref{it:tworemark-2} above
ensures that $ M $ is even since the $2$-cell reverses
orientation under
$ \begin{psmallmatrix*} 0 & 1\\ 1 & 0 \end{psmallmatrix*} $.
Because in this remark we are also assuming that $ k \ne \k -2 $,
we know from Corollary~\ref{cor:lcm} that the greatest common divisor of the orders of the $ \C_{P_i} $ divides~$12$.
So this method could lead to an element $ \bgeo $ as in Theorem~\ref{mthm}(i)
but with $24$ replaced by either~$6$, $12$, or $24$.
\end{remark}

\begin{remark} \label{sadremark} In our calculations, we find the following for all fields $ k$ that differ from $\k -1 $, $ \k -2 $
and $ \k -3 $:
\begin{itemize}
\item[$\bullet$]
 the $ \gcd $ of the orders of the stabilisers $ \C_{P_i} $ is~6 or~12, i.e., 3 does not occur;

\item[$\bullet$]
 the sum $ \sum_{i = 1}^m 12\cdot |\C_{P_i}|^{-1}\cdot\sum_{j = 1}^{n_i} ( [r_{i,j}] - [\ol{r_{i,j}}] ) $ belongs to $ \BB k {} $,
so that by Remark~\ref{mpropnew} and Corollary~\ref{Bcor}(iii),
this must be another expression for $ \bgeo $ in $ \BB k {} $.
\end{itemize}
Unfortunately, we have not been able to prove either of these statements
in general.
\end{remark}

\subsection{An explicit example}\label{minus5example}

With the same notation as before Theorem~\ref{mthm}, we consider the element $ \bb = \sum_{i = 1}^m M\cdot |\C_{P_i}|^{-1}\cdot \sum_{j = 1}^{n_i} [r_{i,j}] $ of~$ \Bbar k $, where $ M$ is the greatest common divisor of the orders $|\C_{P_i}|$. As
in Remark~\ref{tworemark}, the proof of Theorem~\ref{mthm}
shows that there is
a positive divisor~$ e $ of $ 2 |\O^\times| $
such that $ \ddR {2,k} {} ( e \bb ) $ is in the image
of the composition $ \Bbar {\Q} \to \twt \Q^\times \to \twt k^\times $,
and it gives a way of computing~$ e $
from the tessellation. But this depends on choices,
for example, on
how one pairs the faces of the $ P_i $ under the action of $ \C $,
so the $ e $ found may not be~optimal.

But one can also do this algebraically,
by computing $ \bb $ and determining a (minimal) positive integer~$ e $
with~$ e \ddR {2,k} {} (\bb) $ in the image of~$ \Bbar {\Q} $.
For this we can use Corollary~\ref{twcomputecor} and Remark~\ref{twcomputerem}:
if $ S $ is a finite set of finite places of~$ k $ such that $ \twt \O_S^\times \subset \twt k^\times $
contains~$ \delta  = \ddR {2,k} {} (\b) $, and $ A = \Q^\times \cap \O_S^\times $, then one
can compute if $ e \delta $ is in the image of $ \twt A $ or
not.
If this is the case then one can find its preimage in $ \twt \Q^\times $ from
Lemma~\ref{kertwmap}, algorithmically determine if an element
in there gives the trivial element in $ K_2(\Q) $, and if so,
express it in terms of $ \ddR {2,\Q} {} ([x]) $ with~$ x $~in~$ \Q^\flat $.

Note that a different subdivision~\eqref{subdiv}
might a priori give rise to a different $ \bb $ and a different~$ e $,
but the other choices are irrelevant in this algebraic approach.

For the reader's convenience we illustrate this, and the methods of the proof of Theorem~\ref{mthm},
in the special case that $ k = \k -5 $. In particular, in this case
we find that both methods give the same element, which generates
$ \BB k {} $.

The lifts to $ \O^2 $ (up to scaling by $ \O^\times $) of the vertices
$ v_1,\dots,v_8 $ in the two elements $ P_1 $ and $ P_2 $ of $ \Sigma_3^* $
are the columns of the matrix
\begin{equation} \label{exlifts}
\left(
\col{\w+1\\-2}\\\col{1\\0}\\\col{2\\\w-1}\\\col{2\\\w}\\\col{0\\1}\\\col{-\w\\2}\\\col{1\\-1}\\\col{-\w+1\\\w+1}
\right)
\,.
\end{equation}
Both polytopes are triangular prisms, which we write as
$ [a,b,c; A,B,C] $, for $ [a,b,c] $ and $ [A,B,C] $
triangles, with $ A $ above $ a $, etc. Such a prism can be subdivided into
orientated tetrahedra as
$ [a,A,B,C] - [a,b,B,C] + [a,b,c, C] $,
and the resulting subdivision of its orientated boundary is
\begin{equation} \label{psubbound}
   [A,B,C]
 + [a,A,C] - [a,c,C]
 + [a,b,B] - [a,A,B]
 + [b,c,C] - [b,B,C]
 - [a,b,c]
\,.
\end{equation}
Here the first and last terms correspond to the triangular faces,
and the middle terms are grouped as pairs of triangles in the
rectangular faces.

Then $ P_1 = [v_3, v_5, v_4;v_1, v_2, v_6] $ with $ \C_{P_1} $
of order~2, generated by $ g_1 = (v_1 v_3) (v_2 v_4) (v_5 v_6) $
in cycle notation on the vertices of $ P_1 $ ($ v_7 $ and $ v_8 $
are mapped elsewhere).
It interchanges the orientated faces $ [v_1,v_3,v_4,v_6] $ and $ [v_3,v_1,v_2,v_5] $ of
$ P_1 $, and those faces have trivial stabilisers. The orientated face
$ [v_2,v_6,v_4,v_5] $ is mapped to itself by~$ g_1 $
but its stabiliser is non-cyclic of order~4, with one of the two orientation
reversing elements acting as $ h_1 = (v_2 v_5)(v_4 v_6) $.
The two triangles $ [v_1,v_6,v_2] $ and $ [v_3,v_5,v_4] $ are interchanged
by $ g_1 $, and both have stabilisers of order~2, with the one for $ [v_1,v_6,v_2] $
generated by~$ (v_2 v_6) $.

We have $ P_2 = [v_1, v_8, v_3; v_2, v_7, v_5] $ with $ \C_{P_2} $
of order~3, generated by $ g_2 = (v_1 v_3 v_8)(v_2 v_5 v_7) $.
The three orientated faces
$ [v_2,v_7,v_8,v_1] $,
$ [v_7,v_5,v_3,v_8] $ and
$ [v_5,v_2,v_1,v_3] $
are all in the same $ \C $-orbit and have trivial
stabilisers.
The two triangles $ [v_1, v_8, v_3] $ and $ [v_2,v_5,v_7] $
are necessarily non-conjugate (even ignoring orientation) as $ \C_{P_2} $ has order~3,
but both have (orientation non-preserving) stabiliser of order~6, which acts as the full permutation
group on their vertices.

We now subdivide both $ P_1 $ and $ P_2 $ as stated just before~\eqref{psubbound}.
This gives an element
\begin{align*} \pi_T' :=\, &\frac 12 \pi_T  = 3 [v_3,v_1, v_2, v_6] - 3 [v_3, v_5, v_2, v_6] + 3 [v_3, v_5, v_4, v_6]\\
                       &\hskip 1.5truein  + 2 [v_1,v_2, v_7, v_5] - 2 [v_1, v_8, v_7, v_5] + 2 [v_1, v_8, v_3, v_5]
\,.\end{align*}

Applying $ f_{3,k} $ to this element results in
\begin{equation*}
\bb =
7[\tfrac 13 \w + \tfrac23] - 3 [- \tfrac 23 \w + \tfrac53] + 3 [ \tfrac16 \w + \tfrac56] - 2 [- \tfrac 16 \w + \tfrac 76]
\end{equation*}
in $ \Bbar k $.
Using~\eqref{psubbound} one can easily compute the boundary $ \partial \pi_T' $.
Because
\begin{equation*}
\Sigma_2^* = \set{ (v_1,v_3,v_4,v_6) , (v_2,v_6,v_4,v_5), (v_3,v_4,v_5) , (v_1,v_8,v_3) , (v_7,v_2,v_5) }
\,,
\end{equation*}
under the action of $ \C $ we can move the resulting triangles to the eight triangles that
result from the elements of $ \Sigma_2^* $, with one of the `squares'
giving rise to four inequivalent triangles due to the two ways
of triangulating a `square', the other to only one inequivalent
triangle.

In $ (\oplus_t \Z[t])_\C $, where $ t $ runs through
the triangles,
the triangular faces in $ \partial(\pi_T') $ coming from those
of $ P_1 $ and $ P_2 $ cancel under the action
of $ \C $ (this uses that the coefficient of $ [P_2] $ in $ \pi_T' $
is even and the triangular faces of $ P_2 $ have orientation reversing elements
in their stabilisers).
Of course, the `internal' triangles created
by the subdivision into tetrahedra always cancel.
Using~$ g_1, g_2, h_1, h_2 $ one moves the triangles coming
from the `square' faces to the five inequivalent triangles coming
from $ [v_1,v_3,v_4,v_6] $ and $ [v_2,v_6,v_4,v_5] $.
This yields the sum of the elements
\begin{align*}
\partial{[v_4,v_3,v_1,v_6]} = \, & 3 ( [v_3,v_1,v_6] - [v_3,v_4,v_6] + [v_1,v_6,v_4] - [v_1,v_3,v_4] )\\
& \hskip 0.3truein + 2 ( [v_3,v_4,v_6] - [v_3,v_1,v_6] + 2( [v_1,v_3,v_4] - [v_1,v_6,v_4] ) )\end{align*}
and
\begin{equation*}
3 [v_5,v_4,v_6] - 3 [v_5,v_2,v_6] = [v_5,v_4,v_6] - [v_5,v_2,v_6] - \partial [v_5,v_2,v_4,v_6]
\,,
\end{equation*}
where we used $ h_1 $.

So for $ \pi_T = 2 \pi_T'= 6 [P_1] + 4 [P_2] $ we get
$ \partial(\pi_T - 2 \partial [v_4,v_3,v_1,v_6] + 3 \partial [v_5,v_2,v_4,v_6] ) = 0 $
modulo the action of~$ \C $. After multiplying by~2 in order to deal with the ambiguity of the lifts
of the cusps to $ \O^2 $, we then find the element
\[ \bgeo = 4 \bb - 4 [3] + 6 [4/5]\in \BB k {} ,\]
with the last two terms arise because $ \crr_3([v_4,v_3,v_1,v_6]) = 3 $
and $ \crr_3([v_5,v_2,v_4,v_6]) = 4/5$.

To see if one could do better, as discussed just before this
example, we, instead, compute $ \ddR {2,k} {} $ of $ \bb $.
This can be done easily using the matching of triangles under the action
of $ \C $ as before, using the commutativity of~\eqref{newCD} for $ \ddR {2,k} {\O^\times} $,
but as for this we have to lift the vertices to the
column vectors in $ \O^2 $ in~\eqref{exlifts}
(and not up to scaling by $ \O^\times $) we pick up some additional
torsion along the way.
Alternatively, one can choose a finite set of finite primes $ S $
for $ k $ such that, for every $ [z] $ occurring in~$ \bb $, both
$ z $ and $ 1-z $ are $ S $-units, and compute in $ \twt k^\times $ as in
Remark~\ref{twcomputerem} and Corollary~\ref{twcomputecor}.
The result is
\begin{equation*}
 - ( -4 ) \tw  5 + ( -2 ) \tw 3 + ( -1 ) \tw 2 + ( -1 ) \tw ( -\w ) - 2 \tw ( -5 )
\end{equation*}
in $ \twt k^\times $.
Note that $  [v_5,v_4,v_6] - [v_5,v_2,v_6] $ under $ f_{2,k} $ is
mapped to $ \twe -\tfrac12 -\frac{\w}2 - 2 \tw (-\w) = \twe -1 {-\w} - 2 \tw (-5) $.
The first three terms are in the image of $ \Bbar {\Q} $, and
if we multiply the last by 2 then we obtain $ - 4 \tw (-5) = - (-4) \tw (-5) = - (-4) \tw 5 + (-1) \tw (-1) $,
with the first again in the image of~$ \Bbar {\Q} $.
If $ \twe -1 -1 $ would come from $ \Bbar {\Q} $ then it would
come from its torsion by Proposition~\ref{BbarQprop}(ii), which
is generated by $ [2] $.
By Lemma~\ref{kertwmap}, the kernel of $ \twt {\Q^\times} \to \twt k^\times $ has
order~2, and is generated by $ \twe -1 -5 $, and one easily checks
using Corollary~\ref{twcomputecor}
that both $ \twe -1 -1 $ and $ \twe -1 -1 + \twe -1 -5 $ in $ \twt \Q^\times $
are neither trivial nor equal to $ \ddR {2,\Q} {} ([2]) = (-1) \tw 2 $.
Therefore $ \twe -1 -1 $ in $ \twt k^\times $ is not in the image of~$ \Bbar {\Q} $.
Multiplying by~2 again kills the term $ \twe -1 -1 $,
hence $ 4 \b - 4 [3] + 6 [5] $ is in~$ \BB k {} $
is the best possible for our choice of subdivision.
(Note this element equals $ \bgeo $ above as $ [\tfrac45] = - [\tfrac15] = [5] $.
Also note that these calculations also show that $ 2 \bb - 2 \ol{\bb} $ is in $ \BB k {} $,
in line with Remark~\ref{sadremark}, and that this element must
also equal $ \bgeo $.)

In fact, $ K_2(\O) $ is trivial by \cite[\S7]{BeGa},
so by Corollary~\ref{gammageozeta}(i),  $ \psi_k(\bgeo) $ is a generator of
the infinite cyclic group $ \Kind k = \Kindtf k $, $ \bgeo $ is a generator
of the infinite cyclic group $ \BB k {} $, and
the map $ \psi_k \colon \BB k {} \to \Kindtf k = \Kind k $
is an isomorphism. 

Hence a different triangulation cannot give a better result.
Also, both factors 2 in Remark~\ref{tworemark} are necessary, so 
the obstruction of `incompatible lifts' under the action of $ \C $ is non-trivial.

\section{Finding a generator of \texorpdfstring{$K_3(k)$}{K3(k)}, and computing \texorpdfstring{$|K_2(\CO_k)|$}{the order of K2(Ok)}} \label{generators}

In this section, we again restrict to the case of an imaginary quadratic field $k$ and set $\O := \O_k$.
As in Section~\ref{tessellations}, we fix an embedding $\sigma: k \to \CC $ and regard $k$ as a subfield of $\CC$.

We explain how to combine
 an implementation of an algorithm of Tate's, which produces a natural number that is known to be divisible by the order of $ K_2(\O) $, with either the result of Corollary \ref{gammageozeta}(i) or just the known validity of the precise form of Lichtenbaum's conjecture for~$ k $ and $ m = 2 $ (cf.~Example~\ref{klichtenbaum}), to deduce our main computational results.

At this stage we have successfully applied the first of these approaches to about twenty fields and the second to hundreds of fields.
In this way, for example, we have, for all imaginary quadratic fields of discriminant bigger
than $ -1000$,
determined the order of $K_2(\CO)$, where not yet known, and a generator of the infinite cyclic group $K_3(k)_\tf^\ind$ that lies in the image of the injective homomorphism $\psi_k$ constructed in Theorem \ref{psiprop} (thereby verifying that $k$ validates Conjecture \ref{psi is iso}) and hence also the Beilinson regulator value $R_2(k)$.

The results are available online~\cite{database}.
In particular, for each of the listed imaginary quadratic number
fields $ k $, the element $ \balg $ is such that its image $ \psi_k(\balg) $
generates $ \Kindtf k $, thus verifying Conjecture~\ref{psi is iso}
for all those fields. The element $ \bgeo $ is the element of
Theorem~\ref{mthm}, obtained in the way described in Remark~\ref{sadremark}.

\subsection{Dividing \texorpdfstring{$ \bgeo $}{betageo} by  \texorpdfstring{$ |K_2(\O_k)| $}{the order of K2}}\label{first app section}

\subsubsection{}The basic approach is as follows.
An implementation by Belabas and Gangl \cite{BeGa} of (a refinement of) an algorithm of Tate gives an explicit natural number ~$ M $
divisible by $|K_2(\O)|$.
Since typically $M = |K_2(\O)|$,
for any element $ \ageo $ of $ \ker(\delta_{2,k}) $
that lifts $ \bgeo $, we try to find an element $ \a $ in this kernel
for which the difference $ M\cdot \a  - \ageo $ lies in the subgroup generated by~\eqref{5-term} and~\eqref{2-2-term}.
If one finds such an~$\a$, then its class $ \b $ in $\overline{B} (k) $
satisfies~$ \bgeo = M\cdot \b $.
From the result of Corollary~\ref{gammageozeta}(i) it then follows that
$ |K_2(\O)| = M $, that $ \b $ generates~$\overline{B} (k) $, that $\psi_k(\b) $ generates $ \Kindtf k $ and hence, by Theorem~\ref{regthm}(iii),
that~$ R_2(k) = |\reg_2(\psi_k(\b))| $.

\subsubsection{}To find a candidate element $ \a $ as above we first use the methods described in \S\ref{S unit section} below to identify an element $ \a$ for which one can verify \emph{numerically} that $ M \cdot\tD_\s(\a) = \tD_\s(\ageo) $, with~$\tD_\s$ the homomorphism from Remark \ref{factorremark}(ii).
We then aim to prove \emph{algebraically} that $ M\cdot \b = \bgeo $
by writing the difference $ M\cdot \a - \ageo $ as a sum of explicit relations of the form ~\eqref{5-term}
and~\eqref{2-2-term}.

To complete this last step we use a strategy that can be used to investigate whether any element of the form $ \sum_i n_i [x_i] $, where the $ n_i $ are integers and
the $ x_i $ are in~$ k^\flat $, can be written as a sum of such relations,
using suitable finite subsets $ U $ of $ k^\flat $.

We let $ U $ consist of all $ x_i $ and their images under
the 6-fold symmetry implied by the relations~\eqref{2-2-term},
i.e., for $ u $ in $ U $ we also adjoin
$ 1-u $, $ u^{-1} $, $ 1-u^{-1} $, $ (1-u^{-1})^{-1} = \frac{-u}{1-u} $
and~$ (1-u)^{-1} $.

Next, for $ u \ne v $ in $ U $, we consider the element in $ \Z[k^\flat] $
obtained by putting $ x = u $ and~$ y = v $ in~\eqref{5-term}.
We use the result only if all five terms are in~$ U $.

We then form a matrix $ A $ of width $ |U| $, as follows.
\begin{itemize}
\item
For the first row we write $\sum_i n_i [x_i] $ in terms of the
$ \Z $-basis $ \{[u] \text{ with } u \text{ in } U \} $
of the subgroup $ \Z[U] $ of $ \Z[k^\flat] $.

\item
For each of the, $n$ say, 5-term relations that we have just generated, we add
a row writing it in terms of the basis.

\item
For each $ u $ in $ U $ we add rows corresponding to the relations
$ [u] + [1-u] $, $ [u] +  [u^{-1}] $, $ [u] - [1-u^{-1}] $, $ [u] + [\frac{-u}{1-u}] $
and $ [u] - [(1-u)^{-1}] $, resulting in, say,~$ m $ rows in
total.
\end{itemize}

Then the kernel of the right-multiplication by $ A $ on $ \Z^{1+n+m} $
(as row vectors)
gives the relations among the various elements that we put into
the rows of $ A $.
An element in this kernel with $1$ as its first entry
encodes a rewriting of~$ \sum _i n_i [x_i] $ as sum
of elements as in~\eqref{5-term} and \eqref{2-2-term}.

Unfortunately, this straightforward method is rarely
successful. Instead, we may have to enlarge~$ U $,
and the computation can simply become too large.
It was, however, done successfully, to some extent by trial and error, for a
several imaginary quadratic number fields.

\begin{example}
The most notable example among those
is the field  $k = \mathbb{Q}(\sqrt{-303})$, for which it is
known from \cite{BeGa} that~$ |K_2(\O)| = 22$. The results for this case are described in Appendix~\ref{huge}.
\end{example}

\begin{remark} We note that the method described above for verifying identities in $ \overline{B}(k) $ only depends on the definition of $ \overline{B}(k) $ in terms of the boundary
map $ \delta_{2,k} $ and the relations~\eqref{5-term} and~\eqref{2-2-term}
on $ \Z[k^\flat] $ that are used to define $ \ol{\pp}(k) $.
In particular, it does not rely on knowing the validity of Lichtenbaum's Conjecture and so, in principle, the same approach could be used to show that an element
is trivial in $ \Bbar F {} $ for any number field $ F $ (although, in practice, the computations are likely to quickly become unfeasibly large).
\end{remark}

\subsection{Finding a generator of \texorpdfstring{$ \Kindtf k $}{K3} directly}\label{second app section}

This approach relies on the effective bounds on $ |K_2(\O)| $ that are discussed above, the known validity of Lichtenbaum's Conjecture
as in Example~\ref{klichtenbaum}, and an implementation of the `exceptional $S$-unit' algorithm  (see \S\ref{S unit section} below) that produces  elements in $\overline{B}(k)$.
In particular, the reliance on Lichtenbaum's Conjecture means that the general applicability of this sort of approach is currently restricted to abelian fields.

To describe the basic idea, we assume to be given an element $\c $ that equals $N_\c$
times a generator of the (infinite cyclic) group $ \Kindtf k $
for some non-negative integer~$ N_\c $.

Then one has $ |\reg_2(\c)| = N_\c \cdot R_2(k) $ and so Example~\ref{klichtenbaum} implies that
\[ - \frac{\reg_2(\c)}{12\cdot \zeta_k'(-1)} = \frac{N_\c}{|K_2(\O)|} .\]

If one also has an explicit natural number ~$ M $ that is known to be divisible by $ |K_2(\O)| $, then
\begin{equation}\label{numerical} - M \frac{\reg_2(\c)}{12\cdot \zeta_k'(-1)} = N_\c  \frac{M}{|K_2(\O)|}
\end{equation}
is a product of a non-negative and a positive integer.
Hence, if the left hand side of this equality is numerically close to a natural number $d_\c$, then~$ N_\c$ and $M/|K_2(\O)|$ are both divisors of $d_\c$.
In particular, if one
has an element $\c$
with $d_\c = 1$, then one concludes both that $ N_\c = 1 $ (so that $\c $ is a generator of $ \Kindtf k $) and that $ |K_2(\O)| = M $.
We would therefore have identified a generator of $ \Kindtf k $ and determined~$ |K_2(\O)| $.

To find suitable elements $ \c $ we proceed as
described in \S\ref{S unit section} below to generate elements~$ \a $ in the subgroup $ \ker(\delta_{2,k})$ of $\Z[k^\flat]$. We then
let~$ \c$
be the image under $\psi_k$ of the image of~$ \a $ in~$\BB k {} $.

Note here that it is not {\em a priori} guaranteed that $\im(\psi_k)$ contains a generator of $ \Kindtf k $. However, if this is the case (as it was in all of the examples we tested), then $ \psi_k $ is surjective and so one has verified that $k$ validates Conjecture \ref{psi is iso}.

\subsection{Constructing elements in \texorpdfstring{$\ker (\delta_{2,k})$}{ker(delta2k)}  via exceptional \texorpdfstring{$S$}{S}-units}
\label{S unit section}

\subsubsection{}In order to find enough elements in $ \ker(\delta_{2,k}) $ we
fix a finite set~$ S $ of finite places of $ k $, and
consider `exceptional $S$-units', where
an $ S $-unit $ x $ is exceptional if $1-x$ is also an $S$-unit.

To compute with such elements it is convenient to fix a basis of the $S$-units
of $ k $, i.e., a set of $ S $-units that give a $ \Z $-basis of the $ S $-units
modulo torsion.
(This is implemented in GP/PARI \cite{PARI} as `bnfsunits'.)
For each exceptional $S$-unit $ x $ we encode $ x $ and $ 1-x $
using the exponents that arise
when they are expressed in terms of the
basis and a suitable root of unity.
Corollary~\ref{twcomputecor} and Remark~\ref{twcomputerem} then
enable us to compute effectively in $\twt k^\times $ with the
elements~$ (1-x) \tw x $.

\begin{example}
   In the case $k=\Q(\sqrt{-11})$ and $S=\{\wp_2, \wp_3,\overline{\wp}_3\}$ where $\wp_2=(2)$ is the unique prime ideal
of norm~4 and $\wp_3$ and $\overline{\wp}_3$ denote the two prime ideals of norm~3 in $\OO$,
PARI provides the $S$-unit basis $\mathcal{B} = \{b_1,b_2,b_3\} $ with $b_1=2$, $b_2= \frac{-1+\sqrt{-11}}2$ and
 $b_3= \frac{-1-\sqrt{-11}}2$. We find the exceptional $S$-unit  $x= \frac 5{36} - \frac{\sqrt{-11}}{36}$
 of norm $\frac {1}{36}$, for which $1-x$ has norm $\frac{3}{4}$, and write
 $$x = - b_1^{-1} b_2^{-2}\,, \qquad 1-x = - b_1^{-1} b_2^{-2} b_3^{3}\,.$$
It follows that
\[  (1-x) \tw x = (-1) \tw (-1) +(-1) \tw b_1 + (-1) \tw b_3 + 3 (b_1 \tw b_3) + 6 (b_2 \tw b_3),\]
which corresponds to the element $(\ol{1}, \ol{1}, \ol{0},\ol{1},0,3,6)$ under the isomorphism in Corollary~\ref{twcomputecor}.
\end{example}

This approach effectively reduces the problem of finding elements in $\ker (\delta_{2,k})$  to a concrete problem in linear algebra.
Of course, one wants to choose a finite set $S$ of finite places for which one can find sufficiently many exceptional
$S$-units in $k$ such that some linear combination of them in $\ker (\delta_{2,k})$ gives a \emph{non-trivial} element in
(and preferably a generator of) the quotient ~$\overline{B}(k)$.

Note that, whilst one can check for non-triviality of an element
$ \b $ in $ \BB k {} $ by simply verifying that its image under the map $ \tD $ is numerically non-trivial, in order to conclude that~$ \b $ is trivial, we need to know an explicit natural number $ M $ that is divisible by $ |K_2(\O)| $. Then the quantity on the left hand side of (\ref{numerical}) is numerically close to zero if and only if $ \c = \psi_k(\b) $ is trivial.
If that is the case, then the injectivity of $\psi_k$ implies that the element $ \b $ is itself trivial.

\subsubsection{}

Since
$ \overline B (k) $ is cyclic of infinite order, there exists a finite
set $ S $ for which the above procedure can lead to a non-trivial
element.
By Remark~\ref{mpropnew} one can take it to comprise
all of the places that divide any of the principal
ideals~$\OO\cdot r_{i,j}$, $\OO\cdot \ol{ r_{i,j}} $, $ \O \cdot (1-r_{i,j}) $ and~$ \O \cdot (1-\ol{ r_{i,j} } ) $
for the elements~$r_{i,j}$ in Theorem~\ref{mthm}.
In general, this set
is far too large to be practical for the exceptional $S$-unit approach. Fortunately, however, in all of the cases investigated in this paper we
found that a much smaller set suffices. In fact,
it is often enough to take $S$ to comprise all places that divide either $2$ or $3$ or any of the first ten (say) primes that split in $k$.

\begin{example}
In the case $k=\Q(\sqrt{-303})$ it suffices to take $S$ to be the set of places that divide either of $2, 11$ and $13$ (all of which split in $k$) or $3$ (which ramifies in $k$). Imposing small bounds on the
exponents with respect to a chosen basis, we already find $ 683$ exceptional $S$-units in $k$.
Setting $\w := (1+\sqrt{-303})/2$, GP/PARI's
\cite{PARI} `bnfsunits' gives as a basis of the $ S $-units the set
\[\{ -20 + 3 \w , 2 , -4-\w , -36 - \w , 4 - \w , 28 - \w , -12 + \w \}\]
of norms $2^{10}$, $2^2$, $2^5\cdot 3$, $2^7\cdot 11$, $2^3\cdot 11$, $2^5\cdot 13$,
and $2^4\cdot 13$, respectively.

Computing the kernel of $ \delta_{2,k} $ on the corresponding subgroup of $ \Z[k^\flat] $
gives a free $ \Z $-module of rank several
hundreds, but most of the elements of a $ \Z $-basis for this kernel turn out to result
in the trivial element of~$\Bbar {k} $, i.e., correspond to relations
of the type~\eqref{5-term} and~\eqref{2-2-term}.
In this case, with a set of exceptional units that differed from
the one used in Appendix~\ref{huge}, but again using that $ |K_2(\O)| = 22 $,
we found using the approach described in \S\ref{second app section}
that

{\tiny
\begin{dmath*}
-46\crrterm{\frac{-\w-27}{64}}+36\crrterm{\frac{-\w+15}{16}}-14\crrterm{\frac{-\w+41}{16}}-48\crrterm{\frac{-\w+8}{4}}-62\crrterm{\frac{-\w+41}{4}}+18\crrterm{\frac{-\w+41}{48}}+34\crrterm{\frac{-11}{2}}+42\crrterm{\frac{-143}{1}}-16\crrterm{\frac{-253\w+2321}{6144}}+28\crrterm{\frac{-253\w-495}{3328}}+158\crrterm{\frac{-3\w+23}{32}}+4\crrterm{\frac{-39\w-221}{512}}+120\crrterm{\frac{-5\w+73}{64}}+18\crrterm{\frac{-5\w+49}{416}}+44\crrterm{\frac{\w+91}{64}}+70\crrterm{\frac{\w+27}{64}}-12\crrterm{\frac{\w+91}{128}}-36\crrterm{\frac{\w+1}{16}}+82\crrterm{\frac{\w-2}{2}}+92\crrterm{\frac{\w+7}{32}}-36\crrterm{\frac{\w+40}{4}}-22\crrterm{\frac{\w-8}{4}}-116\crrterm{\frac{11}{2}}+58\crrterm{\frac{11}{8}}-34\crrterm{\frac{13}{2}}-84\crrterm{\frac{13}{4}}+34\crrterm{\frac{3\w+20}{8}}+14\crrterm{\frac{9\w-17}{352}}
\end{dmath*}
}%
\noindent{}in
$ \Z[k^\flat]$ belongs to $ \ker(\delta_{2,k}) $, and that its image in $\overline{B}(k)$ is sent by $\psi_k$ to a generator of $ \Kindtf k $.
\end{example}

\begin{example} We now set $ k = \k -4547 $, so that $ \O = \Z[\w] $ with $ \w = \frac{1+\sqrt{-4547}}2 $.
We recall that in \cite{BrG} it was conjectured that $|K_2 (\O)| \, = 233 $.
In fact, while the program developed in loc.~cit.~showed that $|K_2 (\O)|$
divides $233$, the authors were unable to verify their conjecture since this would have required them
to work in a cyclotomic extension of too high a~degree.

By using the approach described in \S\ref{second app section}, we were now able to verify that $|K_2 (\O)|$ is indeed equal to~$233$ and, in addition, that the  element
{\tiny
\begin{dmath*}
132\crrterm{\frac{-2\w+5}{117}}-2\crrterm{\frac{-1}{12}}+8\crrterm{\frac{-2\w-3}{1404}}-6\crrterm{\frac{-2\w+5}{18}}-14\crrterm{\frac{-2\w+1752}{19683}}-2\crrterm{\frac{-1}{2}}+8\crrterm{\frac{-2\w-3}{27}}+2\crrterm{\frac{-1}{288}}+2\crrterm{\frac{-1}{3}}-24\crrterm{\frac{-2\w+1752}{3159}}-2\crrterm{\frac{-1}{36}}+128\crrterm{\frac{-2\w+5}{39}}+24\crrterm{\frac{-2\w-3}{4212}}+74\crrterm{\frac{-2\w+5}{52}}+12\crrterm{\frac{-2\w+5}{54}}-12\crrterm{\frac{-2\w-840}{6591}}-54\crrterm{\frac{-\w+421}{351}}+40\crrterm{\frac{-2\w-3}{72}}-16\crrterm{\frac{-2\w-3}{78}}+12\crrterm{\frac{-2\w+5}{8}}-10\crrterm{\frac{-\w+421}{468}}+2\crrterm{\frac{-13}{16}}+4\crrterm{\frac{-13}{18}}-6\crrterm{\frac{-13}{24}}+2\crrterm{\frac{-13}{243}}-2\crrterm{\frac{-13}{3}}-2\crrterm{\frac{-13}{4}}-8\crrterm{\frac{-169}{324}}-2\crrterm{\frac{-4\w-1680}{177957}}-6\crrterm{\frac{-208}{81}}-2\crrterm{-26}-4\crrterm{\frac{-26}{3}}+2\crrterm{\frac{-27}{4}}+42\crrterm{\frac{-3\w+1263}{2197}}-12\crrterm{\frac{-31\w+162}{13182}}+22\crrterm{\frac{-31\w-131}{2197}}-4\crrterm{\frac{-5\w+6}{64}}-8\crrterm{\frac{-16\w+6736}{2197}}+2\crrterm{\frac{-16\w-4523}{2197}}-26\crrterm{\frac{2\w-5}{1404}}+2\crrterm{\frac{2\w+3}{18}}+2\crrterm{\frac{1}{18}}-4\crrterm{\frac{\w+875}{1053}}-4\crrterm{\frac{\w+875}{1296}}-24\crrterm{\frac{2\w-5}{27}}+12\crrterm{\frac{2\w+1750}{3159}}-2\crrterm{\frac{1}{32}}+14\crrterm{\frac{2\w-5}{351}}-14\crrterm{\frac{2\w-5}{4212}}-50\crrterm{\frac{2\w+3}{52}}-78\crrterm{\frac{2\w+3}{54}}+4\crrterm{\frac{2\w-842}{6591}}+38\crrterm{\frac{\w+420}{351}}-30\crrterm{\frac{2\w-5}{72}}-2\crrterm{\frac{2\w-5}{78}}-14\crrterm{\frac{2\w+3}{8}}+14\crrterm{\frac{\w+420}{4056}}-6\crrterm{\frac{\w+420}{468}}-6\crrterm{117}-2\crrterm{\frac{13}{256}}-2\crrterm{\frac{13}{81}}-4\crrterm{\frac{169}{16}}-14\crrterm{\frac{169}{243}}+4\crrterm{\frac{169}{256}}-10\crrterm{\frac{3494\w-13298}{2197}}-16\crrterm{\frac{4\w-1684}{177957}}-42\crrterm{\frac{16\w+4523}{8788}}+2\crrterm{\frac{243}{256}}+2\crrterm{\frac{26}{27}}-2\crrterm{\frac{26}{9}}+4\crrterm{\frac{3}{32}}-16\crrterm{\frac{3\w+1260}{2197}}+10\crrterm{\frac{31\w+131}{13182}}-24\crrterm{\frac{31\w-162}{2197}}-6\crrterm{\frac{39}{2}}+6\crrterm{\frac{39}{8}}+8\crrterm{\frac{8\w+3360}{351}}+20\crrterm{\frac{8\w+3360}{9477}}-30\crrterm{\frac{10\w+2}{1053}}-38\crrterm{\frac{5\w+1}{54}}+6\crrterm{\frac{5\w+1}{64}}-12\crrterm{\frac{5\w-6}{78}}-4\crrterm{\frac{5\w+1}{1053}}+24\crrterm{\frac{5\w+1}{27}}+8\crrterm{\frac{10\w-12}{729}}-6\crrterm{52}-4\crrterm{\frac{52}{81}}-2\crrterm{\frac{64}{81}}-14\crrterm{\frac{16\w+4523}{4563}}+4\crrterm{\frac{841\w-176104}{177957}}
\end{dmath*}
}
\noindent{}
of $ \Z[k^\flat]$ belongs to $ \ker(\delta_{2,k}) $, and that its image in $\overline{B}(k)$ is sent by $\psi_k$ to a generator of $ \Kindtf k $.
\end{example}

\appendix
\section{Orders of finite subgroups}\label{orders}
In this appendix we again consider a fixed imaginary quadratic
field $ k $, embedded into $ \CC $.
The main aim of this subsection is to prove, in Corollary~\ref{cor:lcm}, that the lowest common multiple of the orders of  finite subgroups of $\PGL_2(\CO)$ is
either 12 or 24, with the latter being the case only
if $ k = \k -1 $ or $ \k -2 $.
The authors are not aware of a suitable reference for this
in the literature, or in fact, of an explicit classification of
types of finite subgroups of $ \PGL_2(\O) $ that do not lie in
 $ \PSL_2(\O) $.  As this is not difficult, we include it for the sake of
completeness.

Using the inclusions $ \PSL_2(\O) \subset \PGL_2(\O) \subset \PGL_2(\CC) = \PSL_2(\CC) $,
 our arguments are based on the following classical result \cite[Chap.~2, Th.~1.6]{EGMbook}
that goes back to Klein~\cite{binaereFormenMathAnn9}.

\begin{proposition} \label{kleinprop}
A finite subgroup of $ \PSL_2(\CC) $ is isomorphic to
a cyclic group of order~$ m \ge 1 $, a dihedral group of order $ 2m $ with $ m \ge 2 $, $ A_4 $, $S_4 $ or $A_5$. Further, all of these possibilities occur. \end{proposition}

It seems that the finite subgroups of $ \PSL_2(\O) $ have been studied more than those of $ \PGL_2(\O) $ even though it is harder to
determine them.

An element $ \c $ in $ \SL_2(\O) $ of
finite order has characteristic polynomial of the form $ x^2+ax+1 $
with $ a $ in $ [-2,2] \cap \O $ as the two roots must be conjugate
roots of unity. Hence $ a = \pm 2 $, $ \pm 1 $ or~0, from
which it follows readily that the image $ \ol{\c} $ in $ \PSL_2(\O) $
of $ \c $ has order 1, 2~or~3.
In view of Proposition~\ref{kleinprop}, this limits the possibilities of a finite subgroup of $ \PSL_2(\O) $
to the cyclic groups of orders 1, 2,~or~3,
the dihedral groups of order 4 or 6,
and~$ A_4 $.
Cyclic groups of order 2~or~3 can be obtained already
in the subgroup $ \PSL_2(\Z) $, generated by
$ \begin{psmallmatrix} 0&-1\\1& 0\hfil \end{psmallmatrix} $
and $ \begin{psmallmatrix} 0&-1\\ \hfil1&-1 \end{psmallmatrix} $
respectively.
The occurrence of the dihedral groups of order~4 and~6, and of~$ A_4 $,
in $ \PSL_2(\O) $ for $ k = \k -d $ with $ d $ a positive square-free integer,
depends on the prime factorisation of $ d $; see~\cite[Satz~6.8.]{Kr15}.

\smallskip

We now consider finite subgroups of $ \PGL_2(\O) $.
An element in $ \PGL_2(\O) $ of odd order is contained in $ \PSL_2(\O) $:
if an odd power of its determinant is a square in $ \O^\times $, then so is its
determinant.
It follows that the only finite subgroups that can occur
in $ \PGL_2(\O) $ but are not necessarily contained in $ \PSL_2(\O) $
are the cyclic groups of even order, the dihedral groups, and~$ S_4 $.
In order to obtain a complete answer, we first look at elements of finite
order in~$ \PGL_2(k) $.

\begin{lemma} \label{finite-order-lemma}
Let $ \c $ be an element
of finite order in $ \GL_2(k) $, and $ \ol{\c} $ its image
in $ \PGL_2(k) $. Write ${\rm ord}(\c)$ and ${\rm ord}(\ol{\c})$ for the respective orders of these elements. 
\begin{itemize}
\item[(i)] For $ k = \k -1 $, one has  ${\rm ord}(\c) \in \{1,2,3,4,6,8,12\}$ and ${\rm ord}(\ol{\c}) \in \{1,2,3,4\}$.

\item[(ii)] For $ k = \k -2 $, one has ${\rm ord}(\c) \in \{1,2,3,4,6, 8\}$ and ${\rm ord}(\ol{\c}) \in \{1,2,3, 4\}$.

\item[(iii)] For $ k = \k -3 $, one has ${\rm ord}(\c) \in \{1,2,3,4,6, 12\}$ and ${\rm ord}(\ol{\c}) \in \{1, 2, 3,6\}$.

\item[(iv)] For all other $ k $, one has ${\rm ord}(\c ) \in \{1,2,3,4,6\}$ and ${\rm ord}(\ol{\c}) \in \{1,2,3\}$.
\end{itemize}
\end{lemma}

\begin{proof} Assume $ \c $ has order~$ n $, and let $ a(x) $ in $ k[x] $ be its minimal polynomial,
so that $ a(x) $ divides $ x^n-1 $.
The statement is clear if $ a(x) $ splits into linear factors,
so we may assume~$ a(x) $ is irreducible in $ k[x] $ and of degree~2.
If $ a(x) $ is in $ \Q[x] $ then it is an irreducible factor
in~$ \Q[x] $ of $ x^n-1 $, necessarily the $ n $th cyclotomic polynomial
as the $ m $th cyclotomic polynomial divides~$ x^m-1 $ if $ m $
divides~$ n $.
So $ \phi(n) = 2 $, and $ n = 3 $, 4 or~6.
If $ a(x) $ is not in $ \Q[x] $ then $ a(x) \ol{a(x)} $
is irreducible in $ \Q[x] $, it must be the $ n $th cyclotomic
polynomial, and $ k $ must be a subfield of~$ \Q(\zeta_n) $.
Then~$ \phi(n) = 4 $, so $ n = 8 $ and $ k = \k -1 $
or $ \k -2 $, or $ n = 12 $ and $ k = \k -1 $
or $ \k -3 $. (Note that $ n = 5 $ is excluded as $ \Q(\zeta_5) $
contains no imaginary quadratic field.)
The statement about the order of $ \ol{\c} $ follows by taking
into account the factorisation of $ x^n-1 $ over $ k[x] $.
In general, the $ 2m $th cyclotomic polynomial divides $ x^m+1 $.
But for $ k = \k -1 $ and $ n = 12 $, so $ m=6 $, we also have $ x^6+1 = (x^3-i)(x^3+i) $
in $ k[x] $, and $ a(x) $ divides one of those factors.
\end{proof}

\begin{remark}
In the cell stabiliser calculation for $ k = \k -3 $ in~\cite{Mendoza}
the symbol $A_4$ should be a dihedral group of order~12
in $\PGL_2(\O)$, where the subgroup in $\PSL_2(\O)$ is dihedral
of order~6.
\end{remark}

We can now determine the types of finite subgroups in $ \PGL_2(\O) $
that do not lie in $ \PSL_2(\O) $.

\begin{proposition} \label{HGLnotSL} Let $ G $ be a finite subgroup of $ \PGL_2(\O) $ that is not contained in $ \PSL_2(\O) $.
\begin{itemize}
\item[(i)] For $ k $ not equal to $ \k -m $ with $ m=1 $,~$ 2$~or~$ 3 $,
$ G $ is isomorphic to a cyclic group of order~$ 2 $, or a dihedral group
of order~$ 4 $ or~$ 6 $.
For $ k = \k -1 $ and $  \k -2 $, $ G $ can also be isomorphic
to a cyclic group of order~$ 4 $, a dihedral group of order~$ 8 $, or $ S_4 $,
and for~$ k = \k -3 $ to a cyclic group of order~$ 6 $ or a dihedral group of order~$ 12 $.

\item[(ii)] All the groups listed in claim (i) occur.
\end{itemize}
\end{proposition}

\begin{proof}
We already observed that an element of finite odd order in $ \PGL_2(\O) $
is contained in~$ \PSL_2(\O) $, so the groups listed in~(i) are those
not ruled out by combining Proposition~\ref{kleinprop}
with Lemma~\ref{finite-order-lemma}.

It remains to show that all such groups occur. Various examples
may of course exist in the literature, but for the sake of completeness
we give some here.
In fact, for $ S_4 $ we use the
stabiliser of the single element in $ \Sigma_3^* $
in our calculations for~$ \k -1 $ respectively~$ \k -2 $ (see
\S\ref{specialfields}).

\emph{Cyclic examples.}
If $ u $ in $ \O^\times $ is not a square, then
$ \begin{psmallmatrix} u&0\\ 0&1 \end{psmallmatrix} $
is not in $ \PSL_2(\O) $ and its order
equals the order of~$ u $. This gives the required subgroups
except for those of order~$ 2 $ for $ \k -1 $, and of order~$ 4 $ for $ \k -2 $.
The former can be obtained by using
$ \begin{psmallmatrix} 0&\rtmo\\ 1&0 \end{psmallmatrix} $,
which is not in $ \PSL_2(\O) $ and has order~$2$, and the latter
by using 
\smash{$ \begin{psmallmatrix} 0&1\\ 1&\sqrt{-2} \end{psmallmatrix} $,}
which is not in $ \PSL_2(\O) $ and has order~$4$.

\emph{Dihedral examples.}
If $ u $ in $ \O^\times $ is not a square, and has order~$ m = 2 $,
$4$~or~$6$, then
$ \begin{psmallmatrix} u&0\\ 0&1 \end{psmallmatrix} $
and~$ \begin{psmallmatrix} 0&1\\ 1&0 \end{psmallmatrix} $
generate a dihedral group of order~$ 2m $.
With $ u=-1 $ this constructs a copy of $ D_4 $, except for $ \k -1 $,
but for this field we can use generators
$ \begin{psmallmatrix} -1&0\\ \hfill 0&1 \end{psmallmatrix} $
and $ \begin{psmallmatrix} 0&1\\ \rtmo&0 \end{psmallmatrix} $.
Taking $ u $ a generator of $ \O^\times $
gives a copy of $ D_8 $ for $ \k -1 $, and a copy of $ D_{12} $ for~$ \k -3 $.

For $ k $ not equal to $ \k -1 $, a suitable copy of $ D_6 $
is generated by
$ \begin{psmallmatrix} 0&-1\\ 1&-1 \end{psmallmatrix} $
and $ \begin{psmallmatrix} 0& 1\\ 1&0 \end{psmallmatrix} $,
whereas for $ \k -1 $ we can use
$ \begin{psmallmatrix} 0&-1\\ 1&-1 \end{psmallmatrix} $
and $ \begin{psmallmatrix} 1& \rtmo\\ 1+\rtmo&-1 \end{psmallmatrix} $.
Finally, a copy of $ D_8 $ for $ \k -2 $ is generated by
\smash{$ \begin{psmallmatrix} 0&1\\ 1&\sqrt{-2} \end{psmallmatrix} $
and $ \begin{psmallmatrix} 0 & -1\\ 1& \hfill 0 \end{psmallmatrix} $.}\vskip .1cm

\emph{$ S_4 $-examples.}
For $ \k -1 $, the orders of
\mm[1,\rtmo+1,0,-1], \mm[\rtmo,-1,-1,0] and  \mm[1,\rtmo,0,-\rtmo] are~2, 3, and~4, respectively, and they generate a subgroup
isomorphic to $ S_4 $.

For $ \k -2 $, the orders of
\mm[-2,-\sqrt{-2}-1,-\sqrt{-2}+1,1], \mm[-\sqrt{-2}-1,-\sqrt{-2}+1,1,\sqrt{-2}] and \mm[2,\sqrt{-2}+1,\sqrt{-2}-1,-2]
are~3,~3, and 2, respectively, and they also generate a subgroup isomorphic to $ S_4 $.
\end{proof}

\begin{corollary} \label{cor:lcm}
The lowest common multiple of the orders of the finite subgroups of
$\PGL_2(\O)$ is ~$24$ if $ k $ is either $ \k -1 $ or $ \k -2 $ and is ~$12$ in all other cases.
\end{corollary}

\begin{proof}
This is true for the groups listed in Proposition~\ref{HGLnotSL},
and the possible finite subgroups of~$ \PSL_2(\O) $
(discussed before Lemma~\ref{finite-order-lemma}) have order dividing~12.
\end{proof}

\section{A generator of \texorpdfstring{$ K_3(k)_\tf^\ind $}{K3} for \texorpdfstring{$ k = \Q(\sqrt{-303}) $}{a particular field}}\label{huge}

For an imaginary quadratic field $ k = \Q(\sqrt{-d})$, with ring of integers $\O$, Browkin \cite{browkin}  has identified conditions under which the order $|K_2 (\O)|$ is divisible by either $2$ or $3$ (for example, he shows that $|K_2 (\O)|$ is divisible by $3$ if $d\equiv 3 \mod 9$).

Moreover, all of the coefficients in the linear combination $ \bgeo $ that occurs in Theorem \ref{mthm}(i) are
divisible by~$2$ if $ k $ is not equal to either $ \k -1 $ of $ \k -2 $.
In addition, if $ k $ is also not equal to $ \k -3 $, then whilst Remark~\ref{tworemark} shows that at least one of the coefficients
in $ \bgeo $ is not divisible by~$3$ one finds, in practice, that most of these coefficients are divisible by $3$.

For these reasons, it can be relatively easy to divide $ \bgeo $ by either $2$ or $3$.
But no such arguments work for division by  primes larger than $3$ and this
requires considerably more work.

It follows that if one uses the approach of  \S\ref{first app section}, then any attempt to obtain a solution $\b$ in~$ \Bbar k $ to the equation $ |K_2(\O)| \cdot\b = \bgeo $,
or equivalently (taking advantage of Remark \ref{mpropnew}) to the equation  $ 2|K_2(\O)| \cdot\b = 2\bgeo $,
in order to find
a generator of $ \BB k {} $, is likely to be much more difficult when $ |K_2(\O)| $ is divisible
by a prime larger than~$3$.

This observation motivates us to discuss the field $ k := \k -303 $,
for which~$\O = \Z[\w]$ with~$\w = (1 + \sqrt{-303})/2$.
We recall that $ k $ was conjectured in \cite{BrG}, and verified in~\cite{BeGa},
to be the imaginary quadratic field of largest
discriminant for which $|K_2(\O)|$ is divisible by a prime
larger than $3$. More precisely, this order was first conjectured and later determined
to equal~22.

We apply the technique described in \S\ref{tessellations}.
The quotient $\PGL_2(\O)\backslash \Hy$ has volume
\[\vol(\PGL_2(\O)\backslash \Hy)
= -\pi \cdot\zeta_k'(-1)
\approx 140.1729768601914879815382141215\dots \,.\]
The tessellation of $\Hy$ consists of $132$ distinct
$\PGL_2(\O)$-orbits of $3$-dimensional polytopes:
\begin{itemize}
\item $87$ tetrahedra 
\item $29$ square pyramids 
\item $13$ triangular prisms 
\item $1$ octahedron 
\item $2$ hexagonal caps---a polytope with a hexagonal base,
$4$ triangular faces, and $3$ quadrilateral faces as shown in Figure~\ref{fig:hexcap}.
\begin{figure}
\newdimen\Rad
\Rad=2cm
\newdimen\r
\r=0.35\Rad
\begin{tikzpicture}
  \draw (0:\Rad) \foreach \x in {60, 120, ...,300} {--(\x:\Rad)} -- cycle;
  \draw (90:\r) \foreach \x in {210, 330} {-- (\x:\r)} -- cycle;
  \foreach \x in {90, 210, 330}
  {
    \draw (\x-30:\Rad) -- (\x:\r)--(\x+30:\Rad);
  };
\end{tikzpicture}
\caption{Hexagonal cap}\label{fig:hexcap}
\end{figure}

\end{itemize}

The stabiliser $\Gamma_P $ in $ \PGL_2(\O) $ of each of these polytopes $ P $
is trivial except for eight polytopes~$ P $.  It
has order $2$ for four triangular prisms, and 
order $3$ for one triangular prism, the octahedron, and both
hexagonal caps.  By Theorem~\ref{mthm},
the tessellation and stabiliser data give rise to an explicit element
$\bgeo$, which we compute by using the `conjugation trick' of Remark \ref{mpropnew} as $\frac{1}{2} (2\bgeo)$ (see Remark~\ref{sadremark} for why we are allowed to divide by 2).
The latter can be written as a sum of 188 terms $  a_i [z_i]$, where $a_i $ is in $ 2 \ZZ$ and $z_i
$ in $ k^\flat$.
By Theorem~\ref{mthm} we have~$ \sum a_i \Di(z_i) = 24 \pi\cdot \vol(\PGL_2(\O)\backslash \Hy) $
with $\Di$ the Bloch-Wigner dilogarithm.

Using the algebraic approach described in \S~\ref{first app section}, we can find
$ \balg = \sum b_j [z_j]$ in $\BB k {} $ with image under the Bloch-Wigner function bounded away from~0,
and it turns out that it suffices to restrict the search to exceptional $S$-units
where $S$ consists of the prime ideals above $\{2,3, 11,13,19\}$.
Here $3$ ramifies in $ \O $, and the other primes are the first four primes that
split in~$\O$.
One of the~$\balg$ found with smallest positive dilogarithm value
has~$110$ terms and all coefficients~$\pm 2$.

By comparing $\sum b_j \Di(z_j)$ with $ \sum_i a_i D(z_i) $ above,
we expect $\bgeo - 22\cdot\balg$
to be trivial in~$ \BB k {} $.
We can prove this by writing a lift to $ \Z[k^\flat] $ explicitly as a sum of the
elements specified in~\eqref{5-term} and~\eqref{2-2-term}.
A linear algebra calculation in Magma \cite{MAGMA} shows that
this can be done
as an integral linear combination of $1648$
$5$-term relations,
plus a good number of $2$-term relations, 
so that, indeed, $ \bgeo = 22 \cdot\balg $.
(Note that $ \BB k {} $ by Corollary~\ref{Bcor}
injects into $ \R $ under the dilogarithm, whereas one has to
contend with torsion if attempting this calculation in $B(k)$ instead.)
The elements $\balg$ and  $\bgeo$ are given below.  The five-term combinations are
available online~\cite{database}.

It follows from Corollary~\ref{gammageozeta}(i) that the resulting element $ \psi_k(\balg) $
generates~$ \Kindtf k $.
This also implies that $ \psi_k \colon \BB k {} \to \Kindtf k $ is bijective (as predicted by Conjecture \ref{psi is iso}).

{\tiny
\begin{dmath*}
\balg=
-2\crrterm{\frac{-\w+41}{52}}-2\crrterm{\frac{-\w+2}{6}}-2\crrterm{\frac{-\w-1}{6}}-2\crrterm{\frac{-\w+92}{64}}-2\crrterm{\frac{-\w+701}{676}}-2\crrterm{\frac{-\w+12}{8}}-2\crrterm{\frac{-\w+4}{8}}-2\crrterm{\frac{-\w-3}{1}}-2\crrterm{\frac{-\w+8}{12}}-2\crrterm{\frac{-\w+12}{16}}-2\crrterm{\frac{-\w+15}{16}}-2\crrterm{\frac{-\w+2}{16}}-2\crrterm{\frac{-\w-11}{16}}-2\crrterm{\frac{-\w-25}{2}}-2\crrterm{\frac{-\w+15}{22}}-2\crrterm{\frac{-\w+26}{22}}-2\crrterm{\frac{-3\w+46}{66}}-2\crrterm{\frac{-\w-14}{22}}-2\crrterm{\frac{-\w+26}{24}}-2\crrterm{\frac{-\w+2}{26}}-2\crrterm{\frac{-\w+8}{32}}-2\crrterm{\frac{-\w+25}{32}}-2\crrterm{\frac{-\w+28}{32}}-2\crrterm{\frac{-\w-4}{32}}-2\crrterm{\frac{-\w+389}{352}}-2\crrterm{\frac{-\w-36}{352}}-2\crrterm{\frac{-\w+8}{4}}-2\crrterm{\frac{-\w-7}{4}}-2\crrterm{\frac{-\w+41}{44}}-2\crrterm{\frac{-\w+8}{48}}-2\crrterm{\frac{-15\w+147}{104}}-2\crrterm{\frac{-4\w+21}{13}}-2\crrterm{\frac{-4\w+21}{33}}-2\crrterm{\frac{-21\w+172}{352}}-2\crrterm{\frac{-21\w+201}{352}}-2\crrterm{\frac{-23\w+211}{256}}-2\crrterm{\frac{-23\w+124}{312}}-2\crrterm{\frac{-2457\w-611}{22528}}-2\crrterm{\frac{-253\w+2321}{6144}}-2\crrterm{\frac{-253\w-495}{3328}}-2\crrterm{\frac{-27\w+535}{676}}-2\crrterm{\frac{-27\w+168}{676}}-2\crrterm{\frac{-27\w+535}{832}}-2\crrterm{\frac{-29\w+1332}{1331}}-2\crrterm{\frac{-29\w+149}{484}}-2\crrterm{\frac{-3\w+84}{64}}-2\crrterm{\frac{-3\w+45}{88}}-2\crrterm{\frac{-3\w+46}{88}}-2\crrterm{\frac{-3\w-9}{11}}-2\crrterm{\frac{-3\w+36}{16}}-2\crrterm{\frac{-3\w+45}{22}}-2\crrterm{\frac{-3\w+46}{22}}-2\crrterm{\frac{-3\w-20}{22}}-2\crrterm{\frac{-3\w-21}{22}}-2\crrterm{\frac{-3\w-17}{256}}-2\crrterm{\frac{-3\w+15}{32}}-2\crrterm{\frac{-3\w+24}{44}}-2\crrterm{\frac{-39\w-221}{512}}-2\crrterm{\frac{-8\w+31}{39}}-2\crrterm{\frac{-51\w+3807}{3328}}-2\crrterm{\frac{-12\w+63}{143}}-2\crrterm{\frac{-9\w+17}{26}}+2\crrterm{\frac{\w+11}{8}}+2\crrterm{\frac{\w+3}{8}}+2\crrterm{\frac{\w-4}{8}}+2\crrterm{\frac{\w+25}{1}}+2\crrterm{\frac{\w+14}{11}}+2\crrterm{\frac{\w+7}{12}}+2\crrterm{\frac{\w+11}{16}}+2\crrterm{\frac{\w+14}{16}}+2\crrterm{\frac{\w+1}{16}}+2\crrterm{\frac{\w+4}{16}}+2\crrterm{\frac{\w-25}{16}}+2\crrterm{\frac{\w+14}{22}}+2\crrterm{\frac{\w+25}{24}}+2\crrterm{\frac{\w-2}{24}}+2\crrterm{\frac{\w+14}{26}}+2\crrterm{\frac{\w-2}{3}}+2\crrterm{\frac{\w+4}{32}}+2\crrterm{\frac{\w+7}{4}}+2\crrterm{\frac{\w-4}{4}}+2\crrterm{\frac{\w+11}{48}}+2\crrterm{\frac{\w+7}{48}}+2\crrterm{\frac{15\w+132}{104}}+2\crrterm{\frac{15\w+204}{176}}+2\crrterm{\frac{21\w-3}{169}}+2\crrterm{\frac{21\w+151}{352}}+2\crrterm{\frac{21\w+180}{352}}+2\crrterm{\frac{23\w+45}{256}}+2\crrterm{\frac{27\w-51}{484}}+2\crrterm{\frac{29\w+335}{512}}+2\crrterm{\frac{29\w-1}{176}}+2\crrterm{\frac{3\w-123}{121}}+2\crrterm{\frac{3\w+33}{13}}+2\crrterm{\frac{3\w+33}{16}}+2\crrterm{\frac{3\w+43}{22}}+2\crrterm{\frac{3\w+3}{26}}+2\crrterm{\frac{3\w+12}{32}}+2\crrterm{\frac{3\w+17}{32}}+2\crrterm{\frac{3\w+20}{32}}+2\crrterm{\frac{3\w+20}{44}}+2\crrterm{\frac{33\w+2651}{5408}}+2\crrterm{\frac{8\w+23}{39}}+2\crrterm{\frac{6399\w+16348}{42592}}+2\crrterm{\frac{7\w-67}{64}}+2\crrterm{\frac{7\w-1}{66}}+2\crrterm{\frac{7\w+181}{312}}+2\crrterm{\frac{9\w+60}{52}}+2\crrterm{\frac{9\w+8}{26}}+2\crrterm{\frac{9\w+8}{44}}.
\end{dmath*}
}

{\tiny
\begin{dmath*}
\bgeo =108\crrterm{\frac{\w+4}{6}}+36\crrterm{\frac{\w+3}{5}}+108\crrterm{\frac{\w+11}{13}}+36\crrterm{\frac{\w+1}{5}}+64\crrterm{\frac{\w+3}{4}}+24\crrterm{\frac{\w+14}{26}}+12\crrterm{\frac{3\w-3}{35}}
+36\crrterm{\frac{3\w}{38}}+64\crrterm{\frac{\w}{4}}+180\crrterm{\frac{\w+3}{8}}+12\crrterm{\frac{\w-4}{22}}+24\crrterm{\frac{\w+36}{32}}+88\crrterm{\frac{\w+5}{10}}+36\crrterm{\frac{21\w+192}{689}}
+24\crrterm{\frac{\w-4}{8}}+136\crrterm{\frac{3\w+17}{32}}+180\crrterm{\frac{\w+3}{11}}+12\crrterm{\frac{5\w+44}{104}}+30\crrterm{\frac{9\w+45}{106}}+12\crrterm{\frac{\w-6}{2}}+88\crrterm{\frac{5\w+23}{48}}
+20\crrterm{\frac{5\w-25}{3}}+12\crrterm{\frac{9\w+38}{38}}+12\crrterm{\frac{\w+19}{20}}+48\crrterm{\frac{\w+4}{24}}+24\crrterm{\frac{5\w-25}{128}}+24\crrterm{\frac{5\w+148}{128}}
+60\crrterm{\frac{\w+24}{26}}+24\crrterm{\frac{42\w+665}{1121}}+24\crrterm{\frac{6\w+95}{77}}+12\crrterm{\frac{2\w+17}{33}}+12\crrterm{\frac{\w-1}{3}}+48\crrterm{\frac{\w-1}{19}}+24\crrterm{\frac{7\w-28}{55}}
+48\crrterm{\frac{35\w+644}{984}}+24\crrterm{\frac{4\w+24}{59}}+12\crrterm{\frac{4\w-28}{7}}+24\crrterm{\frac{7\w+76}{55}}+48\crrterm{\frac{7\w-28}{40}}+60\crrterm{\frac{\w+3}{7}}+12\crrterm{\frac{\w-5}{7}}
+24\crrterm{\frac{5\w-25}{56}}+24\crrterm{\frac{\w+11}{9}}+24\crrterm{\frac{9\w+99}{208}}+60\crrterm{\frac{4\w+8}{41}}+24\crrterm{\frac{15\w-4}{26}}+36\crrterm{\frac{\w+27}{26}}+36\crrterm{\frac{\w+27}{32}}
+12\crrterm{\frac{2\w+10}{53}}+12\crrterm{\frac{2\w+41}{45}}+12\crrterm{\frac{\w+75}{76}}+6\crrterm{\frac{\w+5}{5}}+6\crrterm{\frac{5\w}{81}}+12\crrterm{\frac{25\w+1123}{1298}}+12\crrterm{\frac{11\w+66}{118}}+60\crrterm{\frac{4\w+29}{41}}
+24\crrterm{\frac{15\w+15}{26}}+28\crrterm{\frac{5\w-28}{1272}}+12\crrterm{\frac{25\w+1175}{1166}}+12\crrterm{\frac{\w+5}{3}}+72\crrterm{\frac{3\w+6}{41}}+24\crrterm{\frac{35\w+440}{1007}}+24\crrterm{\frac{3\w+20}{104}}
+24\crrterm{\frac{\w+58}{54}}+24\crrterm{\frac{9\w+522}{583}}+24\crrterm{\frac{3\w+20}{11}}+12\crrterm{\frac{\w+37}{39}}+12\crrterm{\frac{6\w+19}{247}}+12\crrterm{\frac{15\w+980}{902}}+12\crrterm{\frac{\w+7}{10}}
+24\crrterm{\frac{3\w-6}{13}}+24\crrterm{\frac{3\w+57}{76}}+24\crrterm{\frac{\w+36}{44}}+12\crrterm{\frac{\w+76}{78}}+12\crrterm{\frac{27\w+27}{130}}+16\crrterm{\frac{3\w-21}{20}}+28\crrterm{\frac{3\w+18}{59}}
+48\crrterm{\frac{5\w+67}{82}}+28\crrterm{\frac{15\w+1370}{1298}}+28\crrterm{\frac{\w+6}{10}}+72\crrterm{\frac{3\w+32}{41}}+12\crrterm{\frac{8\w+57}{209}}+12\crrterm{\frac{\w+2}{3}}+24\crrterm{\frac{\w+2}{7}}
+12\crrterm{\frac{7\w+64}{78}}+12\crrterm{\frac{16\w+133}{53}}+36\crrterm{\frac{3\w+15}{53}}+24\crrterm{\frac{\w+4}{7}}+48\crrterm{\frac{2\w+2}{39}}+12\crrterm{\frac{\w-1}{2}}+12\crrterm{\frac{18\w+684}{779}}
+12\crrterm{\frac{9\w-72}{308}}+12\crrterm{\frac{\w+7}{7}}+48\crrterm{\frac{\w+1}{4}}+24\crrterm{\frac{175\w+1600}{4134}}+30\crrterm{\frac{9\w+52}{106}}+6\crrterm{\frac{25\w+867}{792}}+24\crrterm{\frac{7\w+64}{50}}
24\crrterm{\frac{36\w+1876}{2173}}+24\crrterm{\frac{\w+40}{44}}+24\crrterm{\frac{3\w+20}{14}}+24\crrterm{\frac{\w+3}{44}}+24\crrterm{\frac{7\w+137}{123}}+6\crrterm{\frac{27\w-28}{26}}+6\crrterm{\frac{\w+79}{82}}
+12\crrterm{\frac{9\w+55}{118}}+12\crrterm{\frac{16\w+155}{779}}+28\crrterm{\frac{\w+6}{6}}+4\crrterm{\frac{3\w+17}{5}}+16\crrterm{\frac{30\w+627}{1007}}+24\crrterm{\frac{4\w-24}{9}}+12\crrterm{\frac{84\w+1480}{2173}}
+24\crrterm{\frac{9\w-36}{28}}+12\crrterm{\frac{27\w+988}{826}}+12\crrterm{\frac{7\w+69}{76}}+24\crrterm{\frac{7\w-28}{198}}+12\crrterm{\frac{21\w+151}{130}}+12\crrterm{\frac{\w+84}{88}}+6\crrterm{\frac{81\w+331}{574}}
+12\crrterm{\frac{9\w+113}{104}}+24\crrterm{\frac{4\w+31}{59}}+16\crrterm{\frac{11\w+40}{106}}+8\crrterm{\frac{165\w-765}{11236}}+8\crrterm{\frac{11\w-51}{15}}+36\crrterm{\frac{21\w+476}{689}}+28\crrterm{\frac{3\w+274}{295}}
+28\crrterm{\frac{5\w+30}{59}}+12\crrterm{\frac{9\w-54}{7}}+8\crrterm{\frac{5\w-20}{9}}+8\crrterm{\frac{9\w-36}{440}}+4\crrterm{\frac{\w-4}{15}}+4\crrterm{\frac{5\w+95}{152}}+4\crrterm{\frac{3\w+35}{50}}+4\crrterm{\frac{\w+19}{25}}
+16\crrterm{\frac{10\w+209}{159}}+20\crrterm{\frac{\w-5}{160}}+28\crrterm{\frac{\w+5}{265}}-12\crrterm{\frac{-3\w}{35}}-12\crrterm{\frac{-\w-3}{22}}-12\crrterm{\frac{-5\w+49}{104}}-12\crrterm{\frac{-\w-5}{2}}
-12\crrterm{\frac{-9\w+47}{38}}-12\crrterm{\frac{-\w+20}{20}}-12\crrterm{\frac{-2\w+19}{33}}-12\crrterm{\frac{-\w}{3}}-12\crrterm{\frac{-4\w-24}{7}}-12\crrterm{\frac{-\w-4}{7}}-12\crrterm{\frac{-2\w+12}{53}}
-12\crrterm{\frac{-2\w+43}{45}}-6\crrterm{\frac{-\w+6}{5}}-6\crrterm{\frac{-5\w+5}{81}}-12\crrterm{\frac{-25\w+1148}{1298}}-12\crrterm{\frac{-11\w+77}{118}}-12\crrterm{\frac{-25\w+1200}{1166}}
-12\crrterm{\frac{-\w+6}{3}}-12\crrterm{\frac{-\w+38}{39}}-12\crrterm{\frac{-6\w+25}{247}}-12\crrterm{\frac{-15\w+995}{902}}-12\crrterm{\frac{-\w+8}{10}}-12\crrterm{\frac{-\w+77}{78}}-12\crrterm{\frac{-27\w+54}{130}}
-16\crrterm{\frac{-3\w-18}{20}}-12\crrterm{\frac{-8\w+65}{209}}-12\crrterm{\frac{-\w+3}{3}}-12\crrterm{\frac{-7\w+71}{78}}-12\crrterm{\frac{-16\w+149}{53}}-12\crrterm{\frac{-\w}{2}}-12\crrterm{\frac{-18\w+702}{779}}
-12\crrterm{\frac{-9\w-63}{308}}-12\crrterm{\frac{-\w+8}{7}}-6\crrterm{\frac{-25\w+892}{792}}-6\crrterm{\frac{-27\w-1}{26}}-6\crrterm{\frac{-\w+80}{82}}-12\crrterm{\frac{-9\w+64}{118}}-12\crrterm{\frac{-16\w+171}{779}}
-4\crrterm{\frac{-3\w+20}{5}}-12\crrterm{\frac{-84\w+1564}{2173}}-12\crrterm{\frac{-27\w+1015}{826}}-12\crrterm{\frac{-7\w+76}{76}}-12\crrterm{\frac{-21\w+172}{130}}-12\crrterm{\frac{-\w+85}{88}}
-6\crrterm{\frac{-81\w+412}{574}}-12\crrterm{\frac{-9\w+122}{104}}-12\crrterm{\frac{-9\w-45}{7}}-4\crrterm{\frac{-\w-3}{15}}
-4\crrterm{\frac{-5\w+100}{152}}-4\crrterm{\frac{-3\w+38}{50}}-4\crrterm{\frac{-\w+20}{25}};
\end{dmath*}
}

\bibliographystyle{amsplain_initials_eprint_doi_url}
\bibliography{k3gens}
\catcode`\_=8  
\end{document}